\documentclass[preprint,10pt]{elsarticle}

\usepackage{hyperref}
\usepackage{url}
\usepackage{graphicx}
\usepackage{subfig}
\usepackage{epsfig}
\usepackage{psfrag}
\usepackage{amsmath,amssymb,amsfonts,amsthm,stmaryrd}
\usepackage{dsfont}
\usepackage{eucal}
\usepackage{mathrsfs}
\usepackage{listings}    
\usepackage{color}
\usepackage{float}
\usepackage{natbib}
\usepackage{booktabs}
\usepackage{tabularx}
\usepackage{setspace}
\usepackage{algorithmicx}
\usepackage{algorithm}
\usepackage{algpseudocode}

\makeatletter
\renewcommand\fs@ruled{\def\@fs@cfont{\bfseries}\let\@fs@capt\floatc@ruled
\def\@fs@pre{\hrule height 1.2pt depth0pt \kern2pt}%
\def\@fs@post{\kern2pt\hrule height 1.2pt depth0pt \kern2pt \relax}%
\def\@fs@mid{\kern2pt\hrule\kern2pt}%
\let\@fs@iftopcapt\iftrue}
\makeatother

\newcommand{\norm}[1]{\left|\left|#1\right|\right|}

\newcommand{\cref}[1]{Chapter~\ref{chapt:#1}}

\newcommand{\D}{\displaystyle}

\newcommand{\di}{\mathrm{d}}

\usepackage{xspace}
\newcommand{\eg}{e.g.,\xspace}

\newcommand{\ie}{i.e.,\xspace}

\newcommand{\etc}{etc.\@\xspace}

\newcommand{\cmatrixb}{\left\{ \begin{matrix}}
\newcommand{\cmatrixe}{\end{matrix} \right\}}

\newcommand{\vm}[1]{\mathbf{#1}}

\newcommand{\bsym}[1]{\boldsymbol{#1}}

\newcommand{\trans}{^\mathrm{T}}

\newcommand{\bc}{\begin{center}}
\newcommand{\ec}{\end{center}}
\newcommand{\bitem}{\begin{itemize}}
\newcommand{\eitem}{\end{itemize}}

\newcommand{\ljump}{\lbrack \! \lbrack } 
\newcommand{\rjump}{\rbrack \! \rbrack } 
\newcommand{\jump}[1]{\ljump {#1} \rjump} 

\newcommand{\pderiv}[2]{ \frac{\partial {#1} }{\partial {#2} } }

\newcommand{\beq}{\begin{equation}}
\newcommand{\eeq}{\end{equation}}
\newcommand{\beqa}{\begin{eqnarray}}
\newcommand{\eeqa}{\end{eqnarray}}

\newcommand{\invisible}[1]{}

\newcommand{\bv}{\begin{verbatim}}
\newcommand{\V}{\verb}                  % EX: \V=-d{#@~}= Expr must

\newcommand{\testpix}[1]{\fbox{\begin{minipage}[c]{\textwidth}
                      #1 \end{minipage} }}

\newcommand{\putpstex}[1]{\includegraphics{#1.pstex_t}}

\newcommand{\grbf}[1]{\mbox{\boldmath{$#1$}}}

\floatstyle{ruled}
\newfloat{Program}{thp}{lop}
\floatname{Program}{Program}

\floatstyle{ruled}
\newfloat{Fbox}{thp}{lop}
\floatname{Fbox}{Box}

\newcolumntype{C}{>{\centering\arraybackslash}X}

\definecolor{darkgray}{rgb}{0.95,0.95,0.95}
\definecolor{mygreen}{rgb}{0,0.6,0}

\lstdefinestyle{Matlab}
{
 basicstyle=\footnotesize, numbers=none, numberstyle=\tiny,%
 showstringspaces=false, language=Matlab, escapechar=|,frame=tb,%
 commentstyle=\color{mygreen}
}

\lstdefinestyle{Matlab-num}
{
 basicstyle=\footnotesize, numbers=left, numberstyle=\tiny,%
 showstringspaces=false, language=Matlab, escapechar=|,frame=tb,%
commentstyle=\color{mygreen}
}

\newcommand{\tty}[1]{\textnormal{\texttt{#1}}}
\newcommand{\sym}[1]{\textnormal{\textit{#1}}}

\lstnewenvironment{snippet}[1][]
{
 \lstset{style=Matlab, xleftmargin=5mm, gobble=4, #1}
}
{}

\lstnewenvironment{snippet1}[1][]
{
 \lstset{style=Matlab-num, xleftmargin=5mm, gobble=4, #1}
}
{}

\makeatletter
\renewcommand*\env@matrix[1][\arraystretch]{%
  \edef\arraystretch{#1}%
  \hskip -\arraycolsep
  \let\@ifnextchar\new@ifnextchar
  \array{*\c@MaxMatrixCols c}}
\makeatother

\def\bibsection{\section*{References}}

\begin{document}


\definecolor{MyDarkBlue}{rgb}{1, 0.9, 1}
\lstset{language=Matlab,
       basicstyle=\footnotesize,
       commentstyle=\itshape,
       stringstyle=\ttfamily,
       showstringspaces=false,
       tabsize=2}
\lstdefinestyle{commentstyle}{color=\color{green}}

\theoremstyle{remark}
\newtheorem{thm}{Theorem}[section]
\newtheorem{rmk}[thm]{Remark}


\definecolor{red}{gray}{0}
\definecolor{blue}{gray}{0}


\begin{frontmatter}

\title{Nitsche's method method for mixed dimensional analysis: conforming and non-conforming continuum-beam and continuum-plate coupling}       

\author[cardiff]{Vinh Phu Nguyen  \fnref{fn1}}
\author[cardiff]{Pierre Kerfriden  \fnref{fn2}}
\author[ucl]{Susanne Claus  \fnref{fn4}}
\author[cardiff]{St\'{e}phane P.A. Bordas \corref{cor1}\fnref{fn3}}

\cortext[cor1]{Corresponding author}

\address[cardiff]{School of Engineering, Institute of Mechanics and Advanced
Materials, Cardiff University, Queen's Buildings, The Parade, Cardiff \\
CF24 3AA}

\address[ucl]{Department of Mathematics, University College London, London, WC1E 6BT, United Kingdom}

\fntext[fn1]{\url nguyenpv@cardiff.ac.uk, ORCID: 0000-0003-1212-8311}
\fntext[fn2]{\url pierre@cardiff.ac.uk}
\fntext[fn4]{\url susanne.claus@ucl.ac.uk}
\fntext[fn3]{\url stephane.bordas@alum.northwestern.edu, ORCID: 0000-0001-7622-2193}


\begin{abstract}
A Nitche's method is presented to couple different mechanical models. They 
include coupling of a solid and a beam and of a solid
and a plate. Both conforming and non-conforming formulations are presented. In a non-conforming formulation,
the structure domain is overlapped by a refined solid model.    
Applications can be found in multi-dimensional analyses in which parts of a structure are 
modeled with solid elements
and others are modeled using a coarser model with beam and/or plate elements. 
Discretisations are performed using both standard Lagrange elements and high order NURBS
(Non Uniform Rational Bsplines) based isogeometric elements.
We present various examples to demonstrate the performance of the method.
\end{abstract}


\end{frontmatter}


\section{Introduction}

Nitsche's method \cite{nitsche} was originally proposed to weakly enforce Dirichlet boundary conditions
as an alternative to equivalent pointwise constraints.
The idea behind a Nitsche based approach is 
to replace the Lagrange multipliers arising in a dual formulation through their physical representation, 
namely the normal flux at the interface. Nitsche also added an extra penalty like term to restore 
the coercivity of the bilinear form. The method can be seen to lie in between the Lagrange multiplier
method and the penalty method.
The method has seen a resurgence in recent years and was applied for interface problems 
\cite{Hansbo20025537,Dolbow2009a},  for connecting overlapping meshes 
\cite{MZA:8203296,MZA:8203286,Sanders2012a,Sanders2011a}, for imposing Dirichlet boundary conditions 
in meshfree methods \cite{FernándezMéndez20041257}, in immersed boundary methods 
\cite{NME:NME4522,NME:NME3339,embar_imposing_2010}, in fluid mechanics 
\cite{Burman-Nitsche-2009,Bazilevs200712}
and for contact mechanics \cite{nitsche-wriggers}. It has also been applied for stabilising constraints
in enriched finite elements \cite{Sanders2008a,Burman-Nitsche-2012}.

In this paper, a Nitsche's method is presented to (1) couple two dimensional (2D) continua and beams
and (2) three dimensional (3D) continua and plates. The continua and the structures are discretised using
either Lagrange finite elements (FEs) or high order B-spline/NURBS isogeometric finite elements.
The need for the research presented in this paper stems from the problem of progressive
failure analysis of composite laminates of which recent studies were carried out by the authors
\cite{Nguyen2013,nguyen_cohesive_2013,nguyen-offset}. It was shown that using NURBS (Non Uniform
Rational B-splines) as the finite element basis functions--the concept coined as
isogeometric analysis (IGA) by Hughes and his co-workers \cite{hughes_isogeometric_2005,cottrel_book_2009}
results in speed up in both pre-processing and processing step of delamination analysis of composite laminates.
However in order for IGA to be applied for real problems such as large composite panels utilized in aerospace
industry Mixed Dimensional Analysis (MDA) \cite{Cuilliere2010} must be employed. In MDA 
some portions of the object of interest can be modeled using reduced-dimensional elements 
(typically beams and shells), while other portions must be modelled using volume (solid) elements due to the need
of accuracy.
In this way, large structures are likely tractable. 
The question is how to model the coupling of different models in a flexible and efficient manner.

Broadly speaking the coupling can be either surface coupling (non-overlapping coupling) or
volume coupling (overlapping coupling) \cite{Guidault2007b}. 
Volume coupling indicates the existence of 
a region in which both model co-exist and is usually realized using the Arlequin method 
\cite{Dhia2005a} of which Abaqus implementation can be found in \cite{Qiao2011a}.
Arlequin method is best suited for coupling different physical models such as continuum-particle
modeling see Refs. \cite{Bauman2008a,Pfaller2011a,Wellmann2012a,Rousseau2009a} among others.
In Arlequin method care must be taken when choosing the space of Lagrange multipliers
and the numerical integration of coupling terms on unstructured meshes is a non-trivial task, particularly
for three dimensional problems. In surface coupling there is
no overlapping of models and the two models can be coupled using one of the following methods

\begin{enumerate}
\item Lagrange multiplier method  as in the mortar method \cite{mortar-1999};
\item Penalty method;
\item Multipoint constraint method;
\item Transition element method
\end{enumerate}

Both Lagrange multiplier method and the penalty method
have their own disadvantages--in the former it is the introduction of extra unknowns 
(thus destroys the positive definiteness of the augmented system of equation) and the proper choice
of the Lagrange multiplier discrete space. In the latter it is the choice of the penalty parameter.
Other method is to use constraint equations \cite{Unger,NME:NME967,Song-MDA} 
(the common multipoint constraint MPC method
in commercial FE packages such as Abaqus, Ansys and MSC Nastran). 
MPC based coupling is a strong coupling method whereas all other coupling methods
are weak couplings.
In \cite{NME:NME967,Shim-MDA} constraint equations are determined by equating the work done by the stresses in each part of the model at the interface between dimensions. 
Example results show that the proposed technique does not introduce any spurious 
stresses at the dimensional interface.
The authors also provided error estimations for their
MDA. The theory of MPC method is described in details in \cite{felippa:note}. In \cite{Unger} a comparative study of three coupling methods: MPC, mortar method and Arlequin method
was presented for the concurrent multiscale modeling of concrete material where the macro domain is
discretised by a coarse mesh and the micro domain (with heterogeneities) is discretised by a very fine mesh.
The authors concluded that the Arlequin method is too expensive to consider.
A somehow related approach is the global-local technique, see \eg \cite{felippa:global-local}
and references therein. 
In a global-local approach, the whole system is first analyzed as a global entity, discarding details 
deemed not to affect its overall behavior. Local details (cutoffs, cracks \etc) are then analyzed 
using the results of the global analysis as boundary conditions. In \cite{allix-non-intrusive} a global-local
technique was adopted as a non-intrusive coupling method.
The authors in \cite{Gabbert} compared the performance, through a simple numerical example, of
the MPC method, the Arlequin method and the global-local method and the conclusion was that the Arlequin method
is a promising coupling method for it results in lower stress jumps at the coupling boundary.
In aforementioned coupling methods either a modified weak formulation has to be used
or a system of linear equations has to be solved with a set of (linear) constraints (for the MPC method).
Transition elements are a yet another method in which standard weak forms can be reused 
to couple solid elements and beam/shell elements 
see \eg \cite{NME:NME1620150704,NME:NME1620360902,transition-gmur,NME:NME938,Garusi2002105} or in an IGA context the bending strip method \cite{kiendl_bending_2010,Raknes2013127}. 
Using transition elements, a mesh-geometry-based solution to mixed-dimensional coupling was presented in 
\cite{Cuilliere2010}.
It should be emphasized that solving mixed dimensional models is much faster than solving full 3D models, but the time required to process these mixed dimensional models might eliminate that advantage.

\begin{figure}
  \centering
  \includegraphics[width=0.6\textwidth]{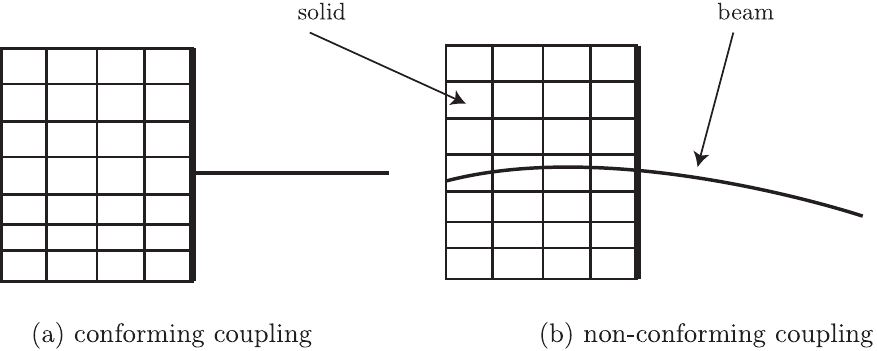}
  \caption{Solid-structure coupling: (a) conforming coupling and (b) non-conforming coupling. 
     The bold lines represent the coupling line/surface.}
  \label{fig:concepts}
\end{figure}

In this paper, we decided to use a surface coupling method and a Nitsche based weak coupling formulation.
Surface coupling was chosen because  we are not solving a multi-physics model and the Nitsche's method
is adopted due to a rich mathematics behind it and the related discontinuous Galerkin methods \cite{Arnold01unifiedanalysis}. The method is symmetric for symmetric problems and do not need additional degrees of freedom. 
However, a user-defined stabilization parameter is required.
Our target is modeling large scale composite laminates and the
long-term goal is to start the analysis from a laminate shell model, build goal-oriented error estimates 
to identify the \textit{hot-spots} where intra or inter-laminar failure is likely to take place, to refine 
these zones with a continuum description, and, once the delaminations and intra-laminar cracks are fully open, replace them by an homogenised shell model. 

The contribution of the paper are listed as follows

\begin{enumerate}
\item Formulation and implementation details for 2D solid and beam coupling;
\item Formulation and implementation details for 3D solid and plate coupling;
\item Conforming and non-conforming coupling are presented, cf. Fig.~\ref{fig:concepts};
\end{enumerate}
The remainder of the paper is organised as follows. Section \ref{problem} presents the problem
description followed by variational formulations given in Section \ref{weak-form}. 
Discretisation is discussed in Section \ref{discretisation} together with implementation details. 
Non-conforming coupling is detailed in
Section \ref{non-conforming} and
numerical examples are provided in Section \ref{sec:examples}.
The presented algorithm applies for both Lagrange basis functions and NURBS basis functions. The latter
have been extensively used in isogeometric analysis (IGA) \cite{hughes_isogeometric_2005,cottrel_book_2009}--a
methodology aims at reducing the gap between FEA (FE Analysis) and CAD (Computer Aided Design) and
facilitates the implementation of rotation-free beam/plate/shell elements 
see \eg \cite{kiendl_isogeometric_2009,benson_large_2011} which is also adopted in this work. 
Additionally, IGA is also the numerical framework
for our recent works on failure analysis of composite laminates 
\cite{Nguyen2013,nguyen_cohesive_2013,nguyen-offset}. 
Coupling formulation for both classical beam/plate model and first order shear deformation 
beam/plate model are presented.
We confine to quasi-static small strain problems 
although Nitsche's method was successfully applied to finite deformation \cite{Sanders2012a} and
free vibration analysis is presented in \cite{nguyen-nitsche1}.
The material is assumed to be homogeneous isotropic but can be applied equally well for composite materials.
We refer to \cite{nguyen-nitsche1} for a similar work in which non-conforming NURBS patches are
weakly glued together using a Nitsche's method. It should be emphasised that 
the bending strip method \cite{kiendl_bending_2010,Raknes2013127} which is used in IGA to join
$C^0$ shell patches can be used to couple a continuum and a shell under a restriction that the
parametrisation is compatible at the interface.

We denote $d_p$ and $d_s$ as the number of parametric directions and spatial directions respectively.
Both tensor and matrix notations are used.
In tensor notation, tensors of order one or greater are written in boldface. Lower case bold-face letters 
are used for first-order tensor whereas upper case bold-face letters indicate high-order tensors. 
The major exception to this rule are the physical second order stress tensor and the strain tensor
which are written in lower case. In matrix notation, the same symbols as for tensors are used to denote
the matrices but the connective operator symbols are skipped and 
second order tensors ($\sigma_{ij}$ and $\epsilon_{ij}$) are written using the Voigt notation
as column vectors; 
$\bsym{\sigma}=[\sigma_{xx}, \sigma_{yy}, \sigma_{zz}, \sigma_{xy}, \sigma_{yz}, \sigma_{xz}]\trans$, 
$\bsym{\epsilon}=[\epsilon_{xx}, \epsilon_{yy}, \epsilon_{zz}, 2\epsilon_{xy}, 2\epsilon_{yz}, 2\epsilon_{xz}]\trans$.

\section{Problem description}\label{problem}

In this section the governing equations of two mechanical systems are presented namely (1) solid-beam
and (2) solid-plate system. For sake of presentation, we use the classical 
Euler-Bernoulli beam theory and Kirchhoff plate
theory. The extension to Timoshenko beam and Mindlin-Reissner plate is provided in Section \ref{shear}.
For subsequent development, superscripts $s$, $b$ and $p$ will be adopted to denote quantities associated with the solid, beam and plate, respectively. 

\subsection{Solid-beam coupling}

We define the domain $\Omega \subset \mathbb{R}^{d_s}$ which is divided into two non-overlapping domains--
$\Omega^s$ for the solid part and $\Omega^b$ for the beam part, cf. Fig.~\ref{fig:domain}. 
The boundary of $\Omega^s$ is partitioned as $\Gamma^s = \overline{\Gamma_u^s \cup \Gamma_t^s}$, 
$\Gamma_u^s \cap \Gamma_t^s = \emptyset$ where $\Gamma_u^s$ and $\Gamma_t^s$ denote the Dirichlet and Neumann boundaries respectively with an overline representing a closed set. Let $\vm{u}^s$, $\bsym{\epsilon}^s$ and $\bsym{\sigma}^s$ be the displacements, strains and stresses in the solid part, respectively. In the beam part, $w$ denotes the transverse displacement. The global coordinate system is denoted by $(x,y)$ and a local coordinate system
$(\bar{x},\bar{y})$ is adopted for the beam. Note that $\Omega_b$ is the mid-line of the beam $\Omega^b_{ext}$.
The coupling interface, also referred to as dimensional interface in the literature, $\Gamma^*$ is defined as $\Gamma^*=\Omega^s\cap\Omega^b_{ext}$.

The governing equations are 

\begin{itemize}

\item For the solid part

\begin{subequations}
\begin{alignat}{2}
-\nabla \;\bsym{\sigma}^s      &=  \vm{b}        &\quad\text{on} \quad \Omega^s \\ 
             \vm{u}^s          &=  \bar{\vm{u}}  &\quad\text{on} \quad \Gamma_u^s \\
  \bsym{\sigma}^s \cdot \vm{n} &=  \bar{\vm{t}}  &\quad\text{on} \quad \Gamma_t^s \label{eq:Neumann}
\end{alignat}
\label{eq:solid}
\end{subequations}
\noindent where the strain is taken as the symmetric part of the displacement gradient $\bsym{\epsilon}^s=\frac{1}{2}(\nabla\vm{u}^s + \nabla\trans \vm{u}^s)$. And the Cauchy stress is a linear function of the strains $\bsym{\sigma}^s=
\vm{C}^s:\bsym{\epsilon}^s$ where $\vm{C}^s$ denotes the fourth order tensor of material properties of the solid 
according to Hooke's law. Prescribed displacements and tractions are denoted by $\bar{\vm{u}}$ and $\bar{\vm{t}}$, respectively.

\item For the beam part

\begin{subequations}
\begin{alignat}{2}
EI \frac{d^4 w}{d\bar{x}^4} &=  p         &\quad\text{on} \quad \Omega^b \\ 
             w         &=  \bar{w}   &\quad\text{on} \quad \Gamma_u^b \\
             \frac{dw}{d\bar{x}}         &=  -\bar{\theta}   &\quad\text{on} \quad \Gamma_\theta^b \\
             EI\frac{d^2w}{d\bar{x}^2}n         &=  \bar{m}   &\quad\text{on} \quad \Gamma_m^b 
\end{alignat}
\end{subequations}
where $I$ is  the moment of inertia of the beam which is, for a rectangular beam, 
given by $I=\frac{bh^3}{12}$ with $h$ being the
beam thickness and $b$ denotes the beam width; $E$ is the Young's modulus and
$p$ denotes the distributed pressure load. $\bar{w}$, $\bar{\theta}$ and $\bar{m}$ are the prescribed deflection, rotation and moment, respectively. $n$ denotes the normal to the natural boundary conditions. 

\item For the coupling part

\begin{subequations}
\begin{alignat}{2}
  \vm{u}^s &=  \vm{u}^b                        &\quad\text{on} \quad \Gamma^* \\
  \bsym{\sigma}^s \cdot \vm{n}^s &= \bsym{\sigma}^b \cdot \vm{n}^s  &\quad\text{on} \quad \Gamma^* \label{eq:tr}
\end{alignat}
\label{eq:coupling-beam}
\end{subequations}
where $\vm{n}^s$ is the outward unit normal vector to the coupling interface $\Gamma^*$;
$\vm{u}^b$ denote the beam displacement field, in the global coordinate system, of a point 
on the coupling interface $\Gamma^*$ which is defined as
\begin{equation}
\vm{u}^b = \vm{R}_v\trans \bar{\vm{u}}^b,\quad
\bar{\vm{u}}^b = \begin{bmatrix} -\bar{y}w_{,\bar{x}} & w(\bar{x}) \end{bmatrix}\trans
\label{eq:beam-disp}
\end{equation}
quantities with a bar overhead denote local quantities defined in the beam coordinate system
and $\bsym{\sigma}^b$ is the beam stress field defined in the global system which is given by
\begin{equation}
\bsym{\sigma}^b = \vm{T}^{-1} \bar{\bsym{\sigma}}^b
\end{equation}
where $\vm{R}_v$ and $\vm{T}^{-1}$ denote the rotation matrices for vector transformation 
and stress transformation, respectively 

\begin{equation}
\vm{R}_v = \begin{bmatrix}
\cos\phi & \sin\phi\\
-\sin\phi & \cos\phi
\end{bmatrix}, \quad
\vm{T}^{-1} = \begin{bmatrix}
\cos^2\phi & \sin^2\phi & -2\sin\phi\cos\phi\\
\sin^2\phi & \cos^2\phi & 2\sin\phi\cos\phi\\
\sin\phi\cos\phi & -\sin\phi\cos\phi & \cos^2\phi-\sin^2\phi
\end{bmatrix}
\label{eq:Rv}
\end{equation}
The beam stresses are defined as
\begin{equation}
\bar{\bsym{\sigma}}^b\equiv
\begin{bmatrix}
\bar{\sigma}_{xx}^b\\
\bar{\sigma}_{yy}^b\\
\bar{\sigma}_{xy}^b
\end{bmatrix}=\underbrace{
\begin{bmatrix}
E &  0 & 0\\
0 &  0 & 0\\
0 & 0 & 0  
\end{bmatrix}}_{\vm{C}^b}
\begin{bmatrix}
\bar{\epsilon}_{xx}^b\\
\bar{\epsilon}_{yy}^b\\
2\bar{\epsilon}_{xy}^b
\end{bmatrix}
\label{eq:beam-stress-strain}
\end{equation}
where the only non-zero strain component is the bending or flexural strain 
$\bar{\epsilon}_{xx}=-\bar{y}w_{,\bar{x}\bar{x}}$. Equations~\eqref{eq:beam-disp} and
\eqref{eq:beam-stress-strain} that transform the one dimensional beam variables (the transverse displacement
$w$) into 2D fields ($\bar{\vm{u}}^b$ and $\bar{\bsym{\sigma}}^b$) are referred to as \textit{prolongation
operators}. Note that these operators are defined only along $\Gamma^*$ where the coupling is taking place.
\end{itemize}

\begin{figure}
  \centering
  \includegraphics[width=\textwidth]{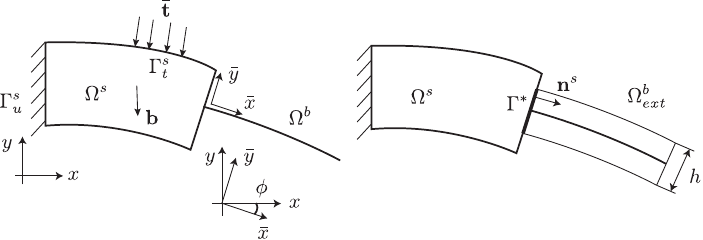}
  \caption{Coupling of a two dimensional solid and a beam.}
  \label{fig:domain}
\end{figure}

\subsection{Solid-plate coupling}

The coupling of a solid and a plate is graphically illustrated by Fig.~\ref{fig:solid-plate1}.
The coupling interface $\Gamma^*$ is the intersection of the solid domain and the three dimensional
plate. The mid-surface of the plate is denoted by $\Omega^p$. 
The global coordinate system is denoted by $(x,y,z)$ and the local coordinate system for the plate
is denoted by $(x_1,x_2,x_3)$ where $(x_1,x_2)$ define the mid-surface. In order to avoid the transformation
back and forth between the two coordinate systems, we assume that they are the same.

\begin{figure}
  \centering
  \includegraphics[width=\textwidth]{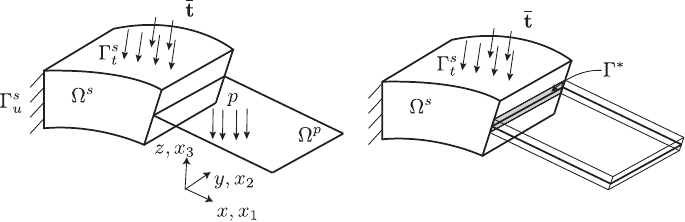}
  \caption{Coupling of a three dimensional solid and a plate.}
  \label{fig:solid-plate1}
\end{figure}

The governing equations for a solid-plate system are as follows

\begin{itemize}

\item For the solid part, cf. Equation~\eqref{eq:solid}.

\item For the plate part

\begin{subequations}
\begin{alignat}{2}
\frac{\partial^4 w}{\partial x_1^4} + 2 \frac{\partial^4 w}{\partial x_1^2\partial x_2^2} +
\frac{\partial^4 w}{\partial x_2^4} &=  \frac{p}{D}         &\quad\text{on} \quad \Omega^p \\ 
             w         &=  \bar{w}   &\quad\text{on} \quad \Gamma_u^p \\
             \pderiv{w}{n}         &=  \bar{\theta}   &\quad\text{on} \quad \Gamma_\theta^p 
\end{alignat}
\end{subequations}
where $w(x_1,x_2)$ denotes the transverse displacement field of the plate;
$D$ is  the plate rigidity which is given by $D=\frac{Eh^3}{12(1-\nu^2)}$ with $h$ being the
plate thickness and $E,\nu$ denote the Young's modulus and Poisson's ratio, respectively.
$p$ denotes the distributed pressure load; $n$ denotes the normal of the boundary $\Gamma_\theta^p$.
The prescribed transverse displacement and rotations are represented by $\bar{w}$ and $\bar{\theta}$, 
respectively. Note that for simplicity we have omitted force/moment boundary conditions since the classical
plate theory we are using is standard, we refer to \cite{hughes-fem-book,taylor-fem-book} for details.

\item For the coupling part

\begin{subequations}
\begin{alignat}{2}
  \vm{u}^s &=  \vm{u}^p                        &\quad\text{on} \quad \Gamma^* \\
  \bsym{\sigma}^s \cdot \vm{n}^s &= \bsym{\sigma}^p \cdot \vm{n}^s  &\quad\text{on} \quad \Gamma^* \label{eq:tr}
\end{alignat}
\label{eq:coupling-plate}
\end{subequations}
where $\vm{n}^s$ is the outward unit normal vector to the coupling interface $\Gamma^*$. And $\vm{u}^p$ is
the displacement field of any point in the plate (not just in the mid-surface). According to the Kirchhoff plate
theory, it is given by

\begin{equation}
\vm{u}^p=
\begin{bmatrix}
-x_3w_{,1}\\
-x_3w_{,2}\\
w(x_1,x_2)
\end{bmatrix}
\label{eq:Kirchhoff-disp}
\end{equation}
where $w_{,1}=\pderiv{w}{x_1}$.
The strain field is then given by

\begin{equation}
\begin{split}
\epsilon_{11}  &= -x_3 \frac{\partial^2 w}{\partial x_1^2} =-x_3w_{,11}\\
\epsilon_{22}  &= -x_3 \frac{\partial^2 w}{\partial x_2^2} =-x_3w_{,22}\\
2\epsilon_{12} &= -2x_3 \frac{\partial^2 w}{\partial x_1 x_2} =-2x_3w_{,12}\\
\end{split}
\label{eq:Kirchhoff-strain}
\end{equation}
and the stress field of the plate $\bsym{\sigma}^p$ is written as
\begin{equation}
\bsym{\sigma}^p\equiv
\begin{bmatrix}
\sigma_{11}\\
\sigma_{22}\\
\sigma_{12}\\
\end{bmatrix}=
\frac{E}{(1-\nu^2)}\begin{bmatrix}
1 & \nu & 0\\
\nu & 1 & 0\\
0 & 0 & 0.5(1-\nu)
\end{bmatrix}
\begin{bmatrix}
\epsilon_{11}\\
\epsilon_{22}\\
2\epsilon_{12}\\
\end{bmatrix} \equiv \vm{C}^p\bsym{\epsilon}^p
\label{eq:Kirchhoff-stress}
\end{equation}

\end{itemize}

\section{Weak forms}\label{weak-form}

In our method, interfacial conditions, Equations~\eqref{eq:coupling-beam} and \eqref{eq:coupling-plate},
are enforced weakly with Nitsche’s method \cite{nitsche}.

\subsection{Solid-beam coupling}

We start by defining the spaces, $\bsym{S}^s$ and $\vm{V}^s$ over the solid domain that will contain 
the solution and trial functions respectively:

\begin{equation}
\begin{split}
\bsym{S}^s&=\{\vm{u}^s(\vm{x})|\vm{u}^s(\vm{x}) \in \bsym{H}^1(\Omega^s), \vm{u}^s=\bar{\vm{u}} \;\;
\text{on $\Gamma_u^s$} \}\\
\bsym{V}^s&=\{\vm{v}^s(\vm{x})|\vm{v}^s(\vm{x}) \in \bsym{H}^1(\Omega^s), \vm{v}^s={\vm{0}} \;\;\text{on $\Gamma_u^s$} \}
\end{split}
\end{equation}
where $\bsym{H}^m(\Omega^{s/b})$ denotes the $m$th order Hilbert space.

In the same manner,  we define the spaces, $\bsym{S}^b$ and $\vm{V}^b$ over the beam domain that will contain 
the solution and trial functions

\begin{equation}
\begin{split}
\bsym{S}^b &=\{w(x)|w(x) \in \bsym{H}^2(\Omega^b), w=\bar{{w}} \;\;
\text{on $\Gamma_u^b$}\;\; \text{and}\; \frac{dw}{dx} =\bar{\theta} \; \text{on $\Gamma^b_\theta$}  \}\\
\bsym{V}^b &=\{v^b(x)|v^b(x) \in \bsym{H}^2(\Omega^b), {v}^b=0 \;\;
\text{on $\Gamma_u^b$}\;\; \text{and}\; \frac{dv^b}{dx} =0 \; \text{on $\Gamma^b_\theta$}  \}\\
\end{split}
\end{equation}

The standard application of Nitsche's method for the coupling is: 
Find $(\vm{u}^s,w) \in \bsym{S}^s \times \bsym{S}^b$ such that

\begin{multline}
\int_{\Omega^s} \bsym{\epsilon}(\vm{v}^s):\bsym{\sigma}^s   \mathrm{d}\Omega +
\int_{\Omega^b} EI w_{,xx} v^b_{,xx}   \mathrm{d}\Omega
-\int_{\Gamma^*}  \left(\jump{\vm{v}} \otimes \vm{n}^s\right) : \{\bsym{\sigma}\} \mathrm{d}\Gamma 
-\int_{\Gamma^*}  \left(\jump{\vm{u}} \otimes \vm{n}^s\right) : \{\bsym{\sigma}(\vm{v})\} \mathrm{d}\Gamma \\+
\int_{\Gamma^*}  \alpha \jump{\vm{v}} \cdot \jump{\vm{u}} \mathrm{d}\Gamma 
=  \int_{\Omega^s} \vm{v}^s\cdot\vm{b} \mathrm{d}\Omega + 
\int_{\Gamma_t^s} \vm{v}^s\cdot \bar{\vm{t}}  \mathrm{d}\Gamma 
+ \int_{\Omega^b} v^b p \di \Omega + \left( \frac{dw}{dx}\bar{m} \right)\bigg\lvert_{\Gamma_m^b}
\label{eq:solid-beam-weakform}
\end{multline}
for all $(\vm{v}^s,v^b) \in \bsym{V}^s \times \bsym{V}^b$.

In Equation~\eqref{eq:solid-beam-weakform}, the jump and average operators $\jump{\cdot}$
and $\{\cdot\}$ are defined as

\begin{equation}
\jump{\vm{u}} = \vm{u}^s - \vm{u}^b, \quad
 \{\bsym{\sigma}\} = \frac{1}{2}(\bsym{\sigma}^s + \bsym{\sigma}^b)
   \label{eq:jump-average}
\end{equation}
which are quantities evaluated only along the coupling interface $\Gamma^*$.

Note that the first two terms in the left hand side of 
Equation~\eqref{eq:solid-beam-weakform} are standard and the last three terms are the extra terms that take into
account the coupling at $\Gamma^*$.
The last term in the left hand side of 
Equation~\eqref{eq:solid-beam-weakform} is the so-called stabilisation term that aims to stabilise the method.
Therefore $\alpha$ is a penalty-like stabilisation term of which  
there exists a minimum  that guarantees stability, see \eg \cite{Griebel}. 
In Section \ref{sec:numerical-analysis}, a numerical analysis is provided to determine this minimum.

\subsection{Solid-plate coupling}

We define the spaces, $\bsym{S}^p$ and $\vm{V}^p$ over the plate domain that will contain 
the solution and trial functions

\begin{equation}
\begin{split}
\bsym{S}^p &=\{w(\vm{x})|w(\vm{x}) \in \bsym{H}^2(\Omega^p), w=\bar{{w}} \;\;
\text{on $\Gamma_u^p$}\;\; \text{and}\; \pderiv{w}{n} =\bar{\theta} \; \text{on $\Gamma^p_\theta$}  \}\\
\bsym{V}^p &=\{v^p(\vm{x})|v^p(\vm{x}) \in \bsym{H}^2(\Omega^p), {v}^p=0 \;\;
\text{on $\Gamma_u^p$}\;\; \text{and}\; \pderiv{v^p}{n} =0 \; \text{on $\Gamma^p_\theta$}  \}\\
\end{split}
\end{equation}

The standard application of Nitsche's method for the coupling is: 
Find $(\vm{u}^s,w) \in \bsym{S}^s \times \bsym{S}^p$ such that

\begin{multline}
\int_{\Omega^s} \bsym{\epsilon}(\vm{v}^s):\bsym{\sigma}^s   \mathrm{d}\Omega +
\int_{\Omega^p} \bsym{\epsilon}({v}^p):\bsym{\sigma}^p   \mathrm{d}\Omega
-\int_{\Gamma^*}  \left(\jump{\vm{v}} \otimes \vm{n}^s\right) : \{\bsym{\sigma}\} \mathrm{d}\Gamma 
-\int_{\Gamma^*}  \left(\jump{\vm{u}} \otimes \vm{n}^s\right) : \{\bsym{\sigma}(\vm{v})\} \mathrm{d}\Gamma \\+
\int_{\Gamma^*}  \alpha \jump{\vm{v}} \cdot \jump{\vm{u}} \mathrm{d}\Gamma 
=  \int_{\Omega^s} \vm{v}^s\cdot\vm{b} \mathrm{d}\Omega+ 
\int_{\Gamma_t^s} \vm{v}^s \cdot \bar{\vm{t}}  \mathrm{d}\Gamma 
\label{eq:solid-plate-weakform}
\end{multline}
for all $(\vm{v}^s,v^p) \in \bsym{V}^s \times \bsym{V}^p$.

\section{Discretisation}\label{discretisation}

In principle any spatial discretisation methods can be adopted, NURBS-based isogeometric finite elements
are, however, utilized in this work because (1) IGA facilitates the construction of rotation-free bending elements
and (2) IGA is also the numerical framework for our recent works on failure analysis of composite laminates 
\cite{Nguyen2013,nguyen_cohesive_2013,nguyen-offset} to which the presented coupling formulation will be
applied in a future work. For simplicity, only B-splines are presented, the formulation is, however, general and 
can be equally applied to NURBS.

\subsection{Basis functions}

We briefly present the essentials of B-splines in this section. More details can be found in the textbook
\cite{piegl_book}.

\subsubsection{B-splines basis functions and solids}

Given a knot vector $\Xi=\{\xi_1,\xi_2,\ldots,\xi_{n+p+1}\}$, $\xi_i \leq \xi_{i+1}$ where $\xi_i$ is the 
\textit{i}th knot, $n$ is the number of basis functions and $p$ is 
the polynomial order.     
The associated set of B-spline basis functions $\{N_{i,p}\}_{i=1}^n$ are
defined recursively by the Cox-de-Boor formula, starting with the zeroth order basis
function ($p=0$)
\begin{equation}
  N_{i,0}(\xi) = \begin{cases}
    1 & \textrm{if $ \xi_i \le \xi < \xi_{i+1}$},\\
    0 & \textrm{otherwise}
  \end{cases}
  \label{eq:basis-p0}
\end{equation}
and for a polynomial order $p \ge 1$
\begin{equation}
  N_{i,p}(\xi) = \frac{\xi-\xi_i}{\xi_{i+p}-\xi_i} N_{i,p-1}(\xi)
               + \frac{\xi_{i+p+1}-\xi}{\xi_{i+p+1}-\xi_{i+1}}
	       N_{i+1,p-1}(\xi)
  \label{eq:basis-p}
\end{equation}
in which fractions of the form $0/0$ are
defined as zero.

Fig.~\ref{fig:bspline-quad-open} illustrates a corresponding set of quadratic basis functions for an open, non-uniform knot vector. Of particular note is the interpolatory nature of the basis function at the two
ends of the interval created through an open knot vector, and the reduced continuity at $\xi = 4$ due to the presence of the location of a repeated knot where $C^0$ continuity is attained. 
At other knots, the functions are $C^1$ continuous ($C^{p-1}$). In an analysis context, non-zero knot spans ($\xi_i,\xi_{i+1}$ is a knot span) play the role of elements. Thus, knots $\xi_i$ are the element boundaries and therefore B-spline basis functions are $C^{p-1}$ across the element boundaries. This is a key difference compared to standard Lagrange finite elements. In this regard, B-spline basis functions are similar to meshless
approximations see \eg \cite{nguyen_meshless_2008}.

\begin{figure}[h!]
  \centering 
  \includegraphics[width=0.7\textwidth]{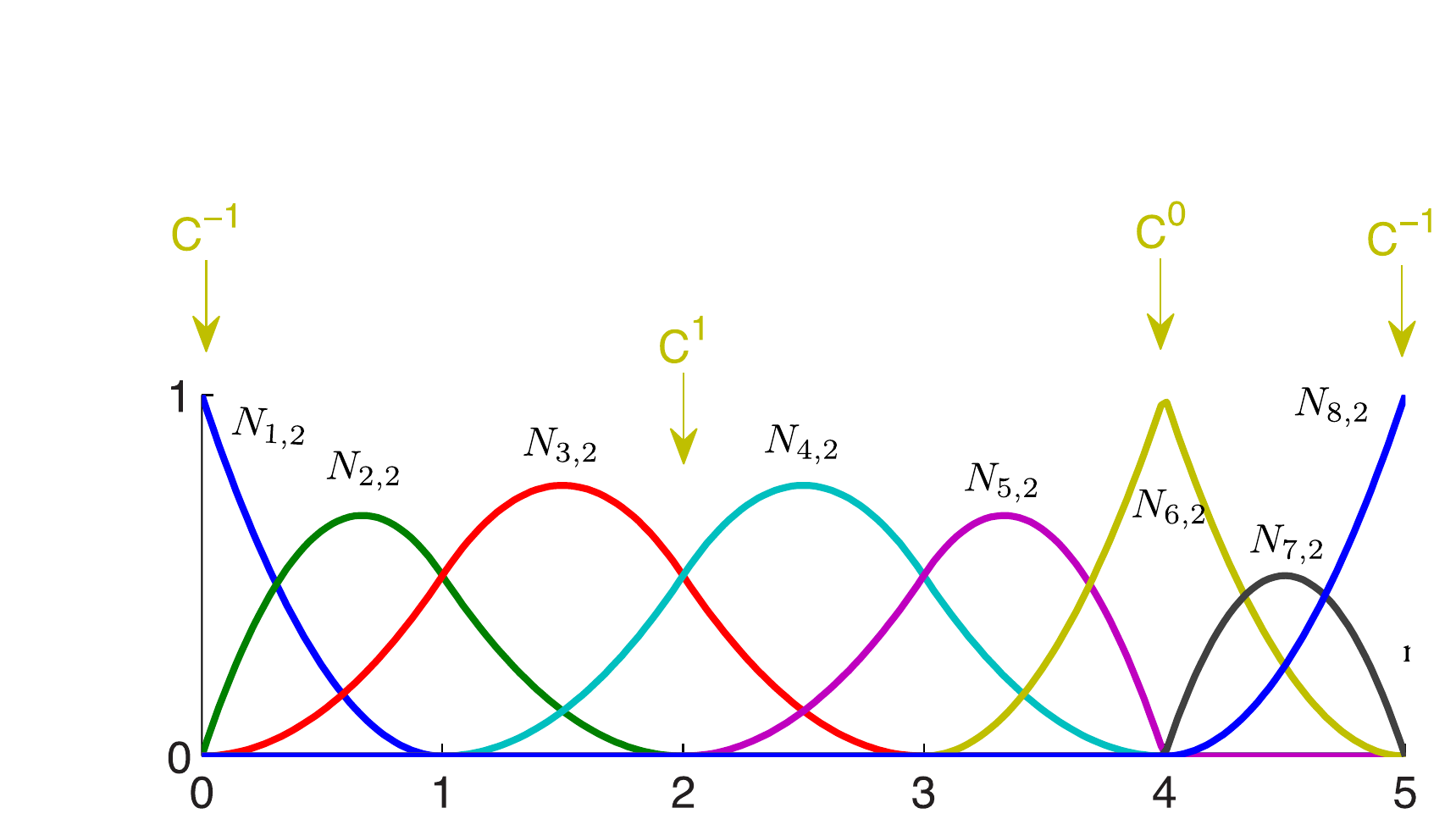}
  \caption{Quadratic B-spline basis functions defined for the open, non-uniform knot vector
  $\Xi=\{0,0,0,1,2,3,4,4,5,5,5\}$. Note the flexibility in the construction of
  basis functions with varying degrees of regularity.} 
  \label{fig:bspline-quad-open} 
\end{figure}

Let $\Xi^1=\{\xi_1,\xi_2,\ldots,\xi_{n+p+1}\}$, 
    $\Xi^2=\{\eta_1,\eta_2,\ldots,\eta_{m+q+1}\}$,
and $\Xi^3=\{\zeta_1,\zeta_2,\ldots,\zeta_{l+r+1}\}$ are the knot vectors
    and a control net $\vm{P}_{i,j,k} \in \mathds{R}^{d_s}$.
A tensor-product B-spline solid is defined as 

\begin{equation}
  \label{eq:bspline_volume}
  	\vm{V}(\xi,\eta,\zeta) =
	\sum_{i=1}^{n}\sum_{j=1}^{m}\sum_{k=1}^l N_{i,p}(\xi)M_{j,q}(\eta) L_{k,r}(\zeta)\vm{P}_{i,j,k}
\end{equation}
or, by defining a global index $A$ through
\begin{equation}
  \label{eq:bspline_volume_mapping}
  A = (n \times m) ( k - 1)  + n( j - 1 ) + i
\end{equation}
a simplified form of Equation~\eqref{eq:bspline_volume} can be written as
\begin{equation}
  \label{eq:bspline_vol_simple}
  	\vm{V}(\boldsymbol{\xi}) = \sum_{A=1}^{n \times m \times l} \vm{P}_A  N_{A}^{p,q,r}(\boldsymbol{\xi} )
\end{equation}

\subsubsection{Isogeometric analysis}

\begin{figure}
  \centering
  \includegraphics[width=0.4\textwidth]{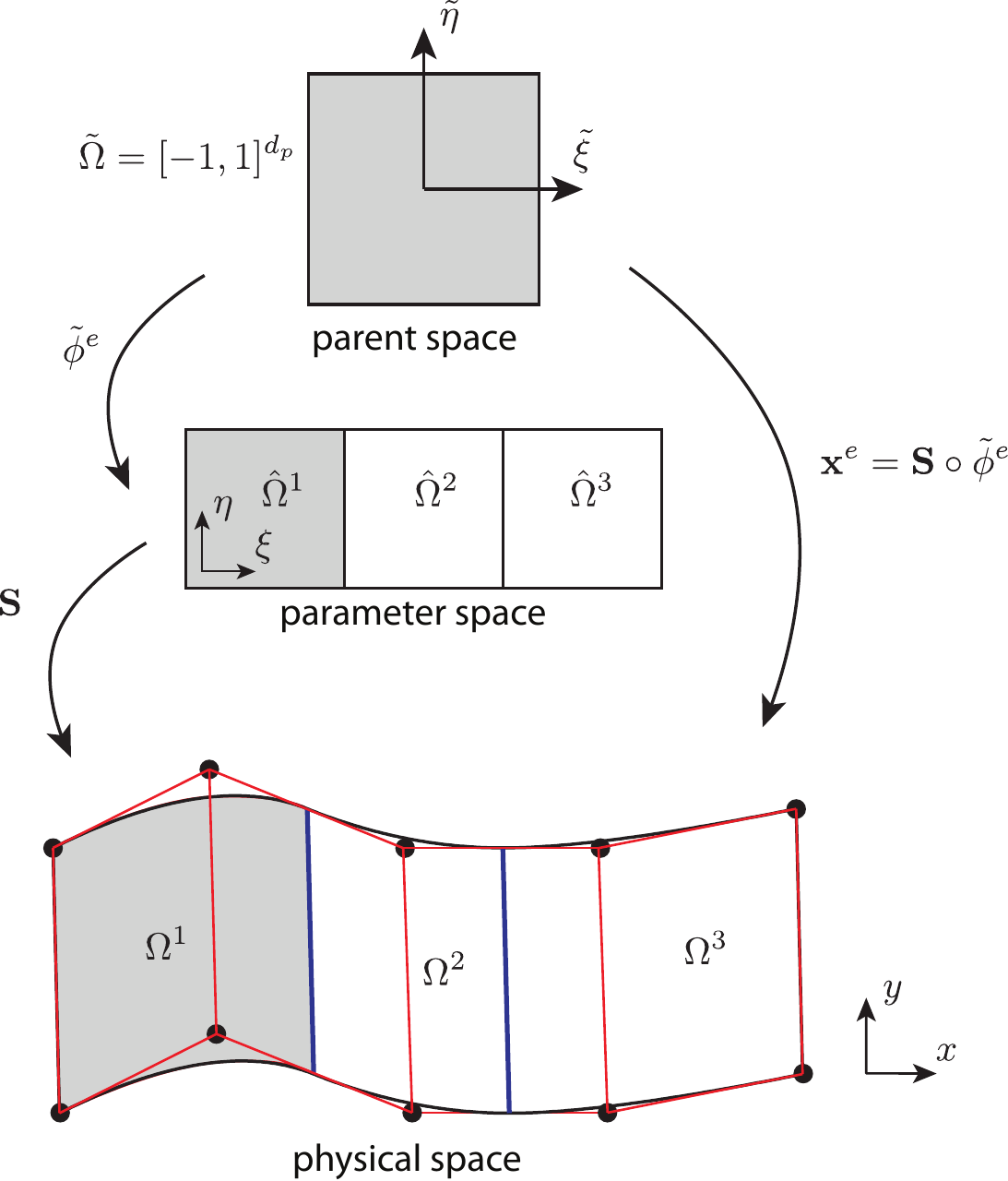}
  \caption{Diagrammatic interpretation of mappings from parent space ($\tilde{\Omega}$) 
     through parametric space ($\hat{\Omega}$) to physical space ($\Omega$). The parent space is where
     numerical quadrature rules are defined.}
  \label{fig:iga_mappings}
\end{figure}

Isogeometric analysis also makes use of an isoparametric formulation, but a key difference over its Lagrangian counterpart is the use of basis functions generated by CAD to discretise both the geometry and unknown fields. 
In IGA, regions bounded by knot lines with non-zero parametric area lead to a natural definition of 
element domains.
The use of NURBS basis functions for discretisation introduces the concept of parametric space which is absent in conventional FE implementations. The consequence of this additional space is that an additional mapping must be performed to operate in parent element coordinates. As shown in Fig.~\ref{fig:iga_mappings}, two mappings are considered for IGA with NURBS: a mapping $\tilde{\phi}^e: \tilde{\Omega} \to \hat{\Omega}^e$ and $\vm{S}: \hat{\Omega} \to \Omega$. The mapping $\vm{x}^e: \tilde{\Omega} \to \Omega^e$ is given by the composition $\vm{S}\circ \tilde{\phi}^e$. 

For a given element $e$, the geometry is expressed as
\begin{equation}
  \label{eq:iga_geometry_discretisation}
  \mathbf{x}^e(\tilde{\boldsymbol{\xi}}) = \sum_{a=1}^{n_{en}} \vm{P}_a^e N_a^e(\tilde{\boldsymbol{\xi}})
\end{equation}
where $a$ is a local basis function index, $n_{en} = (p+1)^{d_p}$ is the number of non-zero basis functions over element $e$ and $\vm{P}_a^e$,$N_a^e$ are the control point and B-spline basis function  associated with index $a$ respectively. We employ the commonly used notation of an element connectivity mapping \cite{hughes-fem-book} which translates a local basis function index to a global index through
\begin{equation}
  \label{eq:element_connectivity_array}
  A = \textrm{IEN}( a, e )
\end{equation}
Global and local control points are therefore related through $\vm{P}_A \equiv \vm{P}_{\textrm{IEN}(a,e)} \equiv \vm{P}_a^e$ with similar expressions for $R_a^e$.  A field $\vm{u}(\mathbf{x})$ which governs our relevant PDE can also be discretised in a similar manner to Equation~\eqref{eq:iga_geometry_discretisation} as
\begin{equation}
  \label{eq:iga_field_discretisation}
  \vm{u}^e(\tilde{\boldsymbol{\xi}}) =  \sum_{a=1}^{n_{en}} \vm{d}_a^e N_a^e(\tilde{\boldsymbol{\xi}})
\end{equation}
where $\vm{d}^e_a$ represents a control (nodal) variable. In contrast to conventional discretisations, these coefficients are not in general interpolatory at nodes. This is similar to the case of meshless
methods built on non-interpolatory shape functions such as the moving least squares
(MLS) \cite{efg-nayroles,NME:NME1620370205,nguyen_meshless_2008}. Using the Bubnov-Galerkin method, an 
expansion analog to Equation~\eqref{eq:iga_field_discretisation} is adopted for the weight function and upon substituting
them into a weak form, a standard system of linear equations is obtained from which $\vm{d}$--the nodal variables
are obtained.

\subsection{Solid-beam coupling}\label{sec:solid-beam}

The solid part and the beam part are discretised into finite elements, cf. Fig.~\ref{fig:domain-mesh}.
In what follows discussion is made on the discretisation of the solid, of the beam and the coupling interface.

Based on the weak form given in Equation~\eqref{eq:solid-beam-weakform}, the discrete equation reads

\begin{equation}
\left( \vm{K}^s + \vm{K}^b + \vm{K}^n + (\vm{K}^n)\trans + \vm{K}^{st} \right) \vm{a} = \vm{f}
\label{eq:general-kuf}
\end{equation}
where $\vm{K}^{s/b}$ are the stiffness matrices of the solid and beam, respectively;
$\vm{K}^{st}$ and $\vm{K}^n$ are the coupling matrices. The superscript $st$ indicates that $\vm{K}^{st}$
is the stabilisation coupling matrix (involving $\alpha$). In the above, $\vm{a}$ represents the nodal
displacement vector. The external force vector is designated by $\vm{f}$. Note that by removing $\vm{K}^n$
and its transpose, Equation~\eqref{eq:general-kuf} reduces to the standard penalty method. Thus, implementing
the Nitsche's method is very identical to the penalty method and usually straightforward.

\subsubsection{Discretisation of the solid}\label{sec:solid}

The solid discretisation is standard and the corresponding stiffness matrix is given by

\begin{equation}
\vm{K}^s=\int_{\Omega^s} (\vm{B}^s)\trans \vm{C}^s \vm{B}^s \di \Omega = 
\bigcup_{e=1}^{nels}\int_{\Omega_e^s} (\vm{B}^s_e)\trans \vm{C}^s \vm{B}^s_e \di \Omega 
\end{equation}
where $nels$ denotes the number of solid elements of $\Omega^s$ and $\bigcup$ denotes the
standard assembly operator; $\vm{B}^s$ is the standard strain-displacement matrix and $\vm{C}^s$ is the elasticity
matrix. For two dimensional element $e$, $\vm{B}^s_e$ is given by

\begin{equation}
\vm{B}_e^s = \begin{bmatrix}
N_{1,x}^s & 0 & N_{2,x}^s & 0 & \ldots\\
0 & N_{1,y}^s & 0 & N_{2,y}^s & \ldots \\
N_{1,y}^s & N_{1,x}^s & N_{2,y}^s & N_{2,x}^s & \ldots 
\end{bmatrix}
\end{equation}
Expressions for three dimensional elements can be found in many FEM textbooks \eg
\cite{hughes-fem-book}. The notation $N_{I,x}$ denotes
the derivative of shape function $N_I$ with respect to $x$. This notation for partial derivatives
will be used in subsequent sections.

\begin{figure}
  \centering
  \includegraphics[width=0.4\textwidth]{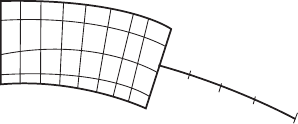}
  \caption{Finite element discretisations of the solid and the beam.}
  \label{fig:domain-mesh}
\end{figure}

The external force vector due to the solid part, if exists, is given by
(for sake of simplicity, body force was omitted)

\begin{equation}
\vm{f}^s = \int_{\Gamma^s_t} (\vm{N}^s)\trans \bar{\vm{t}} \di \Gamma
\end{equation}
where $\vm{N}^s$ denotes the matrix of shape functions. For 2D, it is given by
\begin{equation}
\vm{N}^s = \begin{bmatrix}
N_1^s & 0 & N_2^s & 0 & \ldots\\
0 & N_1^s & 0 & N_2^s & \ldots       
\end{bmatrix}
\end{equation}

\subsubsection{Discretisation of the beam}

Since NURBS facilitate the implementation of rotation free beam elements, we use NURBS to discretise
the beam transverse displacement field 

\begin{equation}
w = N_I w_I
\label{eq:beam-w}
\end{equation}
where $w_I$ denotes the control point transverse displacement.
The resulting stiffness matrix reads \cite{taylor-fem-book}

\begin{equation}
\vm{K}^b = \int_{\Omega^b} EI \vm{B}^b (\vm{B}^b)\trans \di \Omega=
\bigcup_{e=1}^{nelb} \int_{\Omega^b_e} EI \vm{B}^b_e (\vm{B}^b_e)\trans \di \Omega
\end{equation}
where $nelb$ denotes the number of beam elements and
$\vm{B}_e^b$ denotes a vector containing the second derivatives
of the shape functions associated with element $e$ with respect to $\bar{x}$. 

\subsubsection{Coupling matrices}

In order to focus on the coupling itself, we assume that the coordinate system of the beam is identical
to the global one. Therefore rotation and transformation matrices are unnecessary. We postpone the treatment
of a general case to Section \ref{sec:rotation}.

The displacement field of the beam defined for a point on the coupling interface is given by, cf. 
Equation~\eqref{eq:beam-disp}

\begin{equation}
\begin{bmatrix}
u_{x}^b\\
u_{y}^b
\end{bmatrix}=
\begin{bmatrix}
-y N_{I,x}\\
N_{I}      
\end{bmatrix}w_I \equiv \vm{N}^{b} \vm{a}^b
\end{equation}
And the strain field defined in Equation~\eqref{eq:beam-stress-strain} can be written as
\begin{equation}
\bsym{\epsilon}^b   =
\begin{bmatrix}
\epsilon_{xx}^b  \\
\epsilon_{yy}^b  \\
2\epsilon_{xy}^b  
\end{bmatrix}=
\begin{bmatrix}
-yN_{I,xx}\\0\\0
\end{bmatrix}w_I \equiv \vm{B}^{b,c} \vm{a}^b
\end{equation}
the superscript $c$ in the strain-displacement matrix $\vm{B}^{b,c}$ is to differentiate it from
$\vm{B}^{b}$ which is also a strain-displacement matrix.

The jump and average operators defined in Equation~\eqref{eq:jump-average} are thus given by

\begin{equation}
\begin{split}
\jump{\vm{u}} &= \vm{N}^s(\vm{x})\vm{a}^s - \vm{N}^b(\vm{x}) \vm{a}^b\\
\{\bsym{\sigma}\} &= \frac{1}{2}\left(  \vm{C}^s \vm{B}^s \vm{a}^s + \vm{C}^b \vm{B}^{b,c} \vm{a}^b  \right)     
\end{split}
\label{eq:jump-average-discrete}
\end{equation}
where $\vm{C}^b$ was defined in Equation~\eqref{eq:beam-stress-strain}.

Using Equation~\eqref{eq:jump-average-discrete} and from the last three terms in the left hand side of the
weak form given in Equation~\eqref{eq:solid-beam-weakform}, one obtains the following expression for the
coupling matrices

\begin{equation}
\vm{K}^n = \begin{bmatrix}[2.5]
-\D\int_{\Gamma^*} \vm{N}^{s\text{T}} \vm{n} \frac{1}{2}\vm{C}^s\vm{B}^s \mathrm{d}\Gamma &
-\D\int_{\Gamma^*} \vm{N}^{s\text{T}} \vm{n} \frac{1}{2}\vm{C}^b\vm{B}^{b,c} \mathrm{d}\Gamma \\
\D\int_{\Gamma^*} \vm{N}^{b\text{T}} \vm{n} \frac{1}{2}\vm{C}^s\vm{B}^s \mathrm{d}\Gamma &
\D\int_{\Gamma^*} \vm{N}^{b\text{T}} \vm{n} \frac{1}{2}\vm{C}^b\vm{B}^{b,c} \mathrm{d}\Gamma 
\end{bmatrix}\label{eq:nitsche-kdg}
\end{equation}
and by
\begin{equation}
\vm{K}^{st} = \begin{bmatrix}[2.5]
\D\int_{\Gamma^*}  \alpha  \vm{N}^{s\text{T}} \vm{N}^s \mathrm{d}\Gamma &
 - \D\int_{\Gamma^*}  \alpha  \vm{N}^{s\text{T}} \vm{N}^b \mathrm{d}\Gamma\\
 - \D\int_{\Gamma^*}  \alpha  \vm{N}^{b\text{T}} \vm{N}^s \mathrm{d}\Gamma&
 \D\int_{\Gamma^*}  \alpha  \vm{N}^{b\text{T}} \vm{N}^b \mathrm{d}\Gamma
\end{bmatrix}\label{eq:nitsche-kpe}
\end{equation}
And the normal vector in a matrix notation is represented as the following matrix 
\begin{equation}
\vm{n} = \begin{bmatrix}
n_x & 0 & n_y \\ 0 & n_y & n_x
\end{bmatrix}
\label{eq:n-matrix}
\end{equation}
\subsubsection{Implementation of coupling matrices}\label{sec:implementation1}

Computation of the coupling matrices involves integrals over the coupling interface $\Gamma^*$
of which integrands are the shape functions and its derivatives which are zero everywhere except
at elements that intersect $\Gamma^*$. We denote $\Omega^s_b$ a set of solid elements that 
intersect with $\Gamma^*$ and $\Omega^b_b$ the beam element that intersects $\Gamma^*$, see
  Fig.~\ref{fig:solid-beam-impl}. The subscript $b$ indicates elements on the coupling boundary $\Gamma^*$.
\begin{figure}[htbp]
  \centering 
   \includegraphics[width=0.55\textwidth]{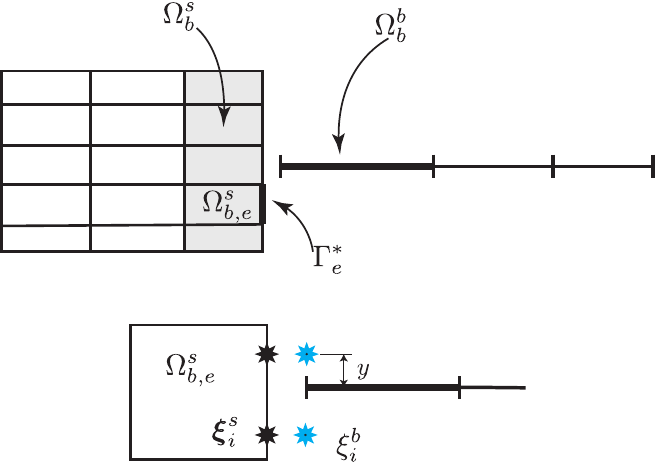}
  \caption{Coupling a beam and a continuum using a Nitsche based method.}
  \label{fig:solid-beam-impl}
\end{figure}

We then use the trace mesh of $\Omega^s_b$ on $\Gamma^*$ to compute the coupling matrices.
For example, the following term taken from $\vm{K}^{st}$ are computed as

\begin{equation}
 \int_{\Gamma^*} \alpha  \vm{N}^{s\text{T}} \vm{N}^b \mathrm{d}\Gamma=
 \bigcup_{e=1}^{nce}\int_{\Gamma^*_e} \alpha  \vm{N}^{s\text{T}} \vm{N}^b \mathrm{d}\Gamma=
 \bigcup_{e=1}^{nce}\sum_{i=1}^{ngp} \alpha  \vm{N}^{s\text{T}}(\bsym{\xi}_i^s) \vm{N}^b (y,\xi_i^b)w_i^s
\end{equation}
where $nce$ (number of coupling elements) denotes the number of $\Omega^s_{b,e}$; 
$ngp$ is the number of Gauss points (GPs) and $w_i^s$
are the weights.
This matrix will be assembled to the rows of the global stiffness matrix using the connectivity array
of solid element $\Omega^s_{b,e}$ and to the columns corresponding to the beam element $\Omega^b_b$.

\subsection{Stabilisation parameter}\label{sec:numerical-analysis}

The Nitsche's bilinear form in the weak form given in Equation~\eqref{eq:solid-beam-weakform} can be written as

\begin{equation}
a(\vm{u},\vm{v})=\tilde{a}(\vm{u},\vm{v}) + \alpha \int_{\Gamma^*} \jump{\vm{u}}\cdot\jump{\vm{v}}\di \Gamma  
-\int_{\Gamma^*}  \left(\jump{\vm{v}} \otimes \vm{n}^s\right) : \{\bsym{\sigma}\} \mathrm{d}\Gamma 
-\int_{\Gamma^*}  \left(\jump{\vm{u}} \otimes \vm{n}^s\right) : \{\bsym{\sigma}(\vm{v})\} \mathrm{d}\Gamma 
\label{eq:eq1}
\end{equation}
where $\tilde{a}(\vm{u},\vm{v})$ denotes the bilinear form of the bulk \ie without the interfacial terms

\begin{equation}
\tilde{a}(\vm{u},\vm{v}) =
\int_{\Omega^s} \bsym{\epsilon}(\vm{v}^s):\bsym{\sigma}^s   \mathrm{d}\Omega +
\int_{\Omega^b} EI v_{,xx} u_{,xx}   \mathrm{d}\Omega
\end{equation}

The problem is to find $\alpha$ such that the bilinear form is coercive \ie $a(\vm{u},\vm{u})\ge c^{te}
\norm{\vm{u}}$. One can write $a(\vm{u},\vm{u})$ as follows

\begin{equation}
a(\vm{u},\vm{u})=\tilde{a}(\vm{u},\vm{u}) + \alpha \int_{\Gamma^*} \jump{\vm{u}}\cdot\jump{\vm{u}}d\Gamma -  
\int_{\Gamma^*} \jump{\vm{u}} \otimes \vm{n}^s : \left[ \vm{C}^s : \bsym{\epsilon}^s  + \vm{C}^b :\bsym{\epsilon}^b\right] \di \Gamma
\label{eq:eq1}
\end{equation}
and the last term can be rewritten as 
\begin{equation}
\int_{\Gamma^*} \jump{\vm{u}}  \otimes \vm{n}^s :  \left[ \vm{C}^s : \bsym{\epsilon}^s + 
\vm{C}^b : \bsym{\epsilon}^b \right] \di \Gamma = (\jump{\vm{u}}, \bar{\vm{t}})_{\Gamma^*}
\end{equation}
where  $(\cdot,\cdot)_{\Gamma^*}$ is the $L_2(\Gamma^*)$ scalar product.

For any scalar product and any parameter $\epsilon$, one has

\begin{equation}
\norm{\vm{u}-\epsilon\vm{v}}^2=\norm{\vm{u}}^2 + \epsilon^2 \norm{\vm{v}}^2 - 2\epsilon (\vm{u},\vm{v})
\end{equation}
where $\norm{\cdot}$ is the norm equipped with $(\cdot,\cdot)$ \ie $\norm{\cdot}=\sqrt{(\cdot,\cdot)}$. 
From this equation, one can write

\begin{equation}
(\vm{u},\vm{v}) = \frac{1}{2\epsilon}\norm{\vm{u}}^2 + \frac{\epsilon}{2}\norm{\vm{v}}^2 - \norm{\vm{u}-\epsilon\vm{v}}^2 \le \frac{1}{2\epsilon}\norm{\vm{u}}^2 + \frac{\epsilon}{2}\norm{\vm{v}}^2
\label{eq:bdt1}
\end{equation}
If we apply Equation~\eqref{eq:bdt1} to $(\jump{\vm{u}},\bar{\vm{t}})_{\Gamma^*}$, we have

\begin{equation}
 (\jump{\vm{u}},\bar{\vm{t}})_{\Gamma^*} \le
\frac{1}{2\epsilon} \norm{\jump{\vm{u}}}^2_{\Gamma^*}
 + \frac{\epsilon}{2} \norm{\bar{\vm{t}}}^2_{\Gamma^*}
\label{eq:eq2}
\end{equation}
Substituting Equation~\eqref{eq:eq2} into Equation~\eqref{eq:eq1} leads to the following

\begin{equation}
a(\vm{u},\vm{u}) \ge \tilde{a}(\vm{u},\vm{u}) + \left(\alpha -\frac{1}{2\epsilon}\right) 
    \norm{\jump{\vm{u}}}^2_{\Gamma^*}
 -  \frac{\epsilon}{2}\norm{\bar{\vm{t}}}^2_{\Gamma^*}
\label{eq:eq3}
\end{equation}
The main idea of the method is to suppose the existence of a configuration-dependent constant $C$ such that
\begin{equation}
\norm{\bar{\vm{t}}}^2_{\Gamma^*} \le C^2 \tilde{a}(\vm{u},\vm{u})
\label{eq:eq4}
\end{equation}
which means that the flux are bounded by the gradient inside the domain.

Substituting Equation~\eqref{eq:eq4} into Equation~\eqref{eq:eq3} leads to the following

\begin{equation}
a(\vm{u},\vm{u}) \ge \left( 1 - \frac{C^2 \epsilon}{2}\right) \tilde{a}(\vm{u},\vm{u}) + 
\left(\alpha -\frac{1}{2\epsilon}\right) \norm{\jump{\vm{u}}}^2
\label{eq:eq5}
\end{equation}
which indicates that the bilinear form is coercive if the following inequalities are satisfied
\begin{equation}
 1 - \frac{C^2 \epsilon}{2} \ge 0, \quad \alpha -\frac{1}{2\epsilon} \ge 0
\label{eq:eq6}
\end{equation}
If we take $\epsilon=1/C^2$, the first inequality in Equation~\eqref{eq:eq6} will be satisfied.
Therefore the bilinear form $a$ is positive definite if $\alpha > C^2/2$. In what follows we present
a way to determine the constant C satisfying Equation~\eqref{eq:eq4}. This is where one has to solve
an eigenvalue problem.

The discrete version of Equation~\eqref{eq:eq4} is
\begin{equation}
\int_{\Gamma^*} \left( \vm{n}^s \vm{C}^s \bsym{\epsilon}^s + \vm{n}^s \vm{C}^b \bsym{\epsilon}^b 
\right)\trans \left( \vm{C}^s \bsym{\epsilon}^s \vm{n}^s + \vm{C}^b \bsym{\epsilon}^b \vm{n}^s \right) \di \Gamma \le C^2 \left(
\int_{\Omega^s} (\bsym{\epsilon}^s)\trans \vm{C}^s \bsym{\epsilon}^s \di \Omega + 
\int_{\Omega^b} (\bsym{\epsilon}^b)\trans EI \bsym{\epsilon}^b \di \Omega 
\right)
\end{equation}
Or
\begin{equation}
\frac{\vm{a}\trans \vm{H} \vm{a}}{\vm{a}\trans \tilde{\vm{K}} \vm{a}} \le C^2 
\end{equation}
where $\tilde{\vm{K}}=\vm{K}^s+\vm{K}^b$ is the total stiffness matrix without coupling terms;
$\vm{a}$ is the unknowns vector and $\vm{H}$ is given by

\begin{equation}
\vm{H}=\begin{bmatrix}[2.5]
\D\int_{\Gamma^*} (\vm{B}^s)\trans (\vm{C}^s)\trans (\vm{n}^s)\trans \vm{n}^s \vm{C}^s \vm{B}^s \di \Gamma &
\D\int_{\Gamma^*} (\vm{B}^s)\trans (\vm{C}^s)\trans (\vm{n}^s)\trans \vm{n}^s \vm{C}^b \vm{B}^b \di \Gamma \\
\D\int_{\Gamma^*} (\vm{B}^b)\trans (\vm{C}^b)\trans (\vm{n}^s)\trans \vm{n}^s \vm{C}^s \vm{B}^s \di \Gamma &
\D\int_{\Gamma^*} (\vm{B}^b)\trans (\vm{C}^b)\trans (\vm{n}^s)\trans \vm{n}^s \vm{C}^b \vm{B}^b \di \Gamma 
\end{bmatrix}
\end{equation}

So we need to find the maximum of the Rayleigh quotient $R$ defined as $R=\frac{\vm{u}\trans \vm{H} \vm{u}}{\vm{u}\trans \tilde{\vm{K}} \vm{u}}$. By theorem this maximum is $\lambda_1$, the largest eigenvalue of 
matrix $\tilde{\vm{K}}^{-1}\vm{H}$. Finally the stabilisation parameter $\alpha$ is chosen 
according to

\begin{equation}
\alpha = \frac{\lambda_1}{2}\label{eq:alpha-eigen}
\end{equation}

\subsection{Solid-plate coupling}\label{sec:plate-continuum}

The discrete equation is identical to Equation~\eqref{eq:general-kuf} in which $\vm{K}^b$
will be replaced by $\vm{K}^p$. The discretisation of the solid was presented in Section \ref{sec:solid}.
Therefore only the discretisation of the plate and the computation of coupling matrices are presented.

\subsubsection{Discretisation of the plate}

The deflection $w$ is approximated as follows

\begin{equation}
w = N_I(\xi,\eta) w_I \label{eq:plate-w}
\end{equation}

\noindent where $N_I(\xi,\eta)$ is the NURBS basis function associated with node I
and $w_I$ is the nodal deflection. 

The plate element stiffness matrix is standard (see \eg \cite{taylor-fem-book}) and given by

\begin{equation}
\vm{K}_e^p = \int_{\Omega_e^p} (\vm{B}_e^p)\trans \vm{D}^p_b \vm{B}_e^p \di \Omega
\end{equation}

\noindent where the constitutive matrix $\vm{D}^p$ reads

\begin{equation}
\vm{D}_b^p=\frac{Eh^3}{12(1-\nu^2)}\begin{bmatrix}
1 & \nu & 0\\
\nu & 1 & 0\\
0 & 0 & 0.5(1-\nu)
\end{bmatrix}
\label{eq:bending-stiffness}
\end{equation}

\noindent and the element displacement-curvature matrix $\vm{B}_e^p$ is given by

\begin{equation}
\vm{B}_e^p = \begin{bmatrix} 
N_{1,xx} & N_{2,xx} & \cdots & N_{n,xx}\\
N_{1,yy} & N_{2,yy} & \cdots & N_{n,yy}\\
2N_{1,xy} & 2N_{2,xy} & \cdots & 2N_{n,xy}\\
\end{bmatrix}
\end{equation}
where $n$ denotes the number of basis functions of element $e$.

\subsubsection{Coupling matrices}

Using Equation~\eqref{eq:plate-w}, the displacement field in Equation~\eqref{eq:Kirchhoff-disp}
is given by

\begin{equation}
\vm{u}^p = \begin{bmatrix}
-x_3 N_{I,1}\\
-x_3 N_{I,2}\\
N_I 
\end{bmatrix}
\begin{bmatrix}
w_I
\end{bmatrix}\equiv\vm{N}_I^p\vm{a}^p_I
\end{equation}
and the strain field given in Equation~\eqref{eq:Kirchhoff-strain} is written as

\begin{equation}
\begin{bmatrix}
\epsilon_{11}\\
\epsilon_{22}\\
2\epsilon_{12}\\
\end{bmatrix}=
\begin{bmatrix}
-x_3 N_{I,11}\\
-x_3 N_{I,22}\\
-2x_3 N_{I,12}  
\end{bmatrix}
\begin{bmatrix}
w_I
\end{bmatrix}\equiv\vm{B}_{I}^{p,c}  \vm{a}_I^p
\end{equation}
The coupling matrices are thus given by Equations~\eqref{eq:nitsche-kdg} and
\eqref{eq:nitsche-kpe} where quantities with superscript $b$ are replaced by ones with superscript $p$.

\subsubsection{Implementation of coupling matrices}

Basically the implementation is similar to the solid-beam coupling presented in
Section \ref{sec:implementation1}, but much more involved due to the fact that one solid
element $\Omega^s_{b,e}$ might interact with more than one plate element $\Omega^p_{b,e}$ or
vice versa, Fig.~\ref{fig:plate-continuum} illustrates the former. In what follows, we present
implementation details for Lagrange finite elements. Extension to B-spline or NURBS elements will
be discussed later.

\begin{figure}[htbp]
  \centering 
   \includegraphics[width=0.55\textwidth]{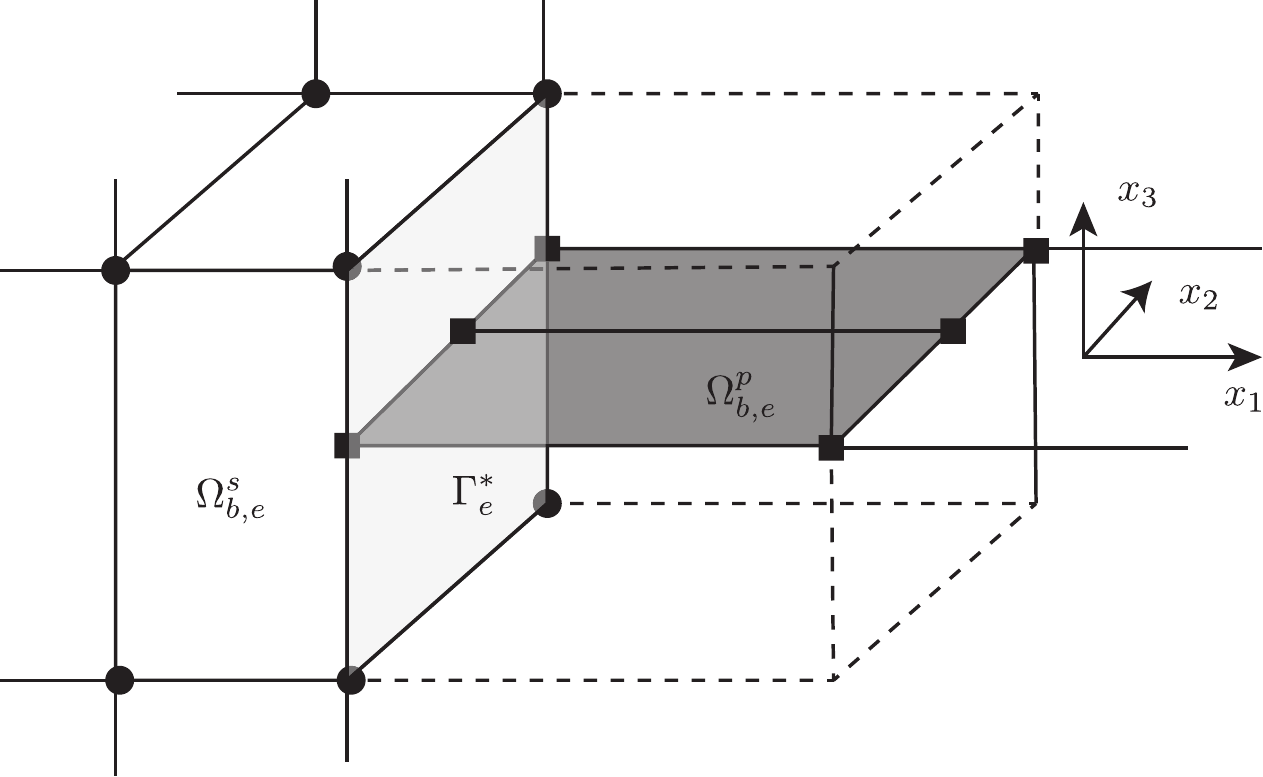}
  \caption{Coupling a plate and a 3D continuum using a Nitsche based method.}
  \label{fig:plate-continuum}
\end{figure}

For each of the solid coupling elements, $\Omega^s_{b,e}$, a number of GPs $\{\xi_i,\eta_i,w_i\}_{i=1}^{ngp}$ 
are placed on the coupling surface $\Gamma^*_e$. Those GPs are transformed to the physical space and then transformed to (1) the parent space of the solid, $[-1,1]^3$ and to (2) the parent space of the plate, $[-1,1]^2$.
The solid GPs are represented by $\bsym{\xi}_i^s=(\xi,\eta,\zeta)_i^s$ and the plate GPs by 
$\bsym{\xi}_i^p=(\xi,\eta)_i^p$.
We refer to Fig.~\ref{fig:plate-continuum-gauss}.
Concretely, one performs the following operations 
\begin{equation}
\begin{split}
\vm{x}_i &= \vm{M}(\xi_i,\eta_i)\vm{x}_{\Gamma^*_e} \\
\vm{x}_i &= \vm{N}^s(\xi_i^s,\eta_i^s,\zeta_i^s)\vm{x}_e^s \rightarrow (\xi_i^s,\eta_i^s,\zeta_i^s) \\
\vm{x}_i^* &= \vm{N}^p (\xi_i^p,\eta_i^p) \vm{x}_e^p \rightarrow (\xi_i^p,\eta_i^p)
\end{split}
\label{eq:ananan}
\end{equation}
where $\vm{x}_{\Gamma^*_e}$ denotes the coordinates of the nodes on the coupling surface $\Gamma^*_e$.
The nodal coordinates of solid and plate elements are designated by $\vm{x}_e^{s/p}$. The last two 
of Equation~\eqref{eq:ananan} are solved using a Newton-Raphson method.

Coupling terms can now be computed, for example

\begin{equation}
\int_{\Gamma^*} \alpha  \vm{N}^{s\text{T}} \vm{N}^p \mathrm{d}\Gamma
= \bigcup_{e=1}^{nbe} 
\int_{\Omega_*^e} \alpha  \vm{N}^{s\text{T}} \vm{N}^p  d \Omega =
\bigcup_{e=1}^{nbe}
\sum_{i=1}^{ngp}  \alpha  \vm{N}^{s\text{T}}(\bsym{\xi}_i^s) \vm{N}^p( \bsym{\xi}_i^p) ) \bar{w}_i
\label{eq:coupling-terms}
\end{equation}
with $\bar{w}_i = w_i \norm{\vm{a}_1 \times \vm{a}_2}$ where $\vm{a}_\alpha$ are  the tangent vectors
of $\Gamma^*_e$ defined by

\begin{equation}
\vm{a}_1 = \vm{M}_{,\xi} \vm{x}_{\Gamma^*_e},\quad
\vm{a}_2 = \vm{M}_{,\eta} \vm{x}_{\Gamma^*_e}
\label{eq:tangents}
\end{equation}
These tangents are needed to compute the outward normal vector 
$\vm{n}=\frac{\vm{a}_1 \times \vm{a}_2}{\norm{\vm{a}_1 \times \vm{a}_2}}$ and put in a matrix notation,
   we write the normal as follows

\begin{equation}
\vm{n} = \begin{bmatrix}
n_x & 0 & 0 & n_y & 0 & n_z \\
0 & n_y & 0 & n_x & n_z & 0\\
0 & 0 & n_z & 0 & n_y & n_x  
\end{bmatrix}
\end{equation}

\begin{figure}[htbp]
  \centering 
   \includegraphics[width=0.65\textwidth]{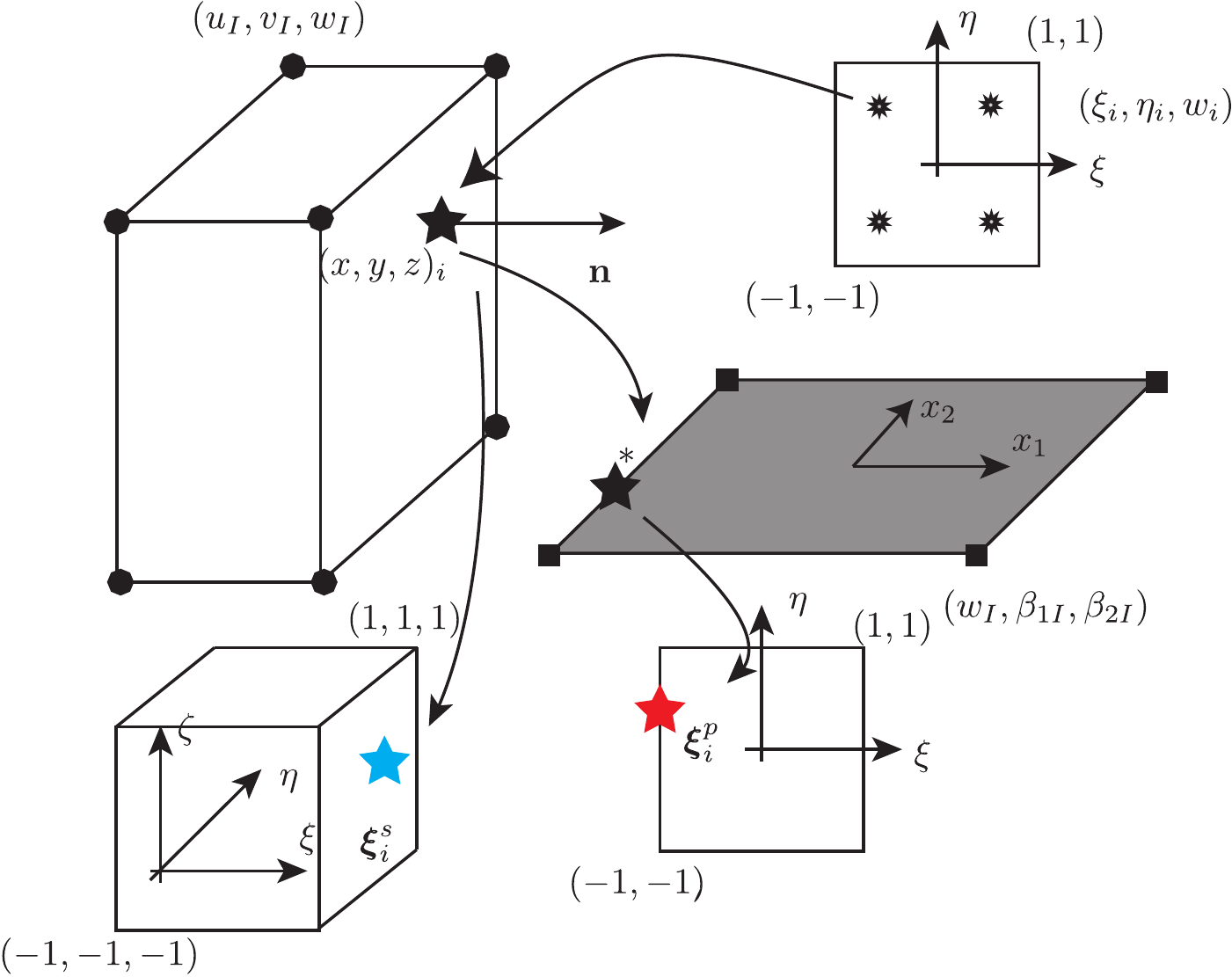}
  \caption{Coupling a plate and a 3D continuum using a Nitsche based method: determination of 
           Gauss points for evaluating the coupling terms.}
  \label{fig:plate-continuum-gauss}
\end{figure}

Care must be taken when assembling the coupling matrix given in Equation~\eqref{eq:coupling-terms}
since it depends on to which plate element $\Omega^p_{b,e}$ the GP $\bsym{\xi}^p_i$ belongs.
We refer to \cite{nguyen-nitsche1} for more details on computer implementation aspects.

\begin{rmk}
Due to the mismatch between the number of stress components in the solid (six) and the plate (three), the $6\times6$ solid constitutive matrix $\vm{C}^s$ used in the coupling matrices are reduced to a $3\times6$ matrix in which
the third, fifth and sixth rows (correspond to zero stresses in the plate theory) are removed. Accordingly, the normal matrix $\vm{n}$ is reduced as
\begin{equation}
\vm{n} = \begin{bmatrix}
n_x &  0 & n_y \\
0 & n_y & n_x \\
0 & 0 & 0 
\end{bmatrix}
\end{equation}
\end{rmk}

\begin{rmk}
The previous discussion was for the standard Lagrange finite elements which are defined in the parent space
where quadrature rules are defined. For NURBS-based isogeometric finite elements, the basis are defined in the
parameter space and thus there is a slight modification to the process. First, the GPs  
$\{\tilde{\xi}_i,\tilde{\eta}_i,w_i\}_{i=1}^{ngp}$ are  transformed to the parameter space 
$\{\xi_i,\eta_i,w_i\}_{i=1}^{ngp}$, we refer to \cite{hughes_isogeometric_2005,cottrel_book_2009} for details. 
Next, Equation~\eqref{eq:ananan} is used as usual to compute
GPs defined now in the parameter space. Note also that if B{\'e}zier extraction is used to implement NURBS-based
IGA, see \eg \cite{borden_isogeometric_2011}, then this remark can be ignored since with B{\'e}zier extraction
the basis are the Bernstein basis, which are defined in the parent space as well, multiplied 
with some sparse matrices.
\end{rmk}

\subsection{Shear deformation models}\label{shear}

The proposed formulation, previously presented, 
can be straightforwardly extended to shear deformation beam/plate theories.
In this section, a treatment of Timoshenko beam (with axial deformation included), 
and Mindlin-Reissner plate theory is given. Also presented is the treatment of the general 
continuum-beam coupling in plane frame problems. High order deformation beam/plate theories
are left out of consideration.

\subsubsection{Timoshenko beam}

The element stiffness matrix associated with a Timoshenko beam element is standard and thus not presented here.
Instead, we focus on the coupling matrices. The displacement field of a Timoshenko beam is given by 

\begin{equation}
\begin{split}
u^b_{\bar{x}} &= u(\bar{x}) -\bar{y}\theta(\bar{x})\\
u^b_{\bar{y}} &= w(\bar{x})
\end{split}
\end{equation}
where $u,w$ and $\theta$ are the axial, transverse displacement and the rotation of the beam mid-line, 
respectively.
The strain field is therefore given by

\begin{equation}
\epsilon_{\bar{x}\bar{x}}^b  = u_{,\bar{x}}-\bar{y} \theta_{,\bar{x}}, \quad
\epsilon_{\bar{y}\bar{y}}^b  = 0, \quad
2\epsilon_{\bar{x}\bar{y}}^b = -\theta + w_{,\bar{x}}
\end{equation}
The stresses are defined as
\begin{equation}
\begin{bmatrix}
\sigma_{\bar{x}\bar{x}}^b\\
\sigma_{\bar{y}\bar{y}}^b\\
\sigma_{\bar{x}\bar{y}}^b
\end{bmatrix}=
\begin{bmatrix}
E &  0 & 0\\
0 &  0 & 0\\
0 & 0 & kG  
\end{bmatrix}
\begin{bmatrix}
\epsilon_{\bar{x}\bar{x}}^b\\
\epsilon_{\bar{y}\bar{y}}^b\\
2\epsilon_{\bar{x}\bar{y}}^b
\end{bmatrix} \equiv \vm{C}^\text{b}\bsym{\epsilon}^\text{b}
\end{equation}
where $k$ denotes the shear correction factor (a value of $5/6$ is usually used) and $G$ is the shear modulus
which is defined as $G=\frac{E}{2(1+\nu)}$.

Without loss of generality, assume that the beam is discretised by two-noded
linear elements. In this case, the displacements of any beam element $e$ are given by

\begin{equation}
\begin{bmatrix}
u^b_{\bar{x}}\\
u^b_{\bar{y}}
\end{bmatrix}=\underbrace{
\begin{bmatrix}
N_1 & 0 & -y N_{1} & N_2 & 0 & - y N_{2}\\
0 & N_{1} & 0 & 0 & N_{2} & 0     
\end{bmatrix}}_{\vm{N}^b}
\begin{bmatrix}
u_1\\w_1\\\theta_1\\u_2\\w_2\\ \theta_2
\end{bmatrix}\label{eq:beam-Nb}
\end{equation}
where $(u_I,w_I,\theta_I)$ are the nodal unknowns at node $I$ of a beam element; $N_I(\xi)$ are the
univariate shape functions, and the strains are given by

\begin{equation}
\bsym{\epsilon}^b=
\begin{bmatrix}
\epsilon_{\bar{x}\bar{x}}^b\\
\epsilon_{\bar{y}\bar{y}}^b\\
2\epsilon_{\bar{x}\bar{y}}^b
\end{bmatrix}=\underbrace{
\begin{bmatrix}
N_{1,\bar{x}} & 0 & -\bar{y} N_{1,\bar{x}} & N_{2,\bar{x}} & 0 & - \bar{y} N_{2,\bar{x}}\\
   0 & 0 & 0 & 0 & 0 & 0\\
0 & N_{1,\bar{x}} & -N_1 & 0 & N_{2,\bar{x}} & - N_2     
\end{bmatrix}}_{\vm{B}^{b,c}}
\begin{bmatrix}
u_1\\w_1\\\theta_1\\u_2\\w_2\\ \theta_2
\end{bmatrix}
\end{equation}
Having obtained $\vm{N}^b$, $\vm{C}^b$ and $\vm{B}^{b,c}$, the coupling matrices given in
Equations~\eqref{eq:nitsche-kdg} and \eqref{eq:nitsche-kpe} can be computed.

\subsubsection{Plane frame analysis}\label{sec:rotation}

Firstly we introduce the rotation matrix $\vm{R}$ that relates the local nodal
unknowns to the global unknowns as

\begin{equation}
\underbrace{
\begin{bmatrix}
u_1 \\ w_1 \\ \beta_1 \\ u_2 \\ w_2 \\ \beta_2
\end{bmatrix}}_{\vm{a}^b}=\underbrace{
\begin{bmatrix}
\cos\phi & \sin \phi & 0 & 0 & 0 & 0\\
-\sin\phi & \cos\phi & 0 & 0 & 0 & 0 \\      
0 & 0 & 1 & 0 & 0 & 0 \\      
0 & 0 & 0 & \cos\phi & \sin\phi & 0 \\      
0 & 0 & 0 & -\sin\phi & \cos\phi & 0 \\      
0 & 0 & 0 & 0 & 0 & 1 \\      
\end{bmatrix}}_{\vm{R}}\underbrace{
\begin{bmatrix}
u_1 \\ w_1 \\ \beta_1 \\ u_2 \\ w_2 \\ \beta_2
\end{bmatrix}_l}_{\vm{a}^b_l}
\end{equation}
where $\phi$ is the angle between $\bar{x}$ and $x$ (cf. Fig.~\ref{fig:domain}) 
and subscript $l$ represents quantities in the local coordinate system. 

The jump operator is rewritten as follows after transforming the beam quantities
to the global system

\begin{equation}
\jump{\vm{u}} = \vm{u}^s - \vm{u}^b = \vm{u}^s - \vm{R}_v\trans\vm{u}^b_l
=\vm{N}^s\vm{a}^s -  \vm{R}_v\trans \vm{N}^b \vm{a}^b_l
=\vm{N}^s\vm{a}^s -  \vm{R}_v\trans \vm{N}^b \vm{R} \vm{a}^b
\end{equation}
where $\vm{N}^b$ is given in Equation ~\eqref{eq:beam-Nb} and 
$\vm{R}_v$ denotes the rotation matrix for vector transformation defined in Equation~\eqref{eq:Rv}.

The averaged stress operator is rewritten as 

\begin{equation}
\{\bsym{\sigma}\}=\frac{1}{2}(\bsym{\sigma}^s + \bsym{\sigma}^b)
                 =\frac{1}{2}(\vm{C}^s\vm{B}^s\vm{a}^s + \vm{T}^{-1} \vm{C}^b \vm{B}^{b,c} \vm{R} \vm{a}^b)
\end{equation}
Therefore, the coupling matrices are given by

\begin{equation}
\vm{K}^n = \begin{bmatrix}[2.5]
-\D\int_{\Gamma^*} \vm{N}^{s\text{T}} \vm{n} \frac{1}{2}\vm{C}^s\vm{B}^s \mathrm{d}\Gamma &
-\D\int_{\Gamma^*} \vm{N}^{s\text{T}} \vm{n} \frac{1}{2} \vm{T}^{-1} \vm{C}^b\vm{B}^{b,c} \vm{R}\mathrm{d}\Gamma \\
\D\int_{\Gamma^*} \vm{R}\trans \vm{N}^{b\text{T}} \vm{R}_v \vm{n} \frac{1}{2}\vm{C}^s\vm{B}^s \mathrm{d}\Gamma &
\D\int_{\Gamma^*} \vm{R}\trans \vm{N}^{b\text{T}} \vm{R}_v \vm{n} \frac{1}{2} \vm{T}^{-1}\vm{C}^b\vm{B}^{b,c}
\vm{R}\mathrm{d}\Gamma 
\end{bmatrix}
\end{equation}
and by
\begin{equation}
\vm{K}^{st} = \begin{bmatrix}[2.5]
\D\int_{\Gamma^*}  \alpha  \vm{N}^{s\text{T}} \vm{N}^s \mathrm{d}\Gamma &
 - \D\int_{\Gamma^*}  \alpha  \vm{N}^{s\text{T}} \vm{R}_v\trans \vm{N}^b \vm{R} \mathrm{d}\Gamma\\
 - \D\int_{\Gamma^*}  \alpha \vm{R}\trans \vm{N}^{b\text{T}} \vm{R}_v \vm{N}^s \mathrm{d}\Gamma&
 \D\int_{\Gamma^*}  \alpha \vm{R}\trans \vm{N}^{b\text{T}} \vm{N}^b \vm{R} \mathrm{d}\Gamma
\end{bmatrix}
\end{equation}

\subsubsection{Mindlin-Reissner plate}

The independent variables of the Mindlin-Reissner  plate theory are the rotation angle $\beta_\alpha$ ($\alpha=1,2$)
and mid-surface (transverse) displacement $w$.
The displacement field of the plate is given by

\begin{equation}
\begin{split}
u_1(x_1,x_2,x_3) &= -x_3 \beta_1\\
u_2(x_1,x_2,x_3) &= -x_3 \beta_2\\
u_3(x_1,x_2,x_3) &= w(x_1,x_2)
\end{split}\label{eq:plate-displacement}
\end{equation}
The strain field is then given by

\begin{equation}
\begin{split}
\epsilon_{11}  &= -x_3 \beta_{1,1}\\
\epsilon_{22}  &= -x_3 \beta_{2,2}\\
2\epsilon_{12} &= -x_3 (\beta_{1,2}+\beta_{2,1})\\
2\epsilon_{13} &= - \beta_{1} + w_{,1}\\
2\epsilon_{23} &= - \beta_{2} + w_{,2}\\
\end{split}\label{eq:plate-strain}
\end{equation}

Finite element approximations of the displacement and rotations are given by

\begin{equation}
w = N_I(\xi,\eta) w_I,\quad
\beta_\alpha = N_I(\xi,\eta) \beta_{\alpha I}
\end{equation}

\noindent where $N_I(\xi,\eta)$ is the shape function associated with node I
and $w_I,\beta_{1I},\beta_{2I}$ denote the nodal unknowns that include the nodal deflection and two rotations.

The element plate stiffness matrix is now defined as, see \eg \cite{taylor-fem-book}

\begin{equation}
\vm{K}_e^p = \int_{\Omega_e^p} \left[ (\vm{B}_{be}^p)\trans \vm{D}^p_b \vm{B}_{be}^p +
(\vm{B}_{se}^p)\trans \vm{D}^p_s \vm{B}_{se}^p 
\right] \di \Omega
\end{equation}

\noindent and the element displacement-curvature matrix $\vm{B}_{be}$ and
displacement-shearing matrix are given by

\begin{equation}
\vm{B}_{be}^p = \begin{bmatrix}
\vm{B}_{b1} & \vm{B}_{b2} & \cdots & \vm{B}_{bn} 
\end{bmatrix},\;\;
\vm{B}_{se}^p = \begin{bmatrix}
\vm{B}_{s1} & \vm{B}_{s2} & \cdots & \vm{B}_{sn} 
\end{bmatrix}
\end{equation}
where components are given by
\begin{equation}
\vm{B}_{bI}=
\begin{bmatrix}
0 & N_{I,1} & 0\\
0 & 0 & N_{I,2}\\
0 & N_{I,2} & N_{I,1}
\end{bmatrix},\quad
\vm{B}_{sI}=
\begin{bmatrix}
N_{I,1} & -N_I & 0 \\
N_{I,2} & 0 & -N_I\\
\end{bmatrix}
\label{plate-curvature1}
\end{equation}
Matrices $\vm{D}^p_b$ and $\vm{D}^p_s$ denote the bending constitutive matrix
and shear constitutive matrix, respectively. The former was given in Equation~\eqref{eq:bending-stiffness}
and the latter is given by

\begin{equation}
\vm{D}^p_s=\frac{kEh}{2(1+\nu)}\begin{bmatrix}
1 & 0\\
0 & 1
\end{bmatrix}
\end{equation}
where $k$ is the shear correction factor.

Next, we focus on the coupling matrices. Firstly 
the displacement field given in Equation~\eqref{eq:plate-displacement} can be written as

\begin{equation}
\vm{u}^\text{plate} = \begin{bmatrix}
0 & -x_3N_I & 0\\
0 & 0 & -x_3 N_I\\
N_I & 0 & 0 
\end{bmatrix}
\begin{bmatrix}
w\\ \beta_1 \\ \beta_2
\end{bmatrix}_I=\vm{N}_I^p\vm{a}^p_I
\end{equation}
and the strain field given in Equation~\eqref{eq:plate-strain} as

\begin{equation}
\begin{bmatrix}
\epsilon_{11}\\
\epsilon_{22}\\
2\epsilon_{12}\\
2\epsilon_{23}\\
2\epsilon_{13}\\
\end{bmatrix}=
\begin{bmatrix}
0 & -x_3 N_{I,1} & 0 \\
0 & 0 & -x_3 N_{I,2}\\
0 & -x_3N_{I,2} & -x_3 N_{I,1}\\  
N_{I,2} & 0 & -N_I \\
N_{I,1} & -N_I & 0      
\end{bmatrix}
\begin{bmatrix}
w_I\\ \beta_{1I}\\ \beta_{2I}
\end{bmatrix}\equiv\vm{B}_{I}^{p,c}  \vm{a}_I^p
\label{plate-curvature1}
\end{equation}
And finally the stress field is given by
\begin{equation}
\begin{bmatrix}
\sigma_{11}\\
\sigma_{22}\\
\sigma_{12}\\
\sigma_{23}\\
\sigma_{13}\\
\end{bmatrix}=
\begin{bmatrix}
\frac{E}{(1-\nu^2)}\begin{bmatrix}
1 & \nu & 0\\
\nu & 1 & 0\\
0 & 0 & 0.5(1-\nu)
\end{bmatrix} & 
\begin{bmatrix}
0 & 0 \\
0 & 0 \\
0 & 0 
\end{bmatrix}\\
\begin{bmatrix}
0 & 0 & 0\\
0 & 0 & 0
\end{bmatrix} &
\frac{k E}{2(1+\nu)}\begin{bmatrix}
1 & 0\\
0 & 1
\end{bmatrix}
\end{bmatrix}
\begin{bmatrix}
\epsilon_{11}\\
\epsilon_{22}\\
2\epsilon_{12}\\
2\epsilon_{23}\\
2\epsilon_{13}\\
\end{bmatrix}\equiv \vm{C}^p\bsym{\epsilon}^p
\label{plate-curvature1}
\end{equation}
Having obtained $\vm{N}^p$, $\vm{C}^p$ and $\vm{B}^{p,c}$, the coupling matrices 
can thus be computed.

\section{Non-conforming coupling}\label{non-conforming}

In this section we present a non-conforming coupling formulation.
In the literature when applied for mesh coupling the method was referred to as composite grids
\cite{MZA:8203286} or embedded mesh method \cite{Sanders2011a}. Due to the great similarity between
solid/beam and solid/plate coupling, only the former is discussed, cf. Fig.~\ref{fig:volume-coupling1}.

\begin{figure}[htbp]
  \centering 
  \includegraphics[width=0.5\textwidth]{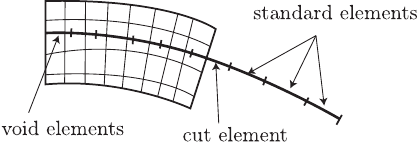}
  \caption{Non-conforming coupling of a solid and a beam.}
  \label{fig:volume-coupling1}
\end{figure}

The weak form and thus the discrete equations are exactly the same (see Section \ref{sec:solid-beam}).
There are however changes in the treatment of beam elements of which some are overlapped by solid elements.
Using the terminology from \cite{Sanders2011a} (which was actually adopted from extended finite element
method--XFEM \cite{Sukumar1}) we divide the beam elements into three groups namely (1) 
   standard elements which do not
intersect with solid elements, (2) cut elements (cut by the boundaries of the embedded solid) and (3) void elements which are completely covered by the solid. Void elements do not contribute to the total stiffness of the system.
Therefore, nodes whose support is completely covered by the solid are considered to be inactive and removed
from the stiffness matrix. Another way that preserves the dimension of the stiffness matrix is to fix those inactive beam nodes.

\begin{figure}[htbp]
  \centering 
  \includegraphics[width=0.5\textwidth]{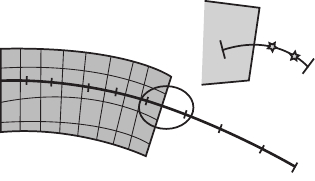}
  \caption{Non-conforming coupling of a solid and a beam: integration of cut element.}
  \label{fig:volume-coupling2}
\end{figure}

For numerical integration of cut elements, the sub-triangulation (for solid/plate problem) as usually used in
the context of XFEM \cite{mos_finite_1999}  can be adopted. However, in the provided examples, we use a simpler
(easy for implementation) technique: a large number of Gauss points is used for cut elements and those points which are in the solid domain are skipped.

For some special configurations as given in Fig.~\ref{fig:volume-coupling3} in which the cut (beam) element
is almost completely covered by the continuum element, care must be taken to avoid a badly scaled singular
stiffness matrix. The remedy is however simple: the node of the cut element whose support is nearly
in the continuum mesh is marked inactive. Practically it is the node/control point that is farthest from the coupling interface.

\begin{figure}[htbp]
  \centering 
  \includegraphics[width=0.5\textwidth]{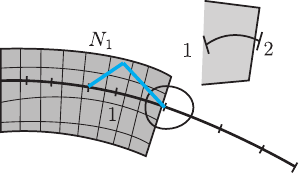}
  \caption{Non-conforming coupling of a solid and a beam: cut element almost falls within the continuum element.
  In order to avoid a badly scaled stiffness matrix, the node of the cut element whose support is nearly
  in the continuum mesh is marked inactive. In the figure, it is node 1 that is inactive.}
  \label{fig:volume-coupling3}
\end{figure}

\begin{rmk}
One might ask why the void elements are necessary although they do not contribute to the global stiffness matrix.
The answer is that they are needed in a model adaptive formulation where at a certain time step, they are
void elements but in other time steps they become active when the refined solid is no longer needed there.
\end{rmk}

\section{Numerical examples}\label{sec:examples}

In accordance to the two topics presented previously, in this section, numerical examples
on beam/continuum and plate/continuum coupling are presented to
demonstrate the performance of the proposed formulation. Specially, the following examples
are provided together with their raison d'etre.

\begin{enumerate}
\item Beam-continuum coupling
   \begin{enumerate}
   \item Cantilever beam: conforming coupling
   \item Cantilever beam: non-conforming coupling
   \item Frame analysis
   \end{enumerate}
\item Plate-continuum coupling
   \begin{enumerate}
   \item Cantilever beam (3D/plate coupling)
   \item Cantilever beam (non-conforming coupling)
   \item Non-conforming coupling of a square plate
   \end{enumerate}
\end{enumerate}

Unless otherwise stated, we use 
MIGFEM--an open source Matlab IGA code\footnote{which is available at \url{https://sourceforge.net/projects/cmcodes/}} for our computations and  the Timoshenko beam theory and the Reissner-Mindlin plate theory.
The visualisation was performed in Paraview \cite{para}.
In all examples domains where reduced models such as beams and plates
are utilized are assumed to be known \textit{a priori}. 
A rule of $(p+1)\times(q+1)$ Gaussian quadrature can be applied for
two-dimensional elements in which $p$ and $q$ denote the orders of
the chosen basis functions in the $\xi$ and $\eta$ direction. The same
procedure is also used for NURBS basis functions in the present work,
although it should be emphasised that Gaussian quadrature is not optimal for IGA.
Research is currently focussed on optimal integration techniques such
as that in \cite{hughes_efficient_2010,Auricchio201215} in which an
optimal quadrature rule, known as the half-point rule, has been applied.

\subsection{Continuum-beam coupling}

\subsubsection{Timoshenko beam: conforming coupling}

Consider a beam of dimensions $L \times D$, subjected to a
parabolic traction at the free end as shown in Fig.~\ref{fig:beam-geo}.
The beam is considered to be of unit depth and is in plane stress
state. The parabolic traction is given by

\begin{equation}
    t_y(y) = -\frac{P}{2I} \biggl ( \frac{D^2}{4} - y^2 \biggr)\label{eq:ty}
\end{equation}

\noindent where $I = D^3/12$ is the moment of inertia. The exact displacement
field of this problem is \cite{elasticity_book}

\begin{equation}
\begin{split}
    u_x(x,y) &=  \frac{Py}{6EI} \biggl [ (6L-3x)x + (2+\nu)\biggl(y^2-\frac{D^2}{4}\biggr) \biggr] \\
    u_y(x,y) &=  - \frac{P}{6EI} \biggl [ 3\nu y^2(L-x) + (4+5\nu)\frac{D^2x}{4} +(3L-x)x^2 \biggr] \\
\end{split}
\label{eq:tBeamExactDisp}
\end{equation}

\noindent and the exact stresses are

\begin{equation}
	\sigma_{xx}(x,y) =  \frac{P(L-x)y}{I}; 
	\quad \sigma_{yy}(x,y) =  0, \quad
\sigma_{xy}(x,y) =  -\frac{P}{2I} \biggl ( \frac{D^2}{4}-y^2\biggr)
\end{equation}

\noindent In the computations, material properties are taken as $E=
3.0 \times 10^7$, $\nu = 0.3$ and the beam dimensions are $D=6$ and
$L=48$. The shear force is $P = 1000$.  
Units are deliberately left out here, given that they can be consistently chosen in any system.
In order to model the clamping condition, 
the displacement defined by Equation~\eqref{eq:tBeamExactDisp} is prescribed as essential boundary 
conditions at $x=0, -D/2 \le y \le D/2$. This problem is solved with bilinear Lagrange elements (Q4 elements)
and high order B-splines elements. The former helps to verify the implementation in addition to
the ease of enforcement of Dirichlet boundary conditions (BCs). For the latter, care must be taken
in enforcing the Dirichlet BCs given in Equation~\eqref{eq:tBeamExactDisp} since the B-spline basis functions
are not interpolatory. 

\begin{figure}[htbp]
  \centering 
  \psfrag{p}{P}\psfrag{l}{$L$}\psfrag{d}{$D$}\psfrag{x}{$x$}\psfrag{y}{$y$}
  \includegraphics[width=0.5\textwidth]{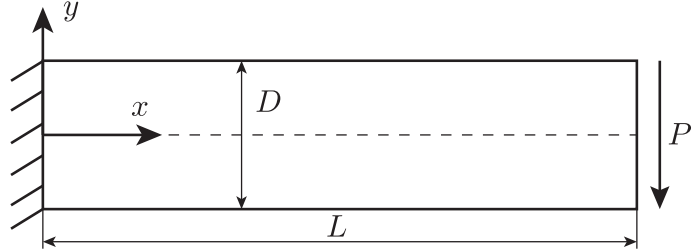}
  \caption{Timoshenko beam: problem description.}
  \label{fig:beam-geo}
\end{figure}

The mixed continuum-beam model is given in Fig.~\ref{fig:beam-continuum-example}.
The end shear force applied to the right end point is $F=P$.

\begin{figure}[htbp]
  \centering 
   \includegraphics[width=0.55\textwidth]{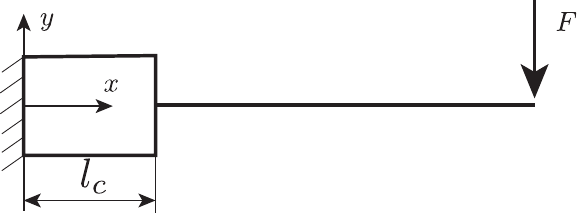}
  \caption{Timoshenko beam:  mixed continuum-beam model.}
  \label{fig:beam-continuum-example}
\end{figure}

\noindent \textbf{Lagrange elements} In the first calculation we take $l_c=L/2$ and a mesh of $40\times10$ Q4 elements (40 elements in the length direction)  was used for the continuum part
and 29 two-noded elements for the beam part. The stabilisation parameter $\alpha$ according to
Equation~\eqref{eq:alpha-eigen} was $4.7128\times10^7$.
Fig.~\ref{fig:beam-continuum-res}a plots the transverse displacement (taken as nodal values) along the beam length at $y=0$ together with the exact solution given in Equation~\eqref{eq:tBeamExactDisp}. An excellent agreement 
with the exact solution can be observed and this verified the implementation. The comparison of the numerical stress field and the exact stress field is given in Fig.~\ref{fig:beam-continuum-res}b with less satisfaction. While the bending stress $\sigma_{xx}$ is well estimated, the shear stress $\sigma_{xy}$ is not well predicted in proximity to the coupling interface. This phenomenon was also observed in the framework of Arlequin method
\cite{Hu2008a} and in the context of MPC method \cite{Gabbert}. Explanation of this phenomenon will be
given subsequently.\\

\begin{figure}[htbp]
  \centering 
   \subfloat[transverse displacement]{\includegraphics[width=0.5\textwidth]{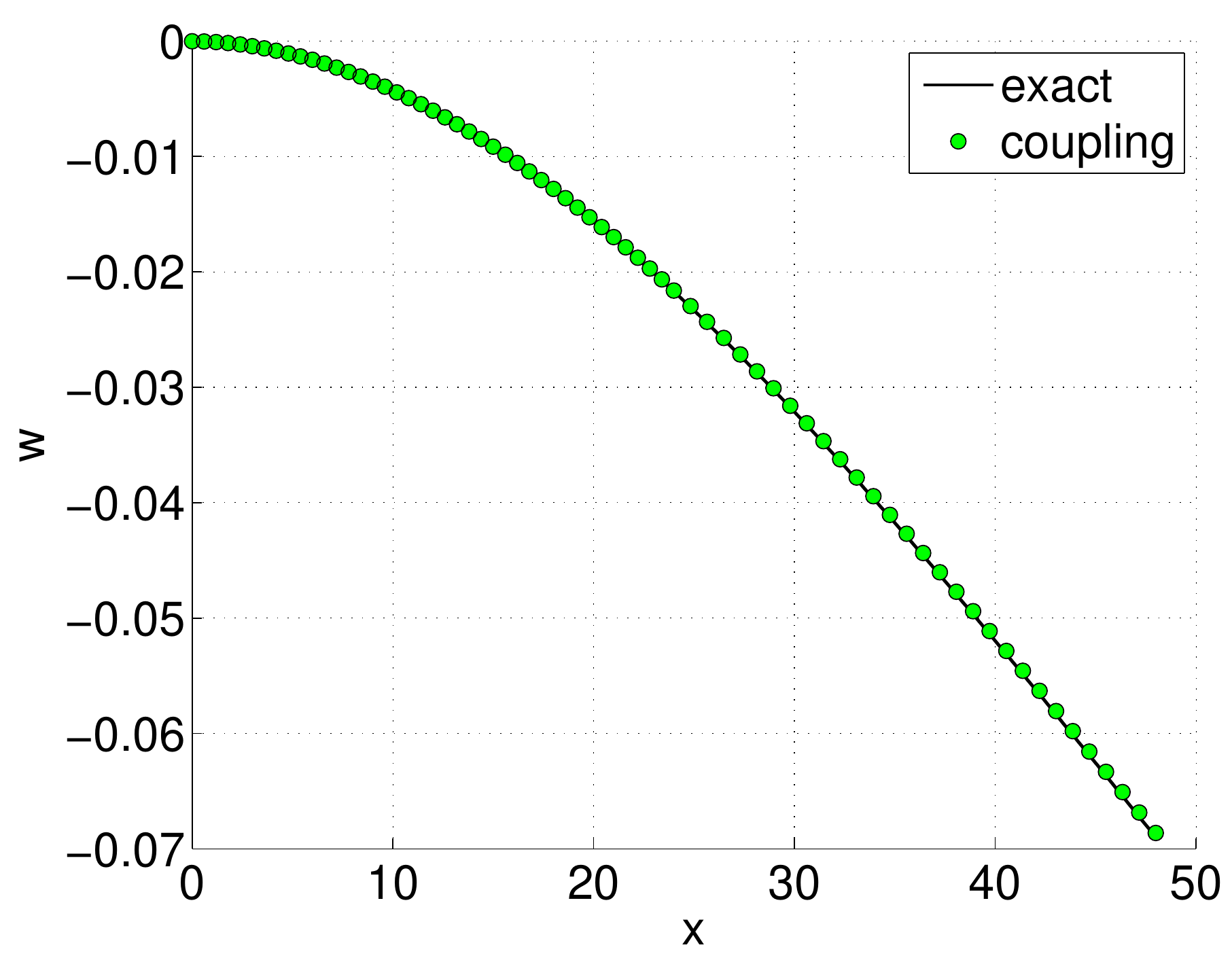}}
  \subfloat[stresses]{\includegraphics[width=0.5\textwidth]{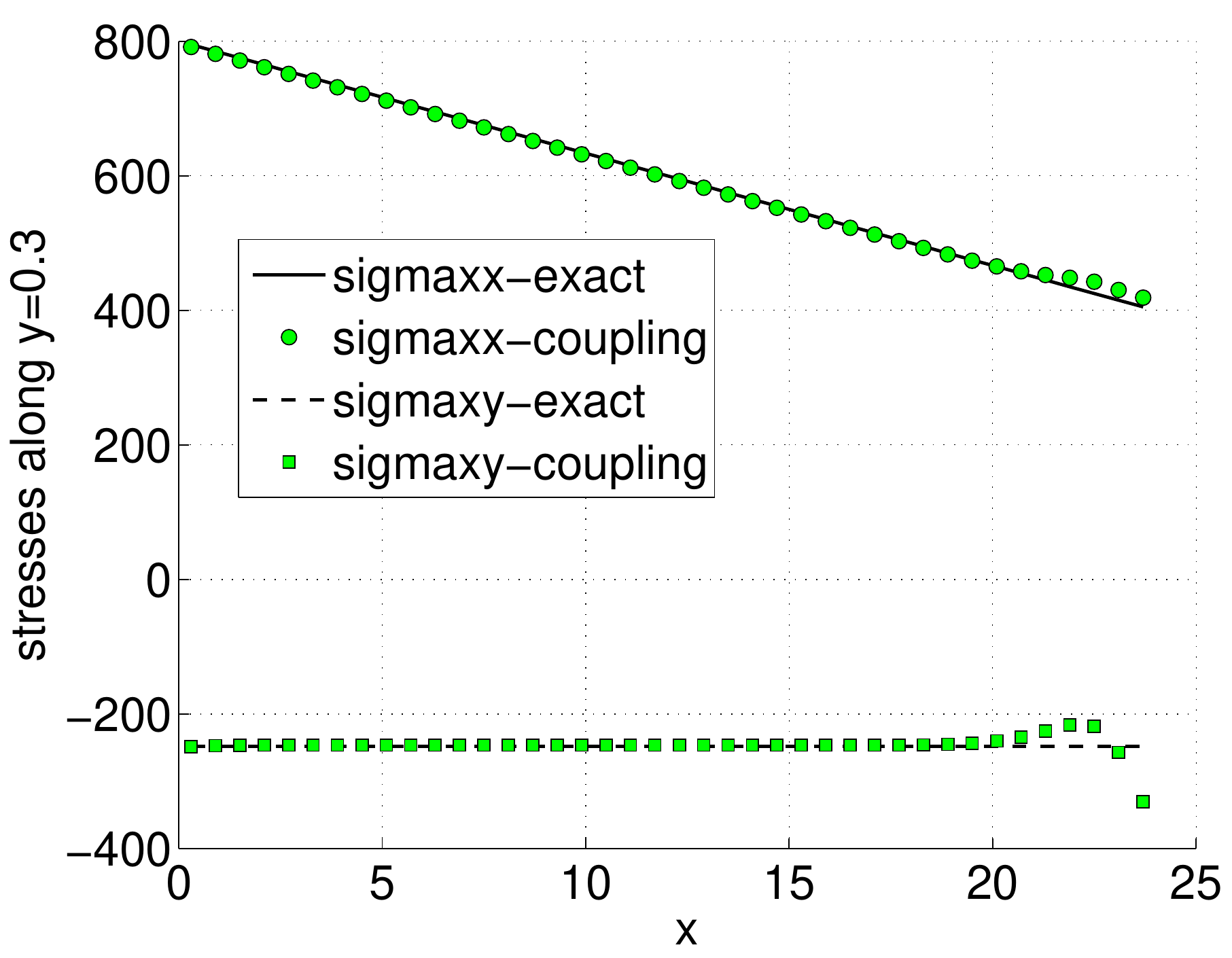}}
  \caption{Mixed dimensional analysis of the Timoshenko beam: comparison of numerical solution
     and exact solution.}
  \label{fig:beam-continuum-res}
\end{figure}

\noindent \textbf{B-splines elements} are used to discretise the continuum part (with bi-variate B-splines
      elements) and the beam part (with uni-variate B-splines elements). Such a mesh is given in 
  Fig.~\ref{fig:beam-2D-spline-mesh}. 
Dirichlet BCs are enforced using the least square projection method see \eg \cite{nguyen_iga_review}.
Note that Nitche's method can also be used to weakly enforce the Dirichlet BCs. However, we use
Nitsche's method only to couple the different models. In what follows, we used the following discretisation--
$16\times4$ bi-cubic continuum elements  and 4 cubic beam elements. 
The stabilisation parameter $\alpha$ according to Equation~\eqref{eq:alpha-eigen} was $5.5\times10^9$.
Comparison between numerical and exact solutions are given in Fig.~\ref{fig:beam-spline-res}. Again, an excellent
estimation of the displacement was obtained whereas the shear stress is not well captured. An explanation for this
behavior is given in Fig.~\ref{fig:beam-2D-explain}. The error of the continuum-beam model consists of two parts:
(1) model error when one replaces a continuum model by a continuum-beam model and (2) discretisation and coupling
errors. When the former is dominant, the coupling method is irrelevant as the same phenomenon was observed 
in the Arlequin method, in the MPC method \cite{Hu2008a,Gabbert} and mesh refinement does not cure the problem. 
In other words, the error originating from the Nitsche coupling is dominated by the model error 
and therefore we cannot draw any conclusions about the "coupling error".
In order to alleviate  the model error,
we computed the shear stress with $\nu=0.0$ and the results are given in Figs.~\ref{fig:beam-spline-poisson}
and \ref{fig:beam-spline-poisson1}. With $\nu=0.0$, the numerical solution is very much better than the one with
$\nu=0.3$.

\begin{figure}[htbp]
  \centering 
   \includegraphics[width=0.55\textwidth]{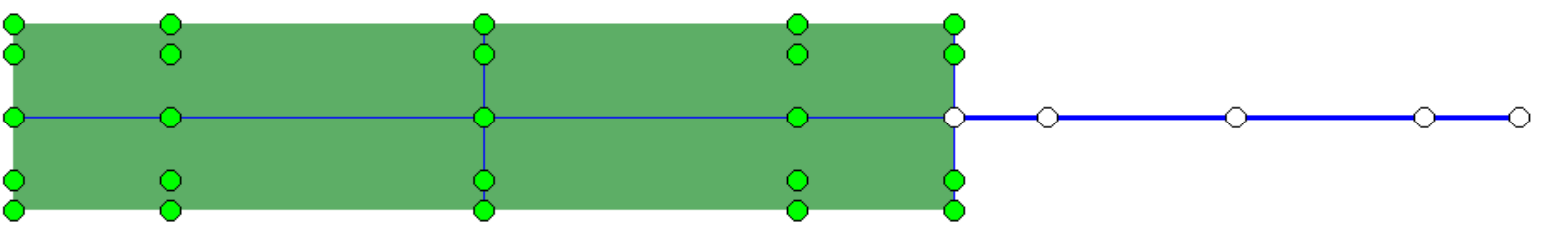}
  \caption{Mixed dimensional analysis of the Timoshenko beam: discretisation with B-splines elements.
  The continuum part is meshed by $2\times2$ cubic B-splines and the beam part is with 2 cubic elements.}
  \label{fig:beam-2D-spline-mesh}
\end{figure}

\begin{figure}[htbp]
  \centering 
  \subfloat[transverse displacement]{\includegraphics[width=0.5\textwidth]{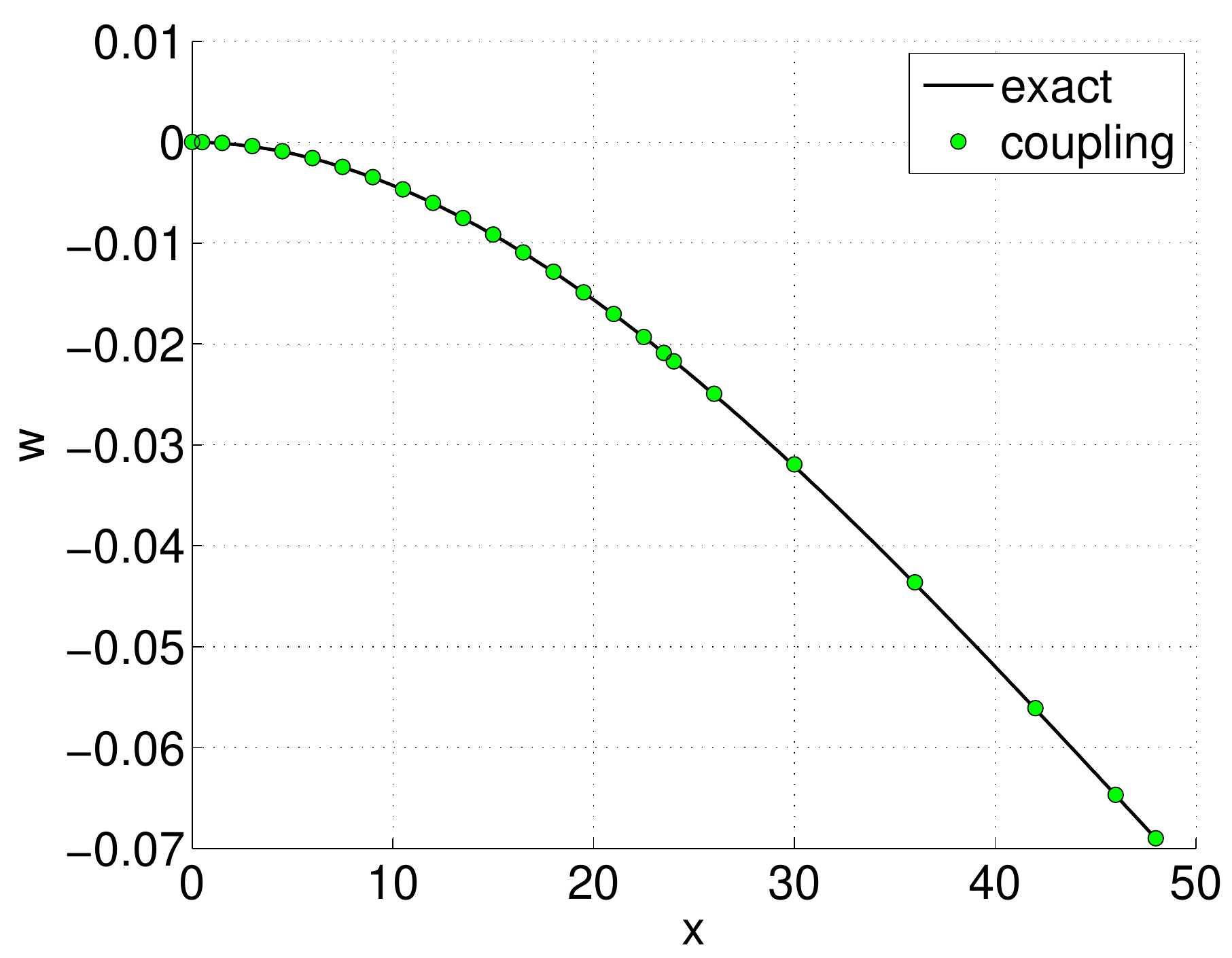}}
  \subfloat[stresses]{\includegraphics[width=0.5\textwidth]{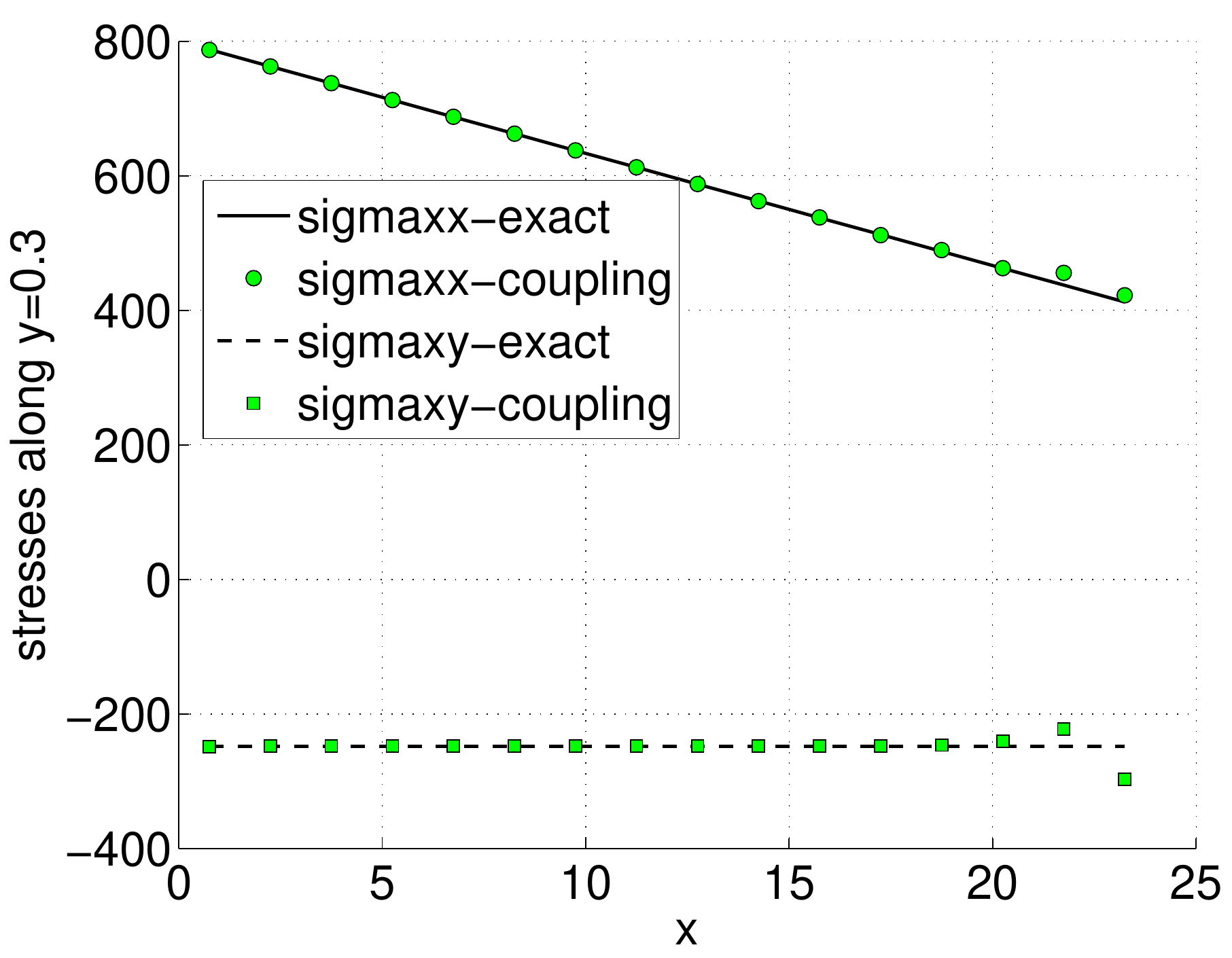}}
  \caption{Mixed dimensional analysis of the Timoshenko beam with B-spline elements: 
     comparison of numerical solution and exact solution.}
  \label{fig:beam-spline-res}
\end{figure}

\begin{figure}[htbp]
  \centering 
   \includegraphics[width=0.55\textwidth]{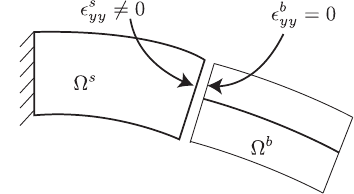}
  \caption{Continuum-beam model: inextensibility assumption employed in the beam model is not consistent with
  the continuum strain state. This introduced an error when one replaces a continuum model by a continuum-beam
  model.}
  \label{fig:beam-2D-explain}
\end{figure}

\begin{figure}[htbp]
  \centering 
  \subfloat[$\nu=0.3$]{\includegraphics[width=0.5\textwidth]{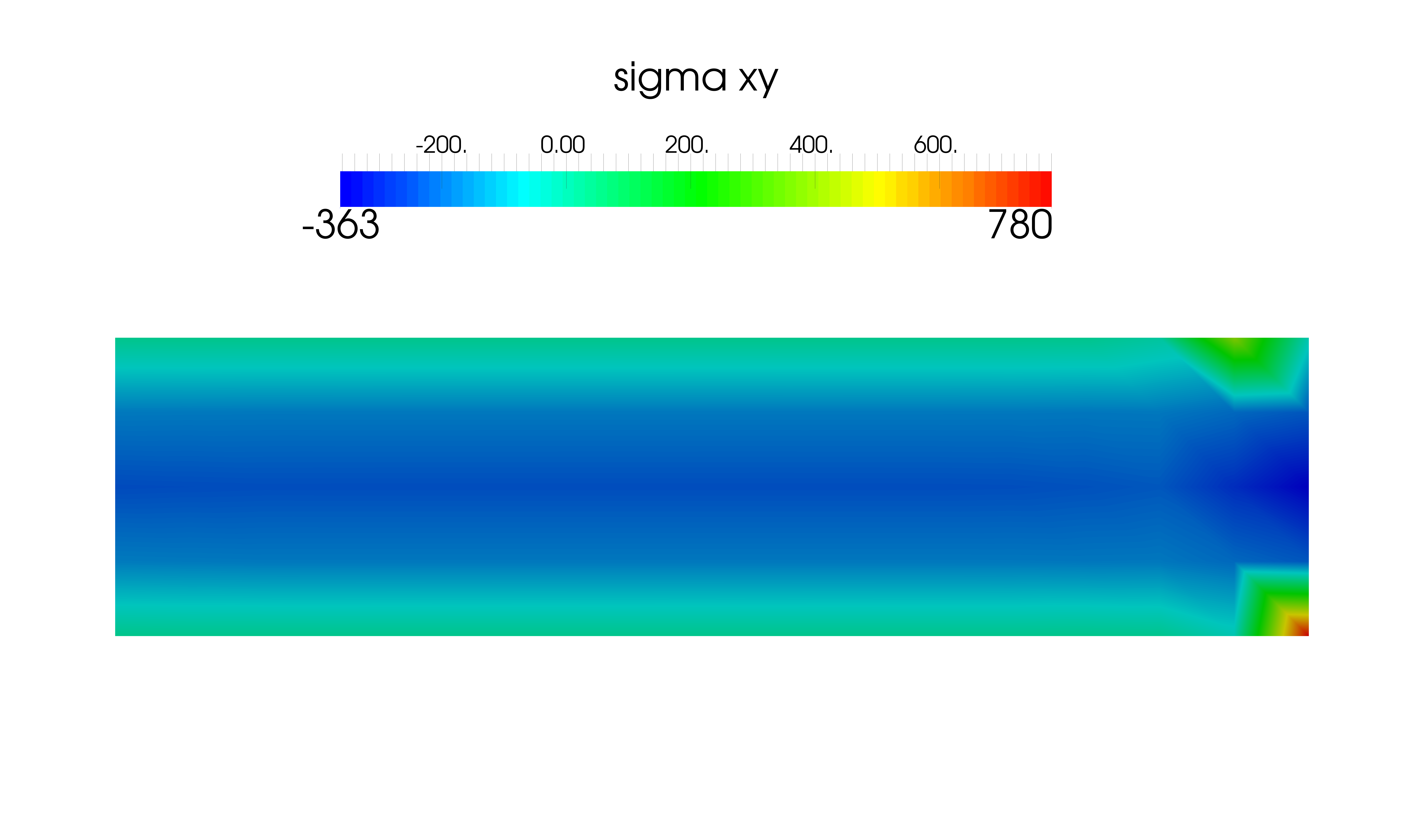}}
  \subfloat[$\nu=0.0$]{\includegraphics[width=0.5\textwidth]{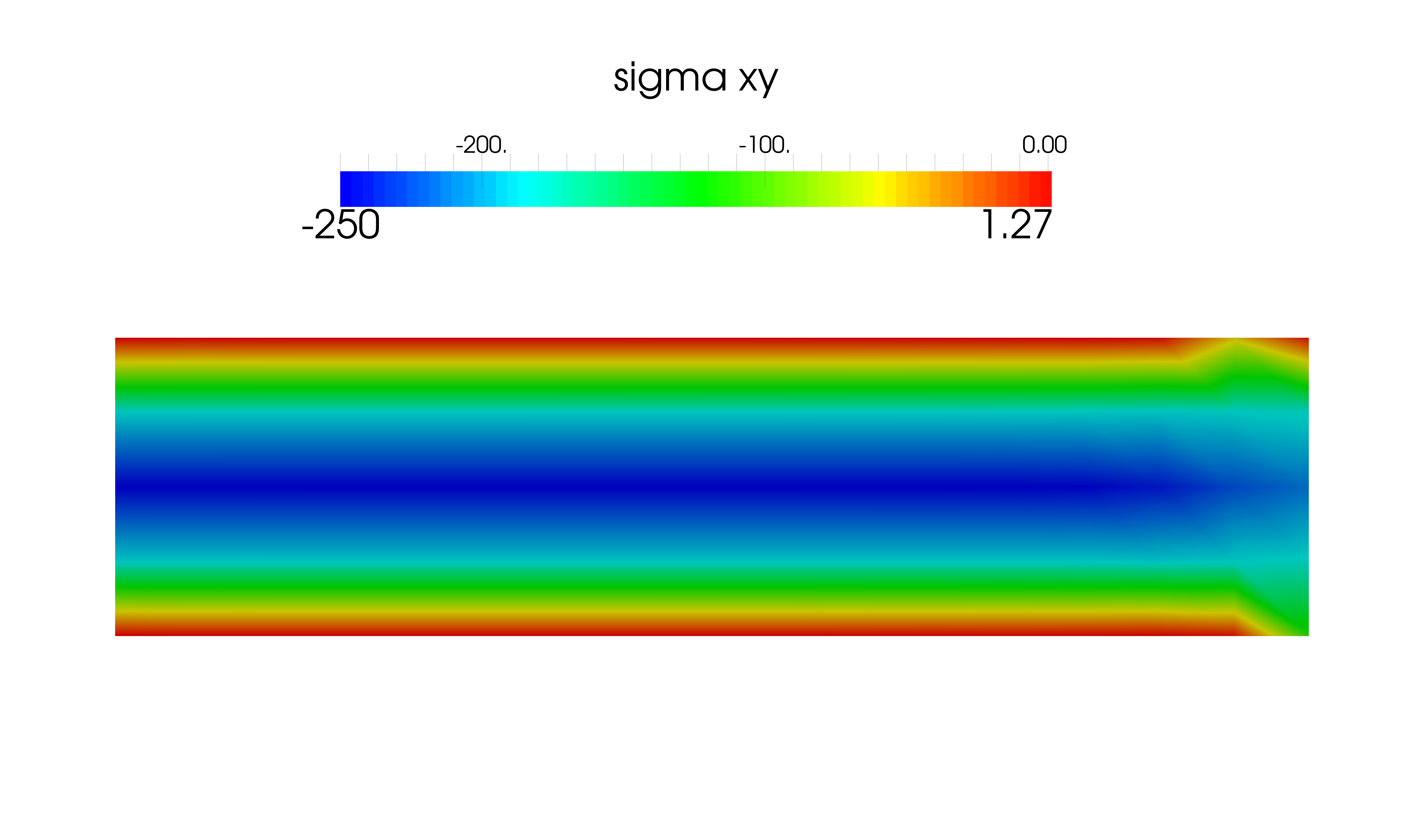}}
  \caption{Mixed dimensional analysis of the Timoshenko beam with B-spline elements: numerical shear stresses
  with $\nu=0.3$ and $\nu=0.0$. } 
  \label{fig:beam-spline-poisson}
\end{figure}

\begin{figure}[htbp]
  \centering 
  \subfloat[$\nu=0.3$]{\includegraphics[width=0.5\textwidth]{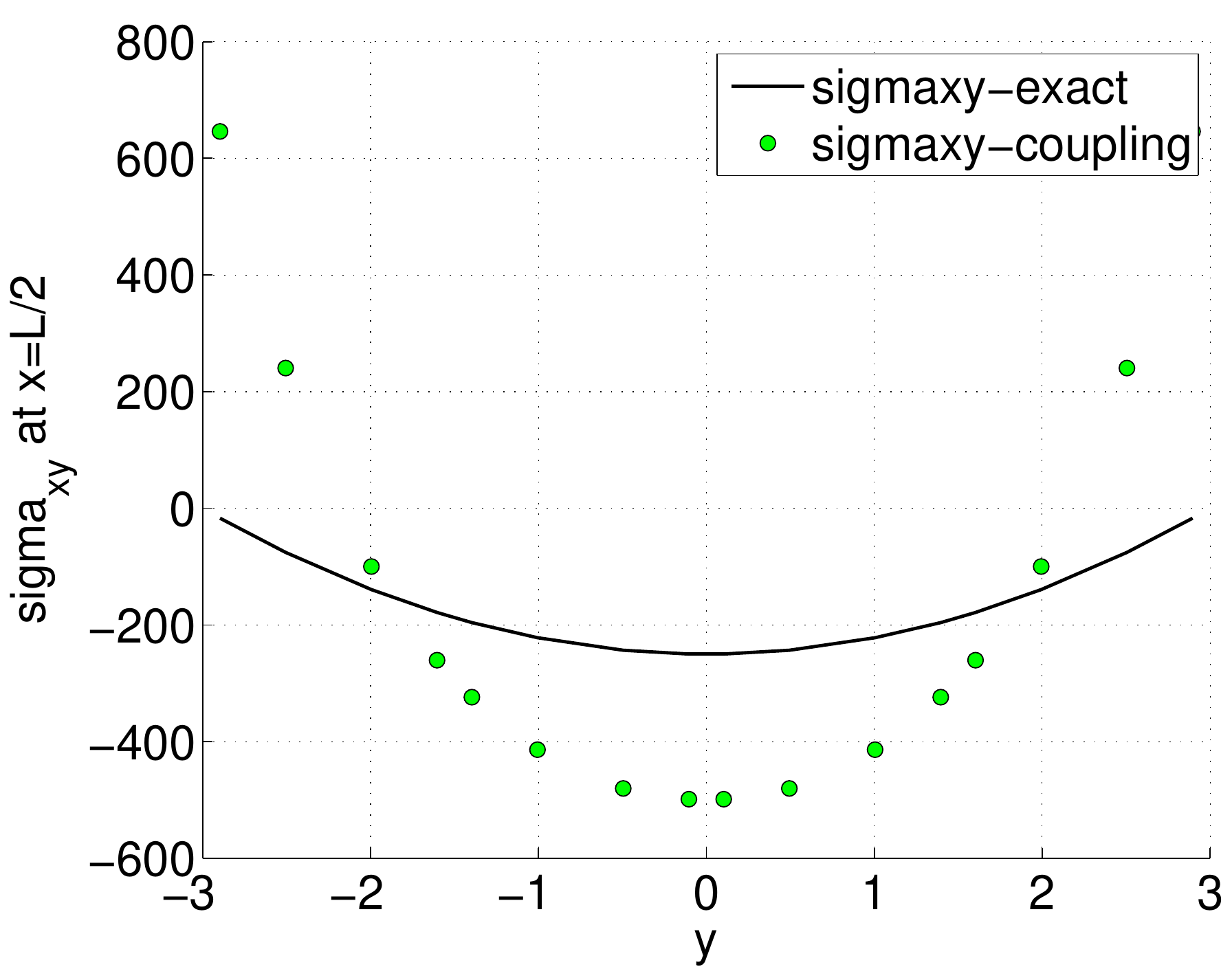}}
  \subfloat[$\nu=0.0$]{\includegraphics[width=0.5\textwidth]{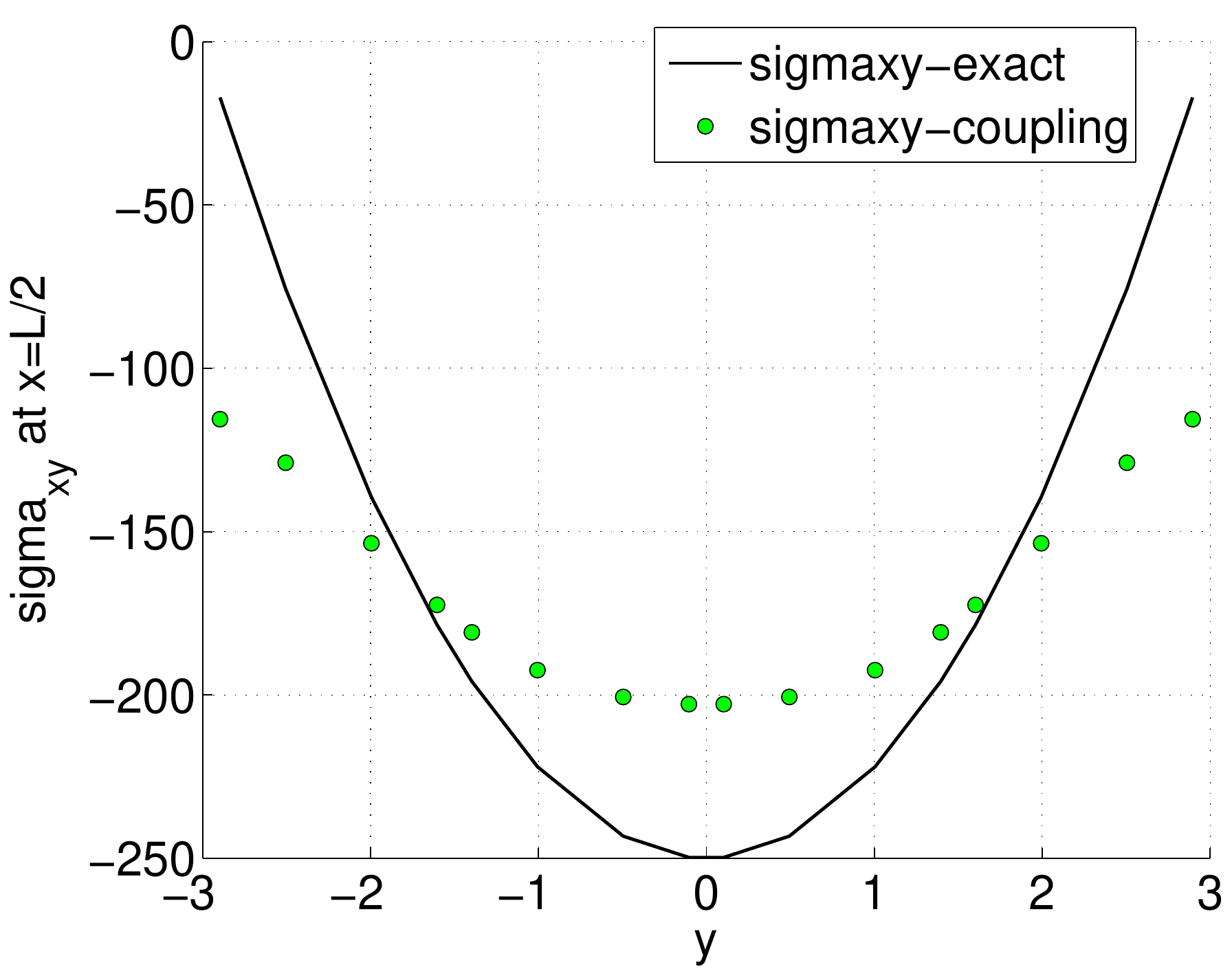}}
  \caption{Mixed dimensional analysis of the Timoshenko beam with B-spline elements: numerical shear stresses
  along the coupling interface with $\nu=0.3$ and $\nu=0.0$.
  The distribution of the shear stress with $\nu=0.3$ is very similar to the result 
  presented in \cite{Gabbert} that uses the MPC (cf. Figure 18 of the referred paper).} 
  \label{fig:beam-spline-poisson1}
\end{figure}

\subsubsection{Timoshenko beam: non-conforming coupling}

In  this section, a non-conforming coupling is considered. The B-spline mesh is given in Fig.~\ref{fig:beam-2D-nonconform-mesh}. Refined meshes are obtained from this one via the knot span subdivision technique. We use
the mesh consisting of $32\times4$ cubic continuum elements and 8 cubic beam elements.
Fig.~\ref{fig:beam-2D-nonconform-disp} gives the mesh and the displacement field in which $l_c=29.97$
so that the coupling interface is very close to the beam element boundary. A good solution was obtained
using the simple technique described in Section \ref{non-conforming}.

\begin{figure}[htbp]
  \centering 
   \includegraphics[width=0.55\textwidth]{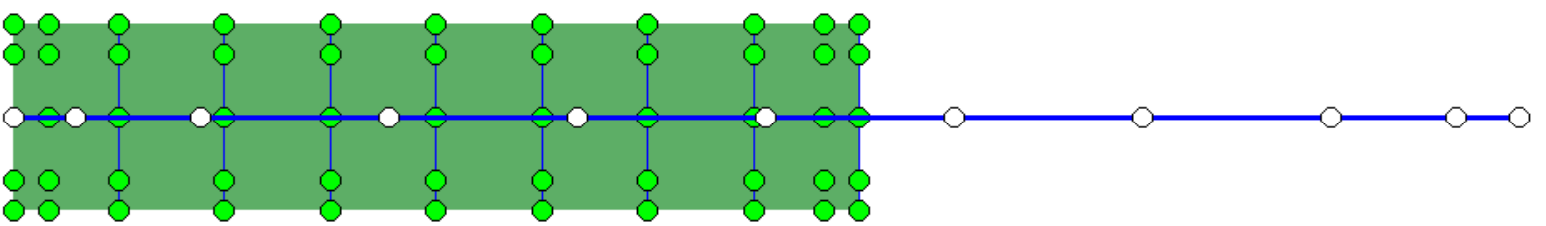}
  \caption{Mixed dimensional analysis of the Timoshenko beam with non-conforming coupling.
  The continuum part is meshed by $8\times2$ bi-cubic B-splines and the beam part is with 8 cubic elements.}
  \label{fig:beam-2D-nonconform-mesh}
\end{figure}

\begin{figure}[htbp]
  \centering 
  \subfloat[mesh]{\includegraphics[width=0.5\textwidth]{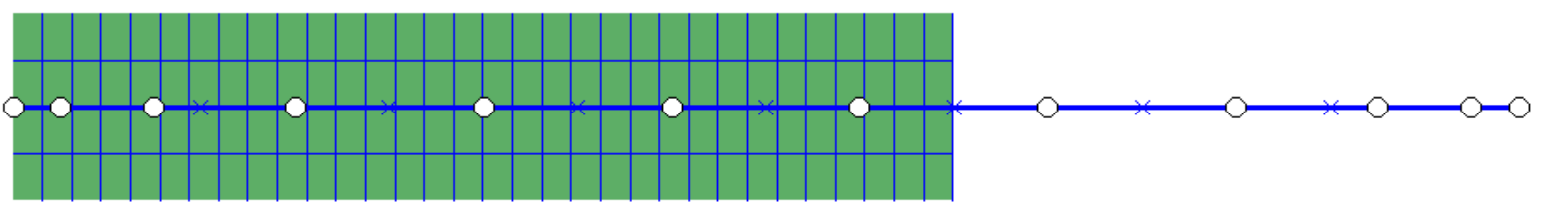}}
  \subfloat[displacement field]{\includegraphics[width=0.5\textwidth]{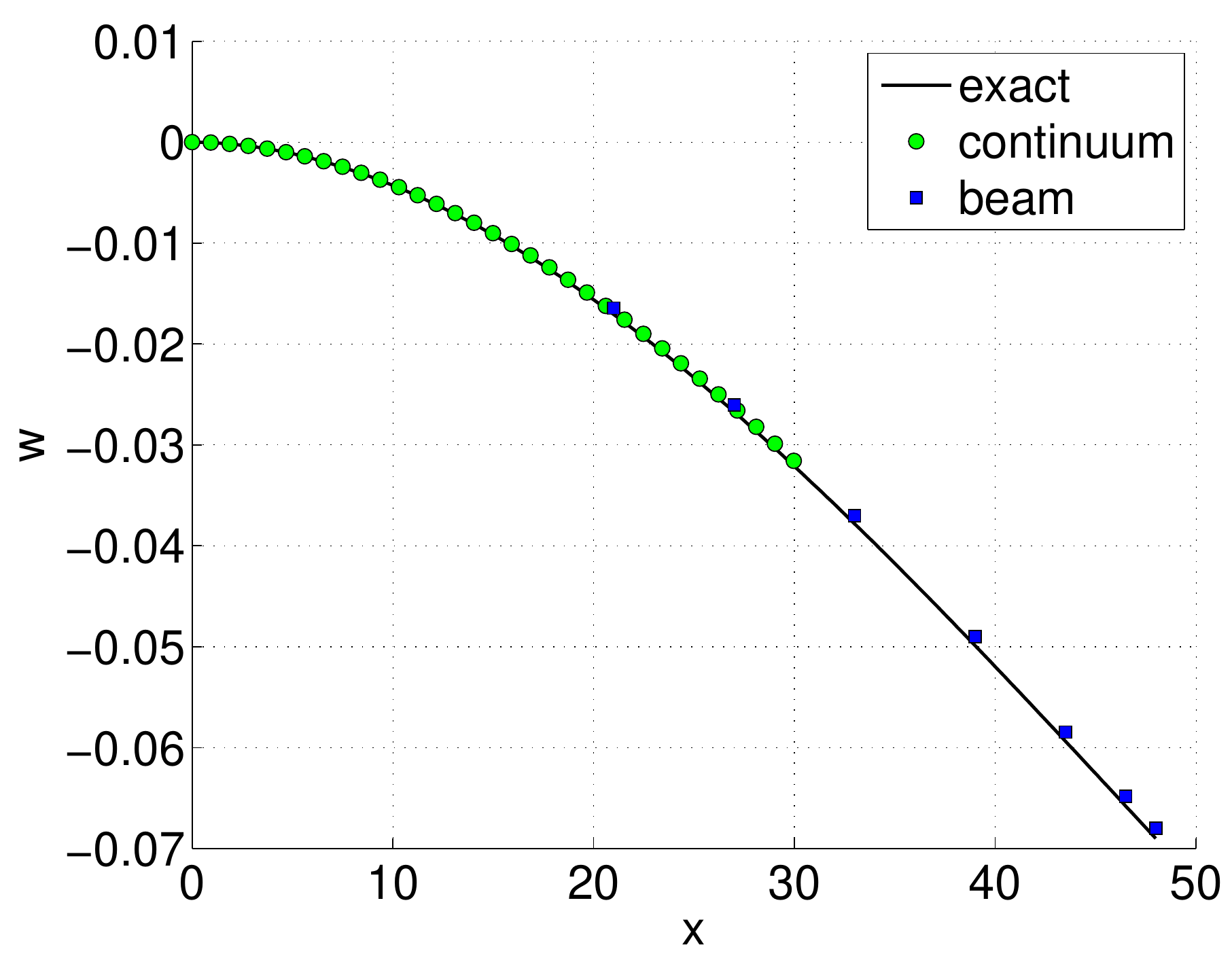}}
  \caption{Mixed dimensional analysis of the Timoshenko beam with non-conforming coupling: 
     (a) $32\times4$ Q4 elements and 8 quartic ($p=4$) beam elements and (b) displacement field.}
  \label{fig:beam-2D-nonconform-disp}
\end{figure}

\subsubsection{Frame analysis}

In order to demonstrate the correctness of the solid-beam coupling in which the beam local coordinate
system is not identical to the global one, we perform a plane frame analysis as shown in 
Fig.~\ref{fig:frame-problem}.
Due to symmetry, only half model is analysed with appropriate symmetric boundary conditions. We solve this model
with (1) continuum model (discretised with 7105 four-noded quadrilateral elements, 7380 nodes, 14760 dofs, Gmsh
\cite{geuzaine2009gmsh} was used) 
and (2) solid-beam model (cf. Fig.~\ref{fig:frame-geo}). 
The beam part are discretised using two-noded frame elements with three
degrees of freedom (dofs) per node (axial displacement, transverse displacement and rotation). 
Note that continuum element nodes have only two dofs. The total number of dofs of the continuum-beam model
is only 5400.
The stabilisation parameter is taken to be $\alpha=10^7$ and used for both coupling interfaces.
A comparison of $\sigma_{xy}$ contour plot obtained with (1) and (2) is given in
Fig.~\ref{fig:frame-res}. A good agreement was obtained.

\begin{figure}[htbp]
  \centering 
   \includegraphics[width=0.4\textwidth]{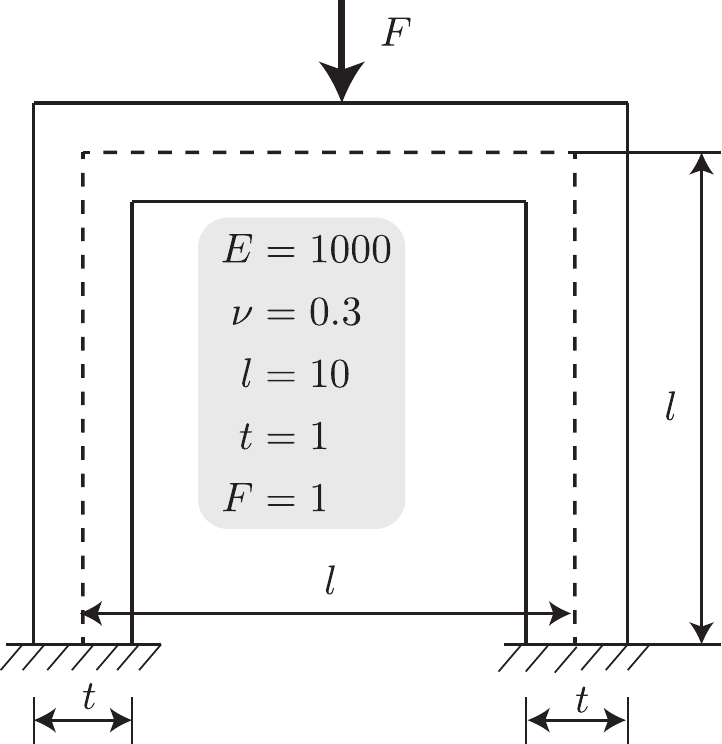}
  \caption{A plane frame analysis: problem description.}
  \label{fig:frame-problem}
\end{figure}

\begin{figure}[htbp]
  \centering 
   \includegraphics[width=0.4\textwidth]{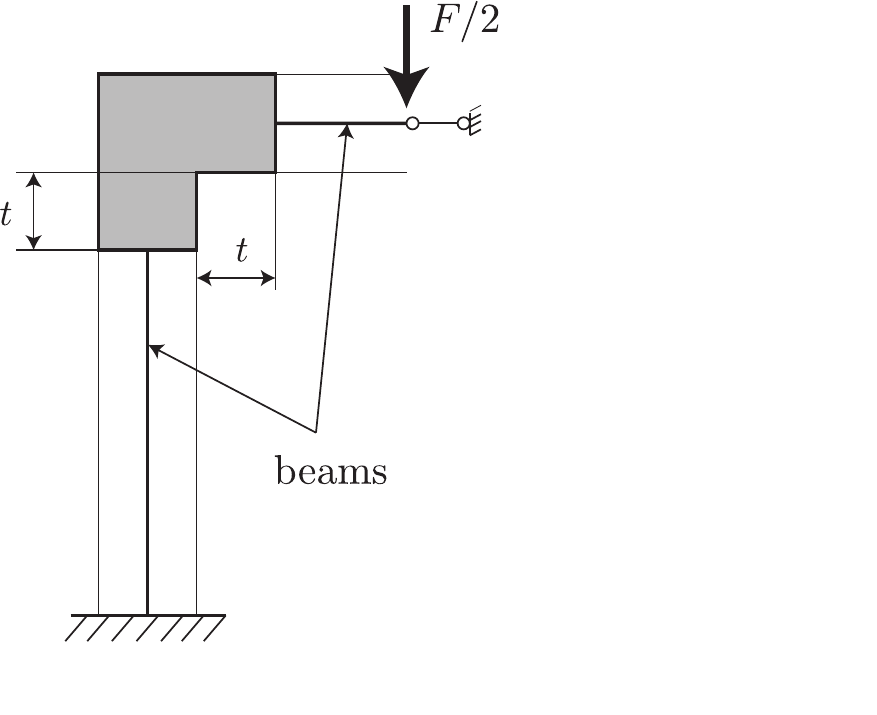}
  \caption{A plane frame analysis: solid-beam model.}
  \label{fig:frame-geo}
\end{figure}

\begin{figure}[htbp]
  \centering 
   \includegraphics[width=0.49\textwidth]{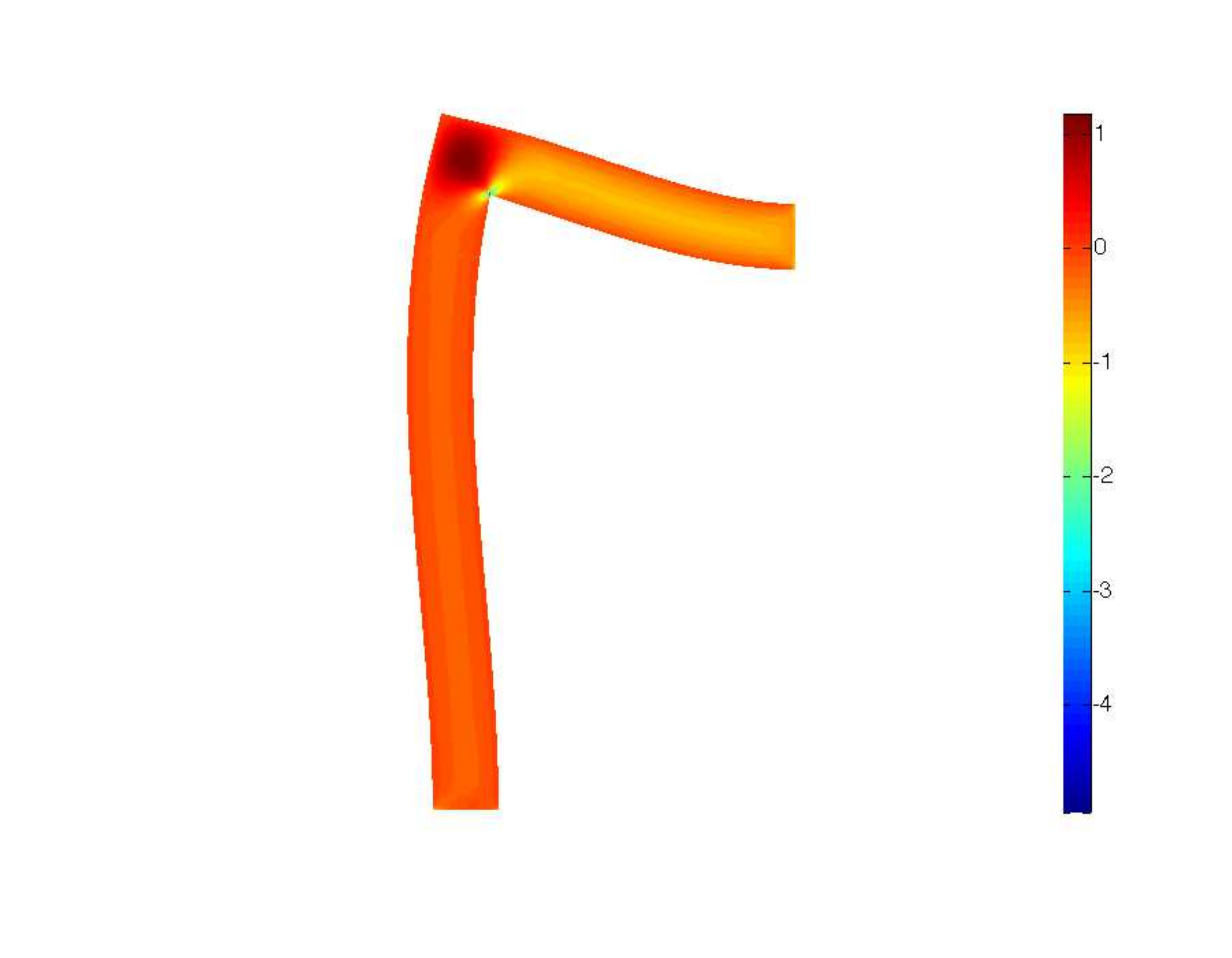}
   \includegraphics[width=0.49\textwidth]{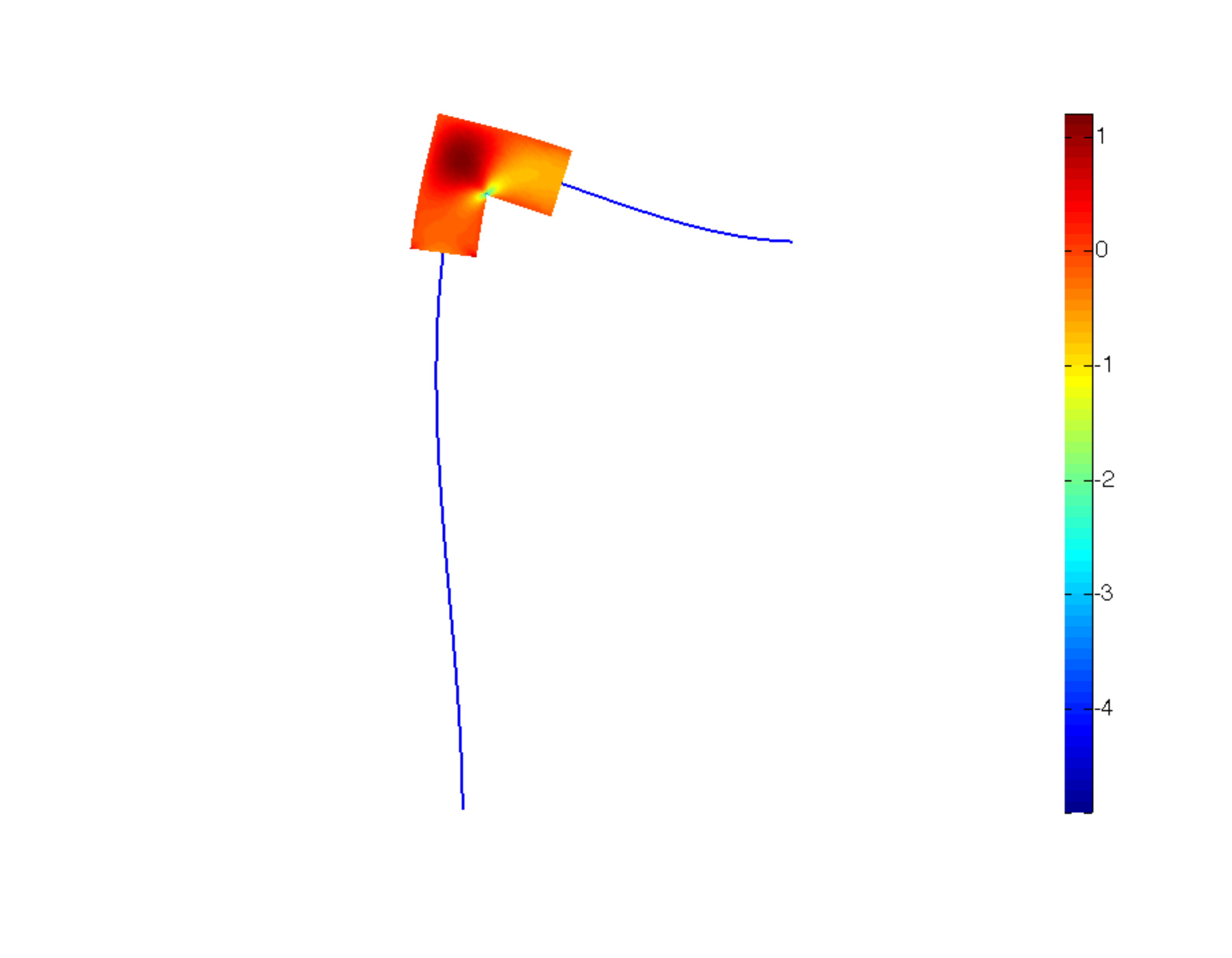}
  \caption{A plane frame analysis: comparison of $\sigma_{xy}$ contour plot obtained with
     solid model (left) and solid-beam model (right).}
  \label{fig:frame-res}
\end{figure}

\begin{rmk}
Although the processing time of the solid-beam model is much less than the one of the solid model, one cannot simply conclude that the solid-beam model is more efficient. The pre-processing of the solid-beam model, if not automatic, can be time consuming such that the gain in the processing step is lost. For non-linear analyses, where the processing time is dominant, we believe that mixed dimensional analysis is very economics.
\end{rmk}

\subsection{Continuum-plate coupling}

\subsubsection{Cantilever plate: conforming coupling}

For verification of the continuum-plate coupling, we consider the 3D cantilever beam given in 
Fig.~\ref{fig:plate-continuum-geo}. The material properties are $E=1000$ N/mm$^2$, $\nu=0.3$.
The end shear traction is $\bar{t}=10$ N/mm in case of continuum-plate model and is
$\bar{t}=10/20$ N/mm$^2$ in case of continuum model which is referred to as the reference model.
We use B-splines elements to solve both the MDA and the reference model. The length of the continuum part
in the continuum-plate model is $L/2=160$ mm.
A mesh of $64\times4\times5$ tri-cubic elements is utilized for the reference model
and a mesh of $32\times4\times5$/ $16\times2$ cubic elements is utilized for the mixed dimensional model,
cf. Fig.~\ref{fig:plate-continuum-meshes}. The plate part of the mixed dimensional model is discretised using
the Reissner-Mindlin plate theory with three unknowns per node and the Kirchhoff plate theory with only one
unknown per node. 
The stabilisation parameter was chosen empirically to be $5\times10^3$. Note that the eigenvalue method
described in Section \ref{sec:numerical-analysis} can be used to rigorously determine $\alpha$. 
However since it would be expensive for large problems, we are in favor of simpler but less rigorous rules
to compute this parameter.
Fig.~\ref{fig:plate-continuum-deform} shows a comparison of deformed shapes of the
continuum model and the continuum-plate model and in Fig.~\ref{fig:plate-continuum-stresses}, the contour plot
of the von Mises stress corresponding to various models is given.

\begin{figure}[htbp]
  \centering 
   \includegraphics[width=0.55\textwidth]{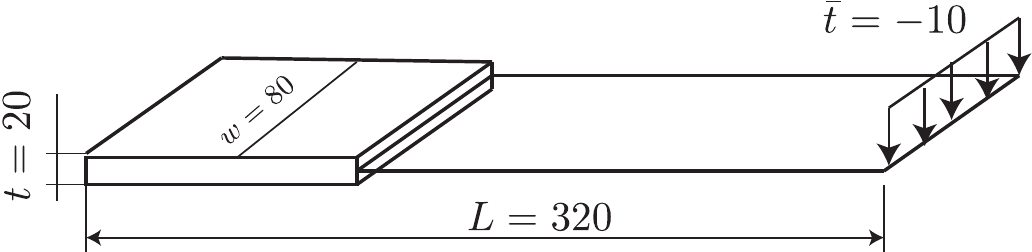}
  \caption{Cantilever beam subjects to an end shear force: problem setup.}
  \label{fig:plate-continuum-geo}
\end{figure}

\begin{figure}[htbp]
  \centering 
   \includegraphics[width=0.55\textwidth]{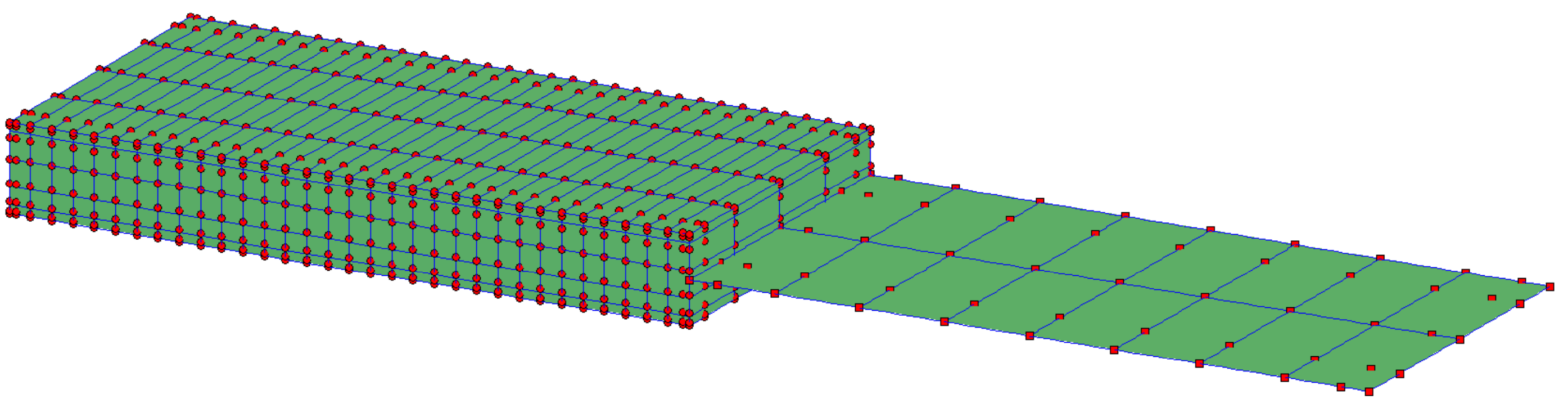}
  \caption{Cantilever beam subjects to an end shear force: typical B-spline discretisation.}
  \label{fig:plate-continuum-meshes}
\end{figure}

\begin{figure}[htbp]
  \centering 
   \includegraphics[width=0.55\textwidth]{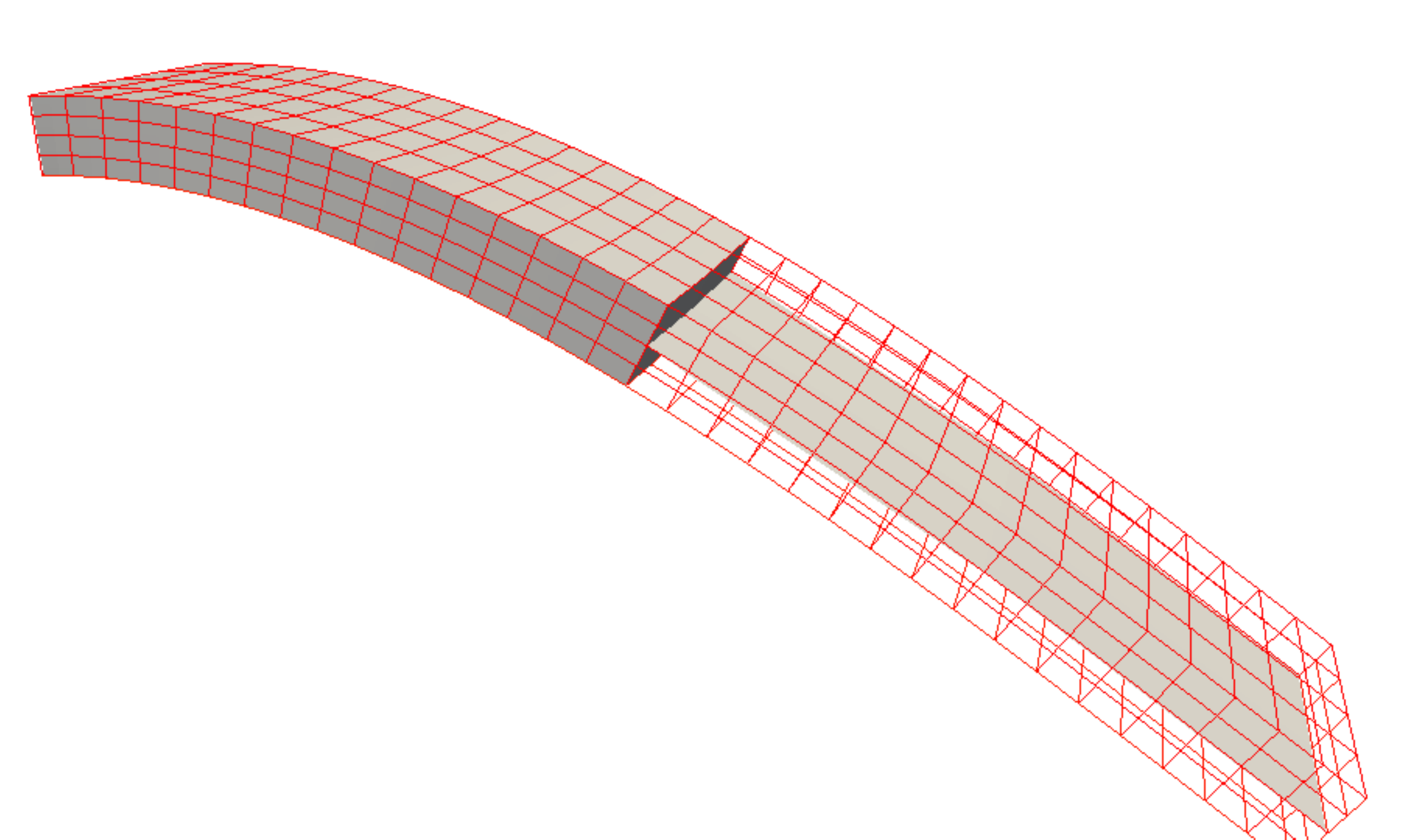}
  \caption{Cantilever beam subjects to an end shear force: comparison of deformed shapes of the
     continuum model and the continuum-plate model.}
  \label{fig:plate-continuum-deform}
\end{figure}

\begin{figure}[htbp]
  \centering 
  \subfloat[reference model]{\includegraphics[width=0.5\textwidth]{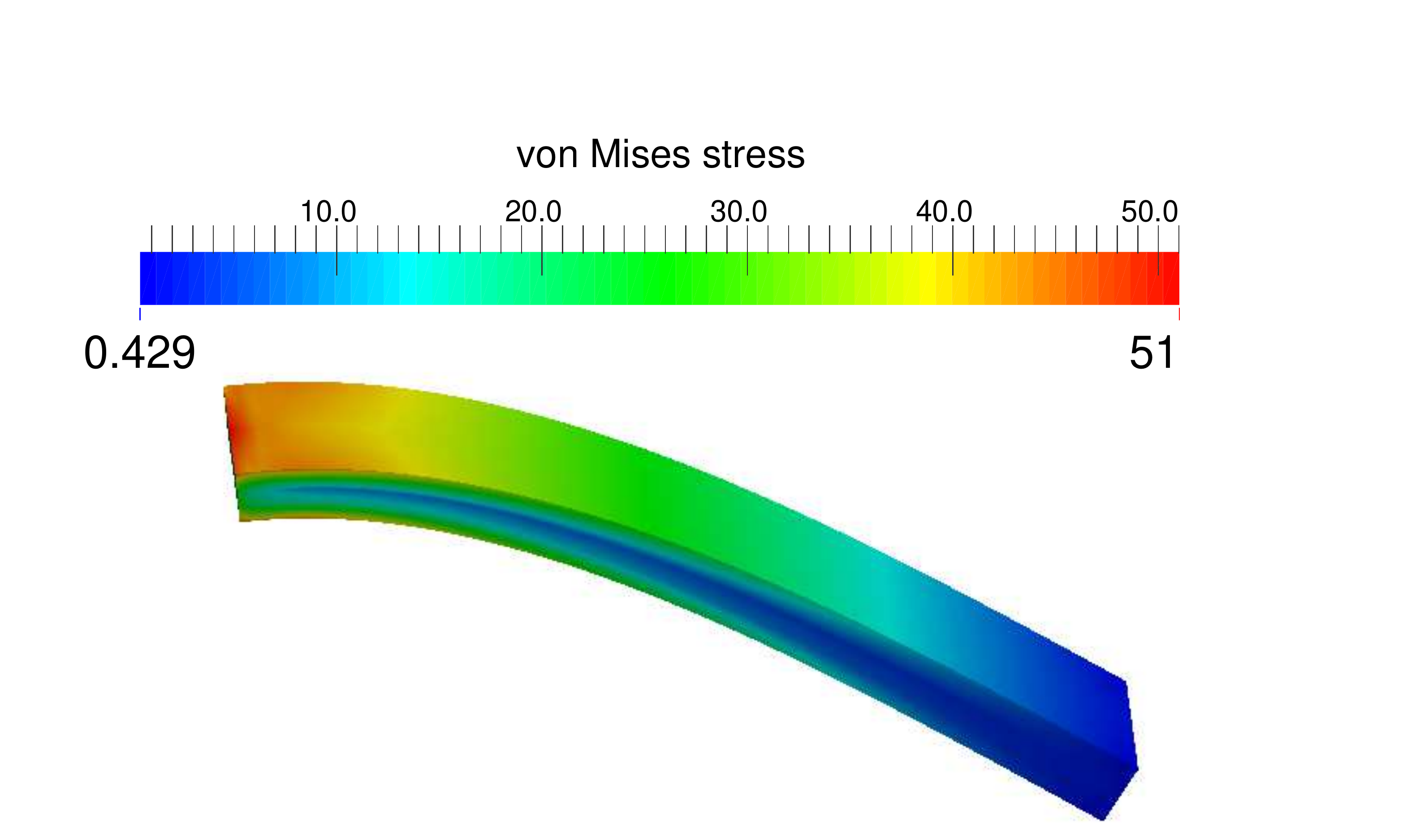}}\\
  \subfloat[mixed dimensional model, Mindlin plate]{\includegraphics[width=0.5\textwidth]{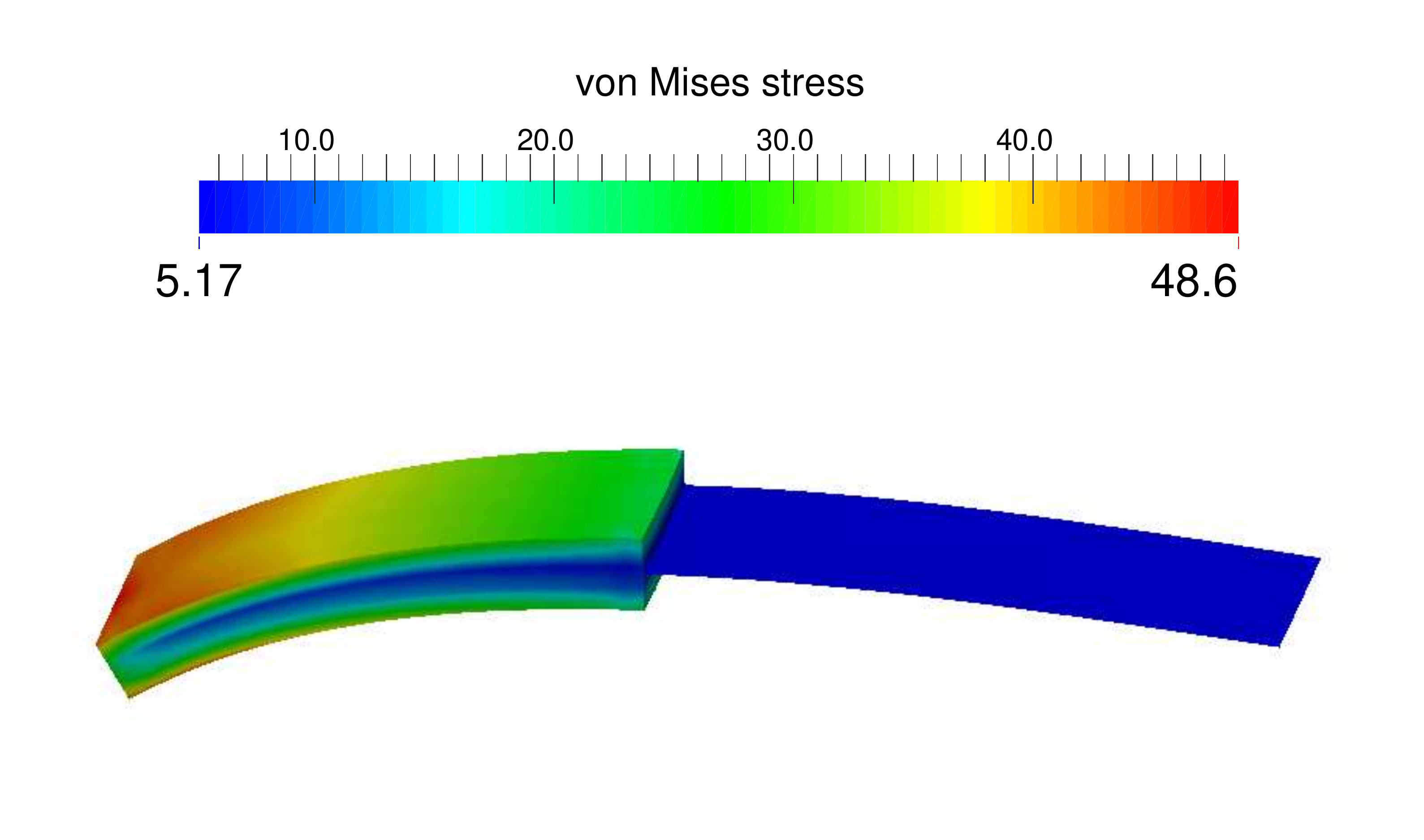}}
  \subfloat[mixed dimensional model, Kirchhoff plate]{\includegraphics[width=0.5\textwidth]{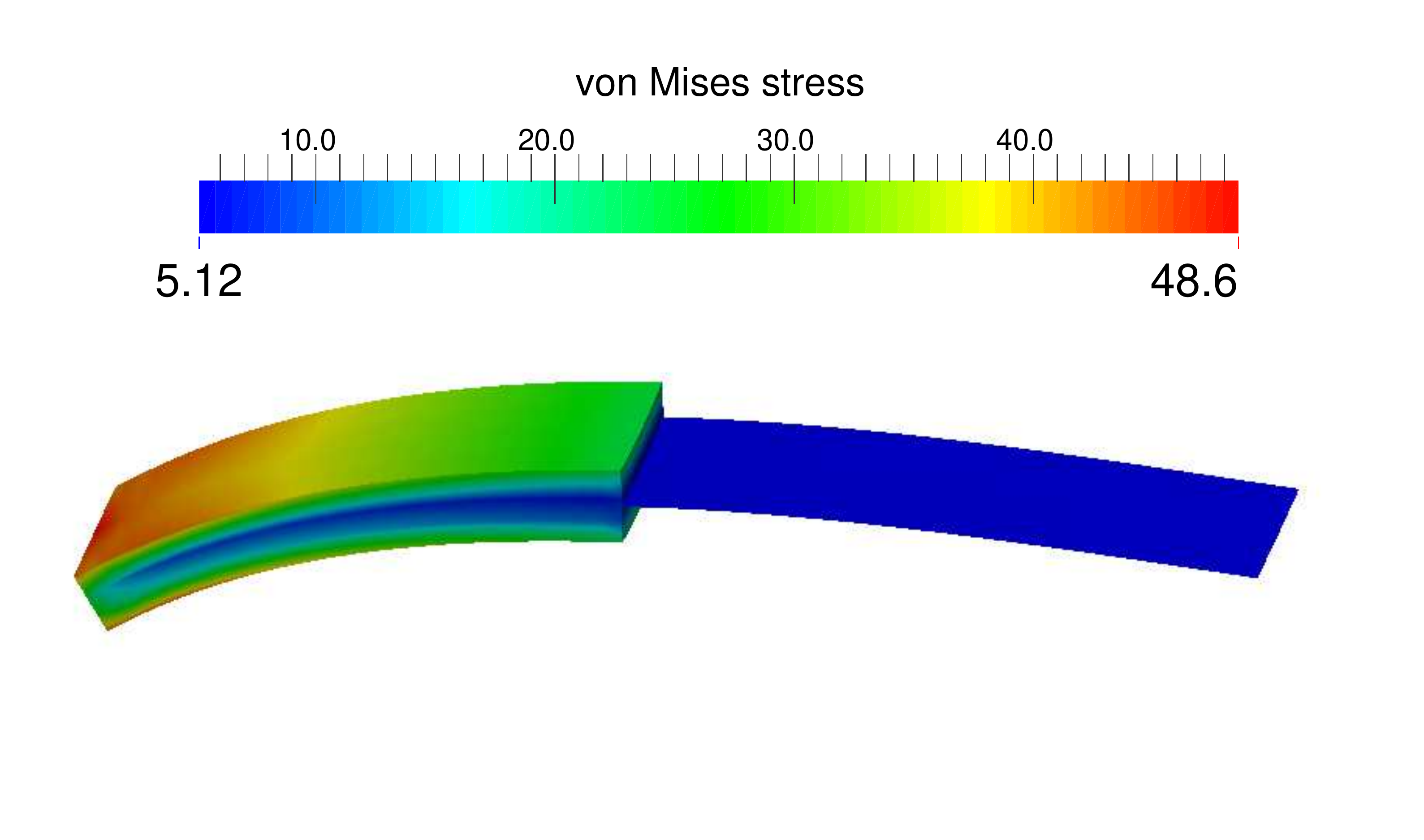}}
  \caption{Cantilever beam subjects to an end shear force}
  \label{fig:plate-continuum-stresses}
\end{figure}

\subsubsection{Cantilever plate: non-conforming coupling}

A mesh of $32\times4\times5$/ $32\times2$ cubic elements is utilized for the mixed dimensional model,
cf. Fig.~\ref{fig:plate-continuum-nonconform-mesh}. The length of the continuum part
in the continuum-plate model is $175$ mm. The contour plot of the von Mises stress is given 
Fig.~\ref{fig:plate-continuum-nonconform-stress} where void plate elements were removed in the visualisation.

\begin{figure}[htbp]
  \centering 
   \includegraphics[width=0.55\textwidth]{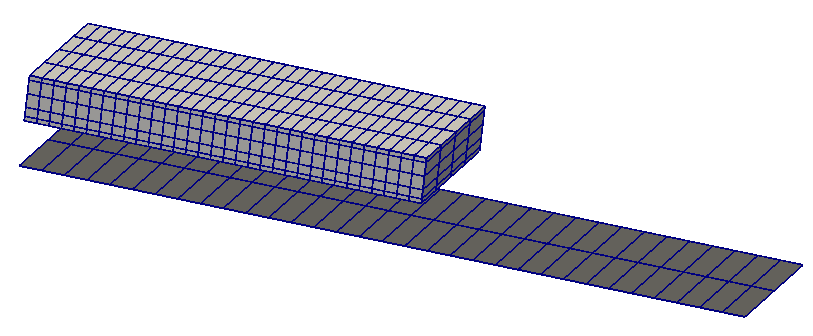}
  \caption{Cantilever beam subjects to an end shear force: discretisation of the solid and the plate.}
  \label{fig:plate-continuum-nonconform-mesh}
\end{figure}

\begin{figure}[htbp]
  \centering 
   \includegraphics[width=0.55\textwidth]{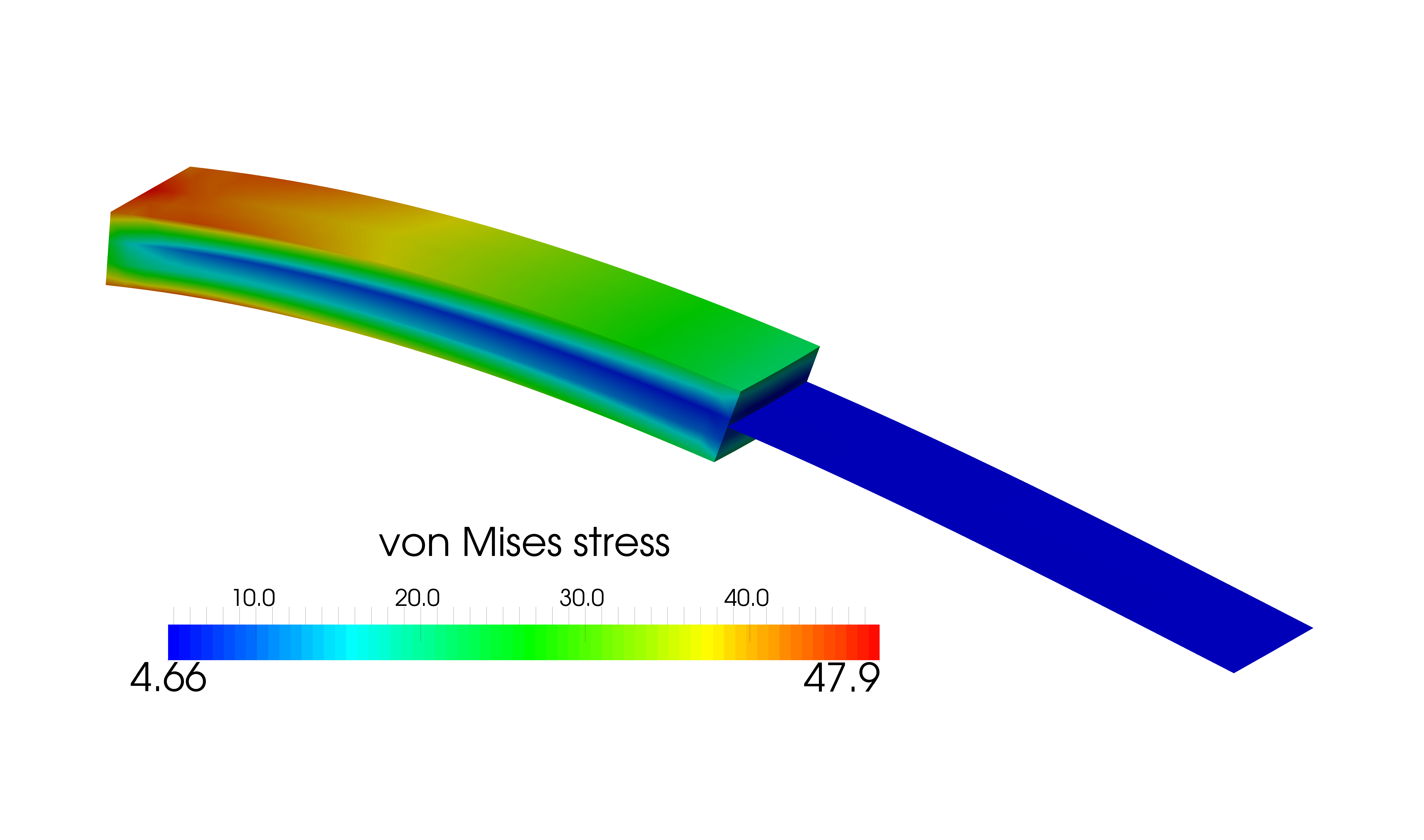}
  \caption{Cantilever beam subjects to an end shear force: von Mises stress distribution.}
  \label{fig:plate-continuum-nonconform-stress}
\end{figure}

\subsubsection{Non-conforming coupling of a square plate}
  
We consider a square plate of dimension $L\times L \times t$ ($t$ denotes the thickness) 
in which there is an overlapped solid of dimension $L_s\times L_s\times t$ as shown 
in Fig.~\ref{fig:squareplate-geo}.  In the computations, material properties are taken as $E=
10^3$, $\nu = 0.3$ and the geometry data are $L=400$, $t=20$ and $L_s=100$. The loading is a gravity
force $p=10$ and the plate boundary is fully clamped. 
The stabilisation parameter was chosen empirically to be $1\times10^6$.
We use rotation free Kirchhoff NURBS plate elements for the plate and NURBS solid
elements for the solid. In order to model zero rotations in a rotation free NURBS plate formulation,
we simply fix the transverse displacement of control points on the boundary and those right next to them
cf. \cite{kiendl_isogeometric_2009}.         

\begin{figure}[htbp]
  \centering 
   \includegraphics[width=0.55\textwidth]{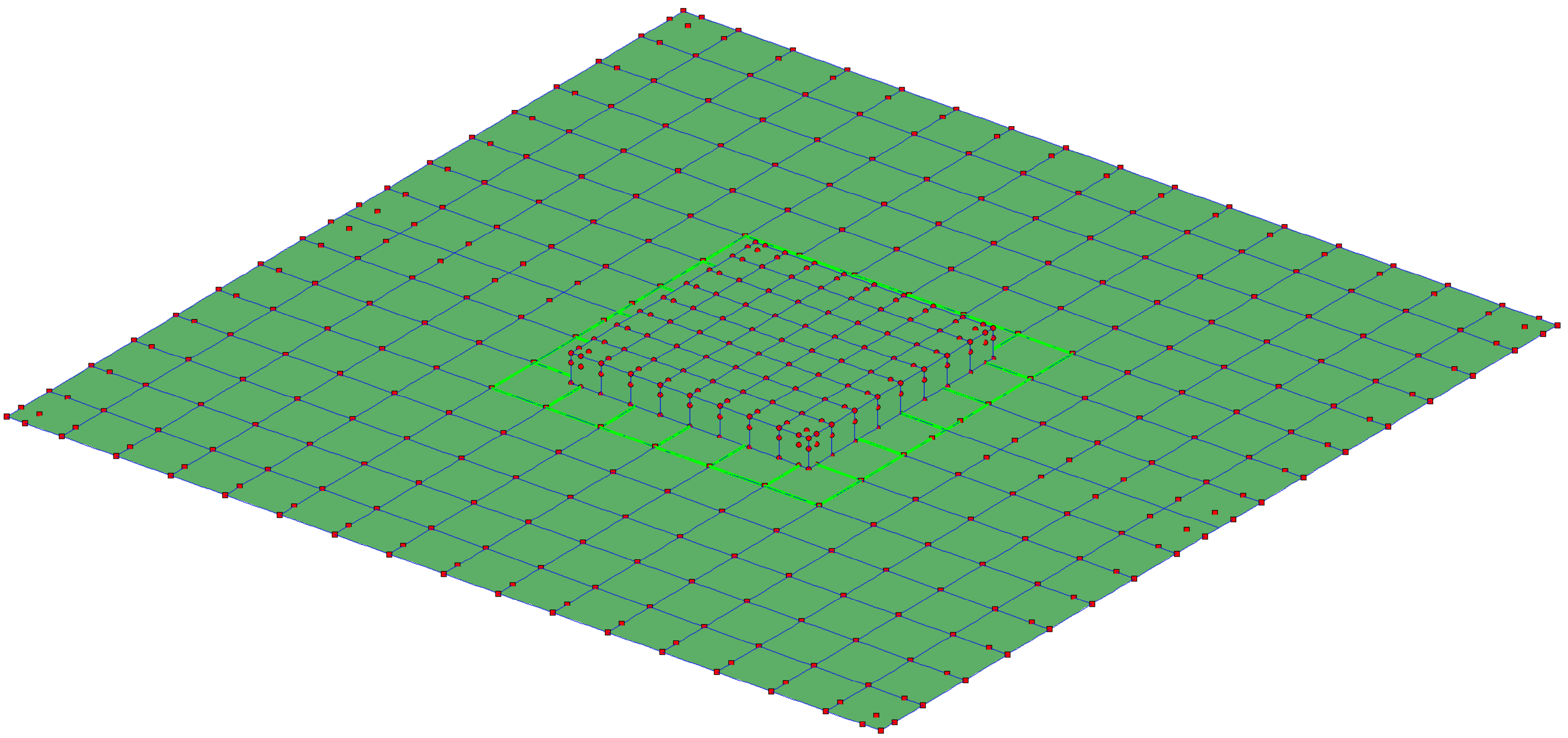}
  \caption{Square plate enriched by a solid. The highlighted elements are those plate elements cut by the 
     solid boundaries. The plate is fully clamped ans subjected to a gravity force.}
  \label{fig:squareplate-geo}
\end{figure}

In order to find plate elements cut by the boundary surfaces of the solid, we use the level sets
defined for the square which is the intersection plane of the solid and the plate. The use of level sets
to define the interaction of finite elements with some geometry entities is popular in XFEM, see \eg  \cite{Sukumar1}. Fig.~\ref{fig:squareplate-deform1} plots the deformed configuration of the solid-plate model and the one obtained with a plate model. A good agreement can be observed. In order to show the flexibility of the non-conforming coupling, the solid part was
moved slightly to the right and the deformed configuration is given in Fig.~\ref{fig:squareplate-deform2}. The same discretisation for the plate is used. This should serve as a prototype for model adaptivity analyses to be presented in a forthcoming contribution.

\begin{figure}[htbp]
  \centering 
   \includegraphics[width=0.49\textwidth]{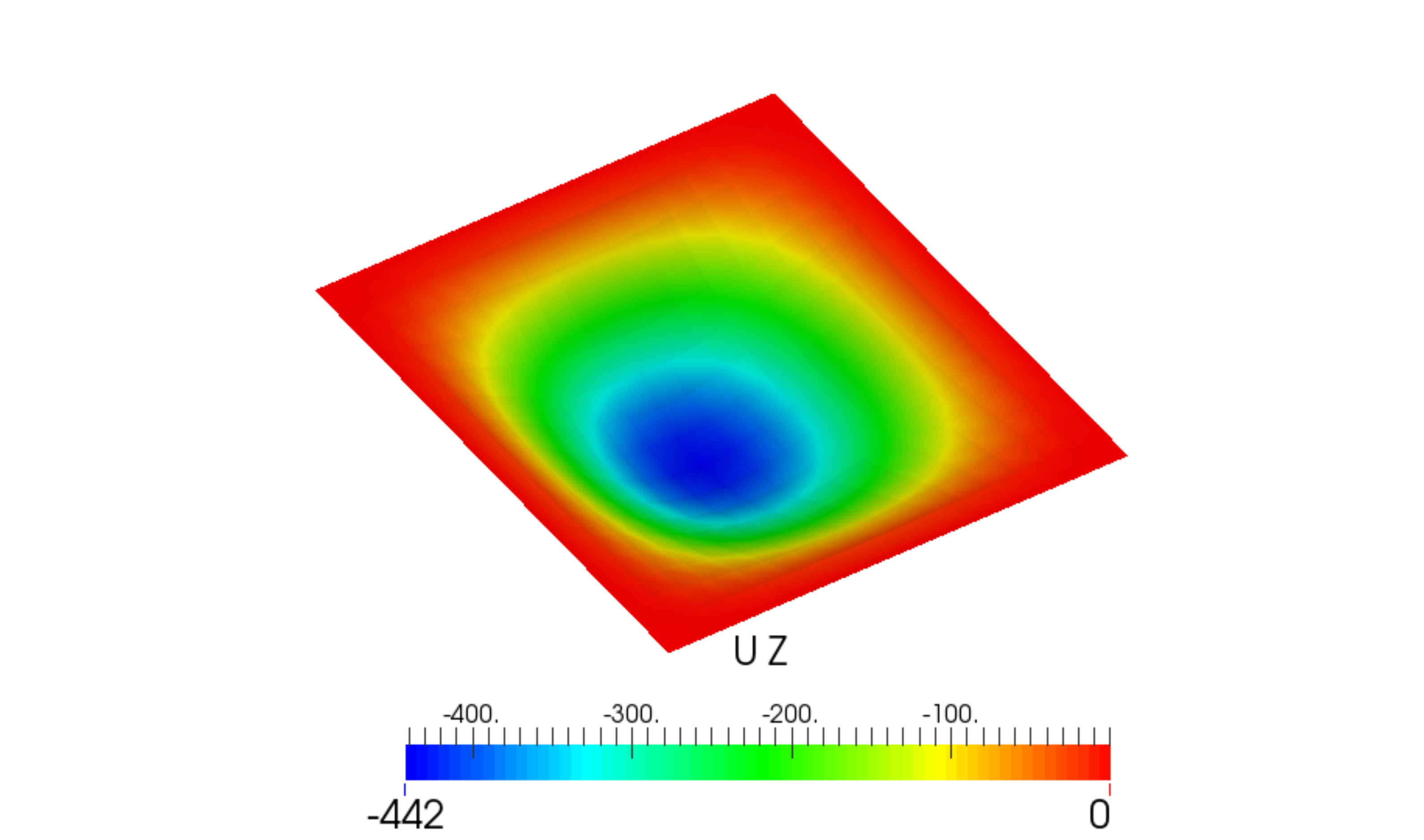}
   \includegraphics[width=0.49\textwidth]{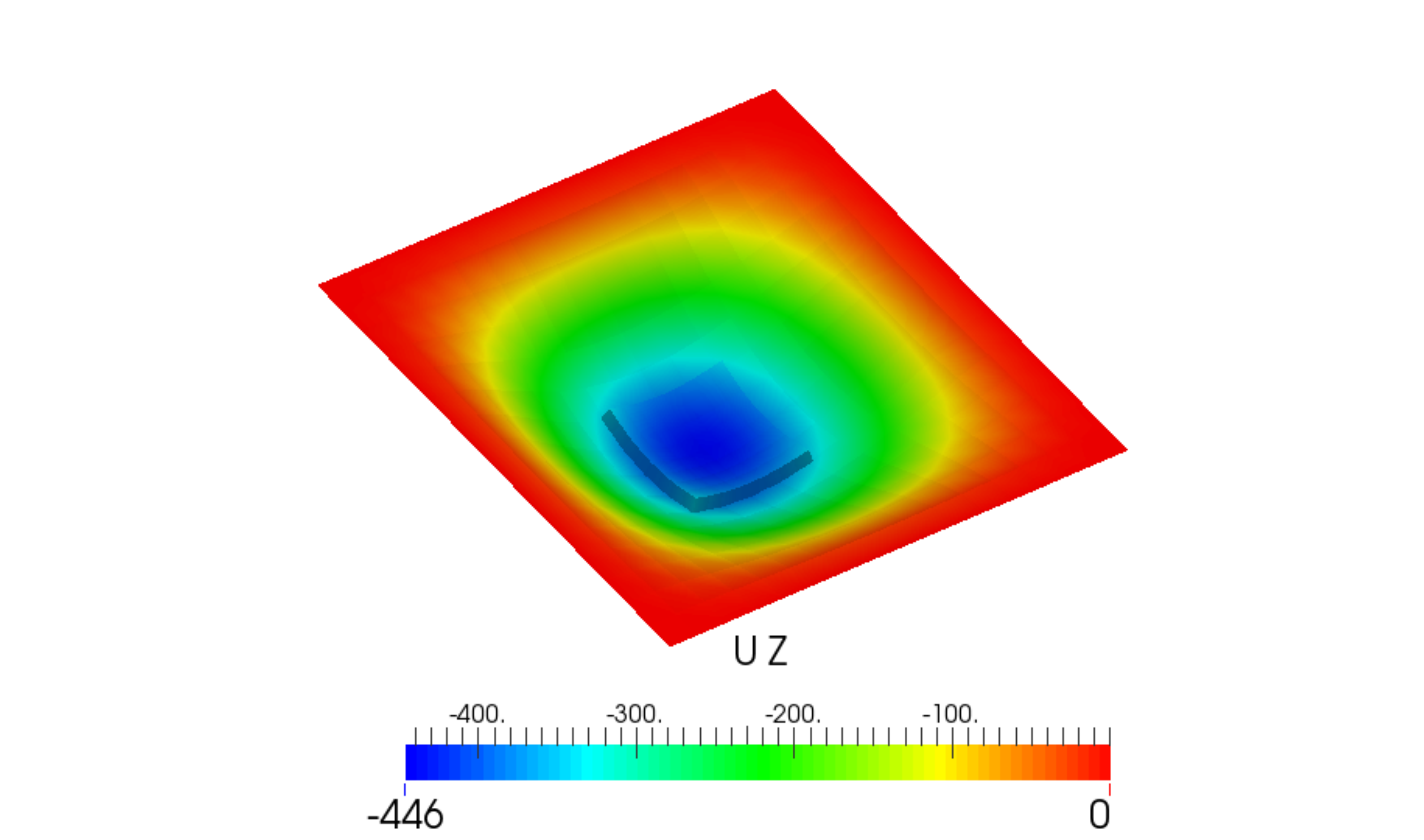}
  \caption{Square plate enriched by a solid: transverse displacement plot on deformed configurations of
     plate model (left) and solid-plate model (right).}
  \label{fig:squareplate-deform1}
\end{figure}

\begin{figure}[htbp]
  \centering 
   \includegraphics[width=0.49\textwidth]{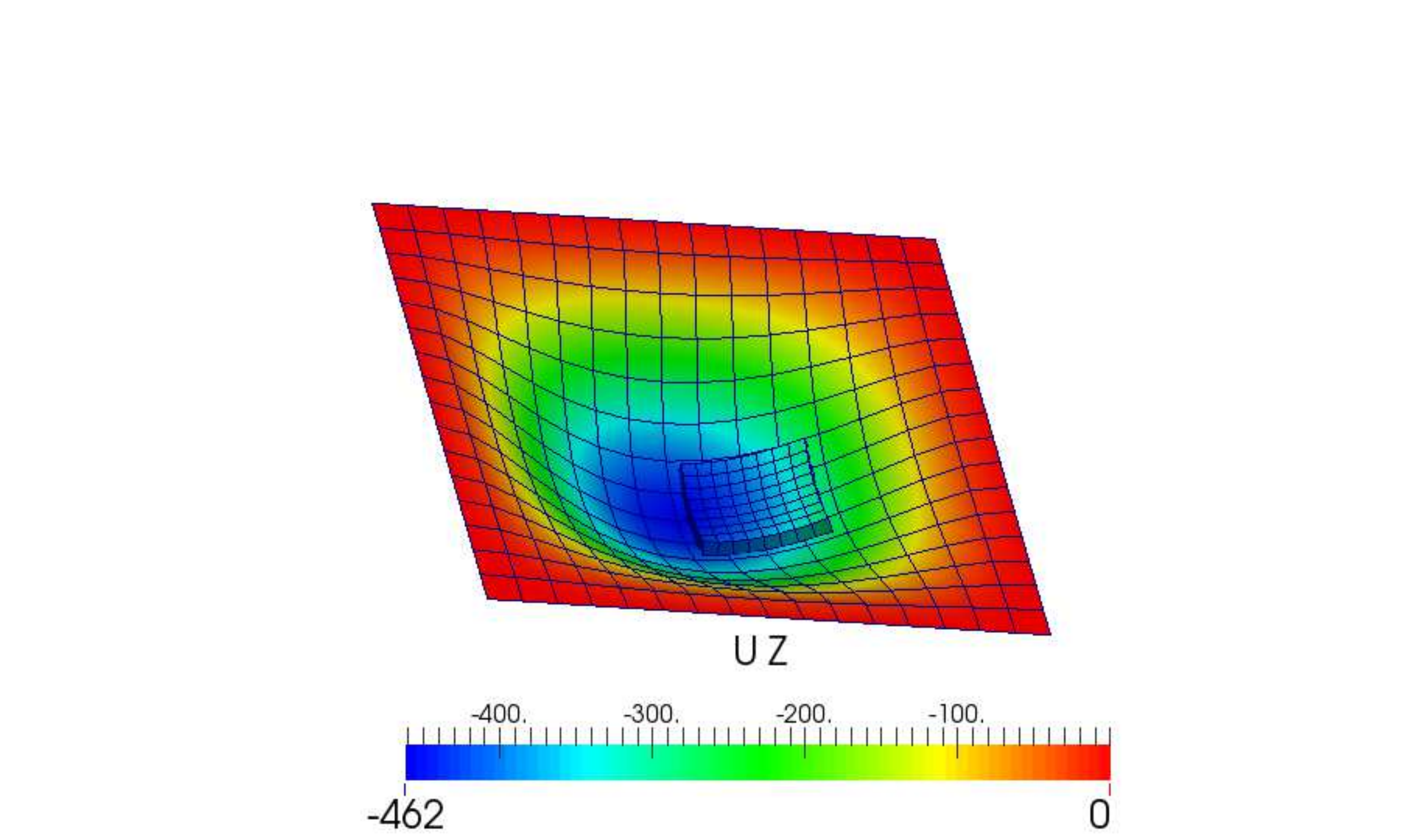}
  \caption{Square plate enriched by a solid: transverse displacement plot where the solid part was
moved slightly to the right .}
  \label{fig:squareplate-deform2}
\end{figure}

\section{Conclusions}\label{sec:conclusions}

We presented a Nitsche's method to couple (1) two dimensional continua and beams 
and (2) three dimensional continua and plates. A detailed implementation of those coupling methods was given. 
Numerical examples using low order Lagrange finite elements and high order
B-spline/NURBS isogeometric finite elements provided demonstrate 
the good performance of the method and its versatility. Both classical beam/plate theories
and first order shear beam/plate models were presented. Conforming coupling where the continuum
mesh and the beam/plate mesh is not overlapped and non-conforming coupling where they are overlapped
are described. The latter provides great flexibility in model adaptivity formulation and the implementation
is much more simpler than the Arlequin method. We also presented a numerical analysis of the bilinear form
that results in technique to compute the minimum value for the stabilisation parameter ensuring the positive
definiteness of the stiffness matrix. The method is, however, expensive. 

The contribution was limited to linear static problems and the extension of the method to (1) 
more complex and detailed analysis of non-linear dynamics problems and (2) nonlinear material problems 
is under way. This will allow to verify the potential of Nitsche coupling for mixed dimensional analysis 
in realistic engineering applications. The nonconforming coupling when combined with an error estimator
will provide an efficient methodology to analyse engineering structures.

\section*{Acknowledgements}

The authors would like to acknowledge the partial financial support of the
Framework Programme 7 Initial Training Network Funding under grant number
289361 ``Integrating Numerical Simulation and Geometric Design Technology".
St\'{e}phane Bordas also thanks partial funding for his time provided by
1) the EPSRC under grant EP/G042705/1 Increased Reliability for Industrially
Relevant Automatic Crack Growth Simulation with the eXtended Finite Element
Method and 2) the European Research Council Starting Independent Research
Grant (ERC Stg grant agreement No. 279578) entitled ``Towards real time multiscale simulation of cutting in
non-linear materials with applications to surgical simulation and computer
guided surgery''.
The authors would like to express the gratitude towards Drs. Erik Jan Lingen
and Martijn Stroeven at the Dynaflow Research Group, Houtsingel 95, 2719 EB Zoetermeer, The Netherlands 
for providing us the numerical toolkit jem/jive.

\bibliography{isogeometric}
\bibliographystyle{unsrt}


\end{document}